\documentclass[11pt]{article}

\usepackage{amsmath,amssymb,amsfonts}
\usepackage{pstricks,pst-node,xy}
\usepackage{xr}

\xyoption{all}

\numberwithin{equation}{section}
\newtheorem{thm}{Theorem}[section] 
\newtheorem{prp}[thm]{Proposition}
\newtheorem{lmm}[thm]{Lemma}  
\newtheorem{crl}[thm]{Corollary}
\newtheorem{dfn}[thm]{Definition} 

\def\e_ref#1{(\ref{#1})}
\def\blr#1{\big\langle{#1}\big\rangle}
\def\Blr#1{\Big\langle{#1}\Big\rangle}
\def\ov#1{\overline{#1}}

\def\ti#1{\tilde{#1}}
\def\wt#1{\widetilde{#1}}
\def\smsize#1{\begin{small}#1\end{small}}
\def\Lgsize#1{\begin{Large}#1\end{Large}}

\renewcommand{\cal}{\mathcal}
\renewcommand{\frak}{\mathfrak}
\renewcommand{\Bbb}{\mathbb}

\def\A{\cal A}
\def\bA{\Bbb A}
\def\B{\cal B}
\def\C{\Bbb C}
\def\cC{\cal C}
\def\D{\frak D}
\def\cD{\cal D}
\def\d{\frak d}
\def\E{\Bbb E}
\def\cE{\cal E}

\def\L{\Bbb L}
\def\cL{\cal L}
\def\M{\frak M}
\def\cM{\cal M}
\def\N{\cal N}
\def\O{\cal O}
\def\P{\Bbb P}
\def\Pf{\P^4}
\def\Pn{\P^n}
\def\Q{\Bbb Q}
\def\S{\cal S}
\def\U{\frak U}
\def\V{\cal V}
\def\W{\cal W}
\def\Z{\Bbb Z}
\def\cZ{\cal Z}

\def\ale{\aleph}
\def\al{\alpha}
\def\be{\beta}
\def\ga{\gamma}
\def\io{\iota}
\def\na{\nabla}
\def\si{\sigma}
\def\ups{\upsilon}
\def\vph{\varphi}
\def\vr{\varrho}

\def\Ga{\Gamma}
\def\Si{\Sigma}

\def\lra{\longrightarrow}
\def\Lra{\Longrightarrow}

\def\eset{\emptyset}
\def\part{\partial}

\def\Aut{\textnormal{Aut}}
\def\st{\textnormal{s.t.}}
\def\Bl{\textnormal{Bl}}
\def\Pr{\textnormal{Pr}}
\def\codim{\textnormal{codim}}
\def\Def{\textnormal{Def}}
\def\Edg{\textnormal{Edg}}
\def\eff{\textnormal{eff}}
\def\ev{\textnormal{ev}}
\def\GW{\textnormal{GW}}
\def\Hom{\textnormal{Hom}}
\def\id{\textnormal{id}}
\def\ide{\textnormal{ide}}
\def\Im{\textnormal{Im}}
\def\Pr{\textnormal{Pr}}
\def\rk{\textnormal{rk}}
\def\val{\textnormal{val}}
\def\Ver{\textnormal{Ver}}

\begin{document}

\title{A Desingularization of the Main Component of\linebreak
the Moduli Space of Genus-One Stable Maps into $\Pn$}

\author{Ravi Vakil\thanks{Partially supported by an NSF grant DMS--0228011}~~and 
Aleksey Zinger\thanks{Partially supported by an NSF Postdoctoral Fellowship}}
\date{\today}

\maketitle

\begin{abstract}
\noindent
We construct a natural smooth compactification of the space of smooth genus-one
curves with $k$ distinct points in a projective space. 
It can be viewed as an analogue of a well-known smooth compactification 
of the space of smooth genus-zero curves, i.e.~the 
space of stable genus-zero maps $\ov\M_{0,k}(\Pn,d)$.
In fact, our compactification is obtained from the singular  
space of stable genus-one maps $\ov\M_{1,k}(\Pn,d)$ through 
a natural sequence of blowups along ``bad'' subvarieties. 
While this construction is simple to describe, 
it requires more work to show that the end result is a smooth space. 
As a bonus, we obtain desingularizations of certain natural sheaves over 
the ``main'' irreducible component $\ov\M_{1,k}^0(\Pn,d)$ of $\ov\M_{1,k}(\Pn,d)$.
A number of applications of these desingularizations in enumerative geometry
and Gromov-Witten theory are described in the introduction.  
\end{abstract}

\tableofcontents

\section{Introduction}
\label{intro_sec}

\subsection{Background and Applications}
\label{background_subs}

\noindent
The space of degree-$d$ genus-$g$ curves with $k$ distinct marked points in $\Pn$ 
is generally not compact, but admits a number of natural 
compactifications\footnote{We call a space $\ov\M$ a {\it compactification} of $\M$ 
if $\ov\M$ is compact and contains~$\M$.
In particular, $\M$ need not be dense in~$\ov\M$.}.
Among the most prominent compactifications is the moduli space of stable genus-$g$
maps, $\ov\M_{g,k}(\Pn,d)$, constructed in~\cite{Gr} and~\cite{FuP}.
It has found numerous applications in classical enumerative geometry and 
is a central object in Gromov-Witten theory.
However, most applications in enumerative geometry and some results
in GW-theory have been restricted to the genus-zero case.
The reason for this is essentially that the genus-zero moduli space 
has a particularly simple structure: it is smooth and contains
the space of smooth genus-zero curves as a dense open subset.
On the other hand, the moduli spaces of positive-genus stable maps 
fail to satisfy either of these two properties.
In fact, $\ov\M_{g,k}(\Pn,d)$ can be arbitrarily singular according to~\cite{murphy}.
It is thus natural to ask whether these failings can be remedied by modifying   
$\ov\M_{g,k}(\Pn,d)$, preferably in a way that leads to a range of applications.
As announced in~\cite{summ} and shown in this paper, the answer is yes if $g\!=\!1$.\\

\noindent
We denote by $\M_{1,k}(\Pn,d)$ the subset of $\ov\M_{1,k}(\Pn,d)$
consisting of the stable maps that have smooth domains.
This space is smooth and contains the space of genus-one curves with $k$ 
distinct marked points in $\Pn$ as a dense open subset, provided $d\!\ge\!3$.
However, $\M_{1,k}(\Pn,d)$ is not compact.
Let $\ov\M_{1,k}^0(\Pn,d)$ be the closure of $\M_{1,k}(\Pn,d)$ in 
the compact space $\ov\M_{1,k}(\Pn,d)$.
While $\ov\M_{1,k}^0(\Pn,d)$ is not smooth, it turns out that 
a natural sequence of blowups along loci disjoint from $\M_{1,k}(\Pn,d)$
leads to a desingularization of $\ov\M_{1,k}^0(\Pn,d)$,
which will be denoted by $\wt\M_{1,k}^0(\Pn,d)$.\\

\noindent
The situation is as good as one could possibly hope.  
A general strategy when attempting to desingularize some space 
is to blow up the ``most degenerate'' locus,
then the proper transform of the ``next most degenerate locus'', and
so on.  This strategy works here, but with a novel twist: 
we apply it to the entire space of stable maps $\ov\M_{1,k}(\Pn,d)$.  
The most degenerate locus is in fact an entire irreducible component, 
and blowing it up removes it\footnote{Blowing up an irreducible component 
of a stack will result in the component being removed (or ``blown out of existence''), 
and the remainder of the stack
is blown up along its intersection with the component in question.}.
Hence one by one we erase the ``bad'' components of $\ov\M_{1,k}(\Pn,d)$.  
Each blowup of course changes the ``good'' component $\ov\M_{1,k}^0(\Pn,d)$, 
and miraculously at the end of the process the resulting space 
$\wt\M_{1,k}^0(\Pn,d)$ is nonsingular.
We note that this cannot possibly be true for an arbitrary $g$,
as $\M_{g,k}(\Pn,d)$ behaves quite badly according to~\cite{murphy}.
The sequential blowup construction itself is beautifully simple.
It is completely described in the part of Subsection~\ref{descr_subs}
ending with the main theorem of the paper, Theorem~\ref{main_thm}.
However, showing that $\wt\M_{1,k}^0(\Pn,d)$ is in fact smooth
requires a considerable amount of preparation (which takes up 
Subsections~\ref{curveprelim_subs}-\ref{map1prelim_subs})
and is finally completed in Subsection~\ref{map1blconstr_subs}.\\

\noindent
The desingularization $\wt\M_{1,k}^0(\Pn,d)$ of $\ov\M_{1,k}^0(\Pn,d)$
possesses a number of ``good'' properties and has a variety of applications  to 
enumerative algebraic geometry and Gromov-Witten theory.
It has already been observed in~\cite{Fo} that the cohomology of 
$\wt\M_{1,k}^0(\Pn,d)$ behaves in a certain respect like the cohomology of 
the moduli space of genus-one curves, $\ov\cM_{1,k}$.
The space $\wt\M_{1,k}^0(\Pn,d)$ can be used to count genus-one curves in~$\Pn$,
mimicking the genus-zero results of~\cite{KM} and~\cite{RT}
(though perhaps not their simple recursive formulas).
Proceeding analogously to the genus-zero case (e.g.~as in~\cite{P}, \cite{enumtang},
and~\cite{genuss0pr}), Theorem~\ref{main_thm} can then be used to count
genus-one curves with tangency conditions and singularities.
In all cases, such counts can be expressed as integrals of
natural cohomology classes on $\ov\M_{1,k}^0(\Pn,d)$ or $\wt\M_{1,k}^0(\Pn,d)$.
Integrals on the latter space can be computed using the localization
theorem of~\cite{ABo}, as $\wt\M_{1,k}^0(\Pn,d)$ is smooth and inherits a torus action 
from $\Pn$ and $\ov\M_{1,k}(\Pn,d)$.\\

\noindent
We next discuss two types of applications of Theorem~\ref{main_thm} in 
Gromov-Witten theory, as well as a bonus result of this paper, Theorem~\ref{cone_thm}. 
It is shown in~\cite{g1comp} and~\cite{g1comp2} that the space $\ov\M_{1,k}^0(\Pn,d)$
has a natural generalization to arbitrary almost Kahler manifolds and gives
rise to new symplectic {\it reduced genus-one GW}-invariants.
These reduced invariants are yet to be constructed in algebraic geometry.
However, the spaces $\wt\M_{1,k}^0(\Pn,d)$ do possess a number of 
``good'' properties and give rise to algebraic invariants of algebraic manifolds;
see the first and last sections of~\cite{summ}.
It is not clear whether these are the same as the reduced genus-one invariants, 
but it may be possible to verify this by using Theorem~\ref{cone_thm}.\\

\noindent
Theorem~\ref{main_thm} also has applications to computing Gromov-Witten
invariants of complete intersections, once it is combined with Theorem~\ref{cone_thm}.
Let $a$ be a nonnegative integer. 
For a general $s\!\in\!H^0(\Pn,\O_{\Pn}(a))$, 
$$Y\equiv s^{-1}(0) \subset \Pn$$
is a smooth hypersurface. 
We denote its degree-$d$ Gromov-Witten invariant by $\GW_{g,k}^Y(d;\cdot)$,
i.e. 
$$\GW_{g,k}^Y(d;\psi)\equiv
\blr{\psi,\big[\ov\M_{g,k}(Y,d)\big]^{vir}}
\qquad \text{for all} \qquad \,  \psi\!\in\!H^*\big(\ov\M_{g,k}(Y,d);\Q\big).$$
Suppose $\U$ is the universal curve over $\ov\M_{g,k}(\Pn,d)$,
with structure map~$\pi$ and evaluation map~$\ev$:
$$\xymatrix{\U \ar[d]^{\pi} \ar[r]^{\ev} & \Pn \\
\ov\M_{g,k}(\Pn,d).}$$
It can be shown that
\begin{equation}\label{genus0_e}
\GW_{0,k}^Y(d;\psi)
=\blr{\psi\cdot e\big(\pi_*\ev^* \O_{\Pn}(a)\big), \big[\ov\M_{0,k}(\Pn,d)\big]}
\end{equation}
for all $\psi\!\in\!H^*(\ov\M_{0,k}(\Pn,d);\Q)$; see~\cite{Bea} for example.
The moduli space $\ov\M_{0,k}(\Pn,d)$ is a smooth orbifold  and
$$\pi_*\ev^* \O_{\Pn}(a) \lra\ov\M_{0,k}(\Pn,d)$$
is a locally free sheaf, i.e.~a vector bundle.
The right-hand side of~\e_ref{genus0_e}
can be computed via the classical localization theorem of~\cite{ABo}.
The complexity of this computation increases quickly with the degree~$d$,
but it has been completed in full generality in a number of different of ways;
see \cite{Ber}, \cite{Ga}, \cite{Gi}, \cite{Le}, and~\cite{LLY}.\\

\noindent
If $n\!=\!4$, so $Y$ is a threefold, then
\begin{equation}\label{genus1_e}
\GW_{1,k}^Y(d;\psi)=
\frac{d(a\!-\!5)\!+\!2}{24}\GW_{0,k}^Y(d;\psi)
+\blr{\psi\cdot e(\pi_*\ev^*\O_{\Pn}(a)),\big[\ov\M_{1,k}^0(\Pf,d)\big]}
\end{equation}
for all $\psi\!\in\!H^*(\ov\M_{1,k}(\Pf,d);\Q)$; 
see \cite[Crl.~\ref{g1gw-cy_crl}]{LZ}.
This decomposition generalizes to arbitrary complete intersections~$Y$
and perhaps even to higher-genus invariants.
The~sheaf
\begin{equation}\label{sheaf_e}
\pi_*\ev^* \O_{\Pn}(a) \lra\ov\M_{1,k}^0(\Pf,d)
\end{equation}
is not locally free.
Nevertheless, its euler class is well-defined:
the euler class of every desingularization of this sheaf is the~same,
in the sense of \cite[Subsect.~\ref{g1cone-appr_subs}]{g1cone}.
This euler class can be geometrically interpreted as the zero set
of a sufficiently good section of the cone 
$$\V_{1,k}^d\lra\ov\M_{1,k}^0(\Pf,d),$$
naturally associated to the sheaf 
\e_ref{sheaf_e}\footnote{$\V_{1,k}^d$ is a variety such that 
the fibers of the projection map to $\ov\M_{1,k}^0(\Pf,d)$ are vector spaces,
but not necessarily of the same dimension.}; 
see the second part of the next subsection and Lemma~\ref{conesheaf_lmm}.\\

\noindent
One would hope to compute the last expression in~\e_ref{genus1_e} by localization.
However, since the variety $\ov\M_{1,k}^0(\Pf,d)$
and the cone $\V_{1,k}^d$ are singular, the localization theorem of~\cite{ABo} is not immediately applicable in the given situation.
Let 
$$\ti\pi\!: \wt\M_{1,k}^0(\Pf,d)\lra\ov\M_{1,k}^0(\Pf,d)$$
be the projection map.
As a straightforward extension of the main desingularization construction of
this paper, we show that  the cone 
$$\ti\pi^*\V_{1,k}^d\lra\wt\M_{1,k}^0(\Pf,d)$$
contains a vector bundle 
$$\wt\V_{1,k}^d\lra\wt\M_{1,k}^0(\Pf,d)$$
of rank $da=\rk\,\V_{1,k}^d|_{\M_{1,k}^0(\Pf,d)}$; see Theorem~\ref{cone_thm}.
It then follows that 
\begin{equation}\label{euler_e}\begin{split}
\blr{\psi\cdot e\big(\pi_*\ev^*\O_{\Pn}(a)\big),
\big[\ov\M_{1,k}^0(\Pf,d)\big]}
&\equiv \blr{\psi\cdot e\big(\V_{1,k}^d),\big[\ov\M_{1,k}^0(\Pf,d)\big]}\\
&=\blr{\ti\pi^*\psi\cdot e(\wt\V_{1,k}^d),\big[\wt\M_{1,k}^0(\Pf,d)\big]}.
\end{split}\end{equation}
The last expression above is computable by localization.\\

\noindent
{\it Remark:} Another approach to computing positive-genus Gromov-Witten invariants
has been proposed in~\cite{MaP}.
In contrast to the approach of~\cite{LZ},  
it applies to arbitrary-genus invariants, but can at present be used to compute
invariants of only low-dimensional and/or low-degree complete intersections.\\

\noindent
The main desingularization construction of this paper is the subject of 
Section~\ref{map1bl_sec}, but its key aspects are presented in the next subsection.
The construction itself and its connections with Sections~\ref{curvebl_sec} 
and~\ref{map0bl_sec} are outlined in Subsection~\ref{outline_subs}.
We suggest that the reader return to Subsections~\ref{descr_subs} and~\ref{outline_subs}
before going through the technical details of the blowup constructions in 
Sections~\ref{curvebl_sec}-\ref{map1bl_sec}.
In the next subsection, we also describe a natural sheaf over $\wt\M_{1,k}^0(\Pn,d)$
which is closely related to the sheaf $\pi_*\ev^* \O_{\Pn}(a)$ over 
$\ov\M_{1,k}^0(\Pn,d)$.
It is shown to be locally free in Section~\ref{cone_sec}.
Finally, all the data necessary for applying the localization theorem of~\cite{ABo} 
to~$\wt\M_{1,k}^0(\Pn,d)$ and~$e(\wt\V_{1,k}^d)$ is given in Subsection~\ref{local_subs}.\\

\noindent
Through this article we work with Deligne-Mumford stacks.
They can also be thought of as analytic orbivarieties.
As we work with reduced scheme structures throughout the paper, 
we will call such objects simply varieties.
Also, all immersions will be assumed to be from smooth varieties.
(The notion of ``immersion'' is often called ``unramified'' in algebraic geometry.)\\

\noindent
The authors would like to thank Jun Li for many enlightening discussions.

\subsection{Description of the Desingularization}
\label{descr_subs}

\noindent
The moduli space $\ov\M_{1,k}(\Pn,d)$ has irreducible components of various dimensions.
One of these components is $\ov\M_{1,k}^0(\Pn,d)$, the closure of 
the stratum $\M_{1,k}^0(\Pn,d)$ of stable maps with smooth domains.
We now describe natural subvarieties of 
$\ov\M_{1,k}(\Pn,d)$\footnote{In fact, these will be substacks of the stack 
$\ov\M_{1,k}(\Pn,d)$. 
They can also be thought of as analytic sub-orbivarieties of
the analytic orbivariety $\ov\M_{1,k}(\Pn,d)$.
As we work with reduced scheme structures throughout the paper, 
we will call such objects simply varieties.}
which contain the remaining components of $\ov\M_{1,k}(\Pn,d)$.
They will be indexed by the~set
\begin{gather*}
\A_1(d,k) \equiv\big\{\si\!=\!(m;J_P,J_B)\!: m\!\in\!\Z^+,\, m\!\le\!d;~ 
[k]\!=\!J_P\!\sqcup\!J_B\big\},\\
\hbox{where}\qquad
[k]=\{1,\ldots,k\}.
\end{gather*}
For each $\si\!\in\!\A_1(d,k)$, let $\M_{1,\si}(\Pn,d)$
be the subset of $\ov\M_{1,k}(\Pn,d)$ consisting of the stable maps $[\cC,u]$
such that $\cC$ is a smooth genus-one curve~$E$ with $m$ smooth rational 
components attached directly to~$E$, $u|_E$ is constant,
the restriction of~$u$ to each rational component is non-constant,
and the marked points on $E$ are indexed by the set~$J_P$.
Here $P$ stands for ``principal component'', $B$ stands for ``bubble
component'', and $\A$ stands for ``admissible set''.
Figure~\ref{m3_fig1} shows the domain of an element of $\M_{1,\si}(\Pn,d)$,
where $\si\!=\!(3;\{2\},\{1\})$,
from the points of view of symplectic topology and of algebraic geometry.
In the first diagram, each shaded disc represents a sphere;
the integer next to each rational component $\cC_i$ indicates the degree of~$u|_{\cC_i}$.
In the second diagram, the components of $\cC$ are represented by curves,
and the pair of integers next to each component $\cC_i$ shows 
the genus of $\cC_i$ and the degree of~$u|_{\cC_i}$.
In both diagrams, the marked points are labeled in bold face.
Let $\ov\M_{1,\si}(\Pn,d)$ be the closure of $\M_{1,\si}(\Pn,d)$
in $\ov\M_{1,k}(\Pn,d)$.
The space $\ov\M_{1,\si}(\Pn,d)$ has a number of irreducible components.
These components are indexed by the splittings of the degree~$d$
into $m$ positive integers and of the set~$J_B$ into $m$~subsets.
However, we do not need to distinguish these components.\\

\noindent
It is straightforward to check that  
$$\ov\M_{1,k}(\Pn,d) =\ov\M^0_{1,k}(\Pn,d)\cup 
\bigcup_{\si\in\A_1(d,k)}\!\!\!\!\!\ov\M_{1,\si}(\Pn,d).$$
Dimensional considerations imply that if $\si\!=\!(m;J_P,J_B)\!\in\!A_1(d,k)$ 
and $m\!\le\!n$, then $\ov\M_{1,\si}(\Pn,d)$ is a union of components
of $\ov\M_{1,k}(\Pn,d)$.
The converse holds as well: 
$\ov\M_{1,\si}(\Pn,d)$ is contained in $\ov\M_{1,k}^0(\Pn,d)$  if $m\!>\!n$
by \cite[Theorem~\ref{g1comp-str_thm}]{g1comp}.
However, we will use the entire collection $\A_1(d,k)$ of subvarieties
of $\ov\M_{1,k}(\Pn,d)$ to construct $\wt\M_{1,k}^0(\Pn,d)$.
The independence of the indexing set $\A_1(d,k)$ of~$n$ leads
to a number of good properties being satisfied by our blowup construction;
see~(2) of Theorem~\ref{main_thm} and the second part of this subsection.
It may also be possible to use this construction to define reduced genus-one 
GW-invariants in algebraic geometry; this is achieved in symplectic topology
in~\cite{g1comp2}.\\

\begin{figure}
\begin{pspicture}(-1.1,-1.8)(10,1.25)
\psset{unit=.4cm}
\rput{45}(0,-4){\psellipse(5,-1.5)(2.5,1.5)\pscircle*(7.5,-1.5){.2}
\psarc[linewidth=.05](5,-3.3){2}{60}{120}\psarc[linewidth=.05](5,0.3){2}{240}{300}
\pscircle[fillstyle=solid,fillcolor=gray](5,-4){1}\pscircle*(5,-3){.2}\pscircle*(4,-4){.2}
\pscircle[fillstyle=solid,fillcolor=gray](6.83,.65){1}\pscircle*(6.44,-.28){.2}
\pscircle[fillstyle=solid,fillcolor=gray](3.17,.65){1}\pscircle*(3.56,-.28){.2}}
\rput(7,.4){\smsize{$\bf 2$}}\rput(5.2,-4.2){\smsize{$\bf 1$}}
\rput(.2,-.9){$d_1$}\rput(3.1,2.3){$d_2$}\rput(7.8,-2.5){$d_3$}
\psarc(15,-1){3}{-60}{60}\pscircle*(16.93,1.3){.2}
\rput(17.6,1.4){\smsize{$\bf 2$}}
\psline(17,-1)(22,-1)\psline(16.8,-2)(21,-3)\psline(16.8,0)(21,1)
\rput(15.2,-3.5){\smsize{$(1,0)$}}\rput(22.4,1){\smsize{$(0,d_1)$}}
\rput(23.4,-1){\smsize{$(0,d_2)$}}\rput(22.4,-3){\smsize{$(0,d_3)$}}
\pscircle*(18.9,-2.5){.2}\rput(18.6,-3){\smsize{$\bf 1$}}
\rput(33,-1){\smsize{$\begin{array}{l}
\si\!=\!(3;\{2\},\{1\})\\  \\ d_1\!+\!d_2\!+\!d_3\!=\!d \\ d_1,d_2,d_3\!>\!0
\end{array}$}}
\end{pspicture}
\caption{The Domain of an Element of $\M_{1,\si}(\Pn,d)$}
\label{m3_fig1}
\end{figure}

\noindent
We define a partial ordering  $\prec$  on the set $\A_1(d,k)$ by
\begin{equation}\label{map1ord_e}
\si'\!\equiv\!\big(m';J_P',J_B'\big)\prec \si\!\equiv\!\big(m;J_P,J_B\big)
\qquad\hbox{if}\quad
\si'\!\neq\!\si,~~~ m'\!\le\!m, ~~\hbox{and}~~ J_P'\!\subseteq\!J_P.
\end{equation}
This relation is illustrated in Figure~\ref{map1ord_fig}, where an element
$\si$ of $\A_1(d,k)$ is represented by an element of the corresponding space
$\M_{1,\si}(\Pn,d)$.
We indicate that the degree of the stable map on every bubble component
is positive by shading the disks in the figure.
We show only the marked points lying on the principal component.
The exact distribution of the remaining marked points between the components is 
irrelevant.\\

\begin{figure}
\begin{pspicture}(-1.1,-1.8)(10,1.25)
\psset{unit=.4cm}
\rput{45}(0,-4){\psellipse(5,-1.5)(2.5,1.5)
\psarc[linewidth=.05](5,-3.3){2}{60}{120}\psarc[linewidth=.05](5,0.3){2}{240}{300}
\pscircle[fillstyle=solid,fillcolor=gray](5,-4){1}\pscircle*(5,-3){.2}
\pscircle[fillstyle=solid,fillcolor=gray](6.83,.65){1}\pscircle*(6.44,-.28){.2}}
\rput(11.5,-1.5){\begin{Huge}$\prec$\end{Huge}}
\rput{45}(15,-4){\psellipse(5,-1.5)(2.5,1.5)\pscircle*(7.5,-1.5){.2}
\psarc[linewidth=.05](5,-3.3){2}{60}{120}\psarc[linewidth=.05](5,0.3){2}{240}{300}
\pscircle[fillstyle=solid,fillcolor=gray](5,-4){1}\pscircle*(5,-3){.2}
\pscircle[fillstyle=solid,fillcolor=gray](3.17,.65){1}\pscircle*(3.56,-.28){.2}}
\rput(22,.4){\smsize{$\bf 2$}}
\rput{45}(25,-4){\psellipse(5,-1.5)(2.5,1.5)
\psarc[linewidth=.05](5,-3.3){2}{60}{120}\psarc[linewidth=.05](5,0.3){2}{240}{300}
\pscircle[fillstyle=solid,fillcolor=gray](5,-4){1}\pscircle*(5,-3){.2}
\pscircle[fillstyle=solid,fillcolor=gray](6.83,.65){1}\pscircle*(6.44,-.28){.2}
\pscircle[fillstyle=solid,fillcolor=gray](3.17,.65){1}\pscircle*(3.56,-.28){.2}}
\end{pspicture}
\caption{Examples of Partial Ordering~\e_ref{map1ord_e}}
\label{map1ord_fig}
\end{figure}

\noindent
Choose an ordering $<$ on $\A_1(d,k)$ extending the partial ordering~$\prec$.
The desingularization
\begin{equation*}\label{desing_e2}
\ti\pi\!: \wt\M_{1,k}^0(\Pn,d)\lra\ov\M_{1,k}^0(\Pn,d)
\end{equation*}
is constructed by blowing up $\ov\M_{1,k}(\Pn,d)$
along the subvarieties $\ov\M_{1,\si}(\Pn,d)$ and their proper transforms
in the order specified by~$<$.
In other words, we first blow up $\ov\M_{1,k}(\Pn,d)$ along $\ov\M_{1,\si_{\min}}(\Pn,d)$,
where 
$$\si_{\min}\equiv(1;\eset,[k])$$ 
is the smallest element of~$\A_1(d,k)$.
We then blow up the resulting space along the proper transform of
$\ov\M_{1,\si_2}(\Pn,d)$, where $\si_2$ is the smallest element 
of $\A_1(d,k)\!-\!\{\si_{\min}\}$.
We continue this procedure until we blow up along the proper transform
of $\ov\M_{1,\si_{\max}}(\Pn,d)$, where 
$$\si_{\max}=(d;[k],\eset)$$ 
is the largest element of~$\A_1(d,k)$.
The variety resulting from this last blowup is the proper transform
$\wt\M_{1,k}^0(\Pn,d)$ of $\ov\M_{1,k}^0(\Pn,d)$,
as all other irreducible components of $\ov\M_{1,k}(\Pn,d)$
have been ``blown out of existence''.\\

\noindent
The first interesting case of this construction, i.e.~for $\ov\M_{1,0}^0(\P^2,3)$,
is described in detail in~\cite{summ}.
The space $\wt\M_{1,0}^0(\P^2;3)$ is a smooth compactification of the space
of smooth plane cubics. 
It has a richer structure than the naive compactification, $\P^9$, does.

\begin{thm}
\label{main_thm}
Suppose $n,d\!\in\!\Z^+$, $k\!\in\!\bar\Z^+$, $<$ is an ordering on the set
$\A_1(d,k)$ extending the partial ordering~$\prec$, and
$$\ti\pi\!: \wt\M_{1,k}^0(\Pn,d)\lra\ov\M_{1,k}^0(\Pn,d)$$ is 
the blowup of $\ov\M_{1,k}^0(\Pn,d)$ obtained by blowing up 
$\ov\M_{1,k}(\Pn,d)$ along the subvarieties $\ov\M_{1,\si}(\Pn,d)$ and 
their proper transforms in the order specified by~$<$.\\
(1) The variety $\wt\M_{1,k}^0(\Pn,d)$ is smooth and is independent of the choice
of ordering~$<$ extending~$\prec$.\\ 
(2) For all $m\!\le\!n$, the embedding $\ov\M_{1,k}^0(\P^m,d)\!\lra\!\ov\M_{1,k}^0(\Pn,d)$
lifts to an embedding
$$\wt\M_{1,k}^0(\P^m,d)\lra\wt\M_{1,k}^0(\Pn,d)$$
and the image of the latter embedding
is the preimage of  $\ov\M_{1,k}^0(\P^m,d)$ under $\ti\pi$.\\
(3) The blowup locus at every step of the blowup construction is a smooth subvariety
in the corresponding blowup of~$\ov\M_{1,k}(\Pn,d)$.\\
(4) All fibers of $\ti\pi$ are connected.
\end{thm}

\noindent {\it Remark:} While in Section~\ref{map1bl_sec} we analyze
the blowup construction starting with the reduced scheme structure on
$\ov\M_{1,k}(\Pn,d)$, Theorem~\ref{main_thm} applies to the standard
scheme structure on $\ov\M_{1,k}(\Pn,d)$ as well.  It is known that
$\M_{1,k}^0(\Pn,d)$ is a smooth stack (as such maps are
unobstructed, see for example \cite[Prop.~5.5(c)]{enumtang}).  Thus, its scheme-theoretic closure,
$\ov\M_{1,k}^0(\Pn,d)$, is reduced.  During the blowup process all
other components of $\ov\M_{1,k}(\Pn,d)$ are ``blown
out of existence'', as is any non-reduced scheme structure.\\

\noindent
In Theorem~\ref{main_thm} and throughout the rest of the paper we denote by 
$\bar\Z^+$ the set of nonnegative integers.
We analyze the sequential blowup construction of Theorem~\ref{main_thm}
in Section~\ref{map1bl_sec} using the inductive assumptions ($I1$)-($I15$)
of Subsection~\ref{map1blconstr_subs}.
One of these assumptions,~($I3$), implies the second part of 
the first statement of Theorem~\ref{main_thm}, 
as different choices of an ordering~$<$ extending the partial ordering~$\prec$
correspond to different orders of blowups along subvarieties that are disjoint.
For example, suppose
$$\si_{\eset}=\big(2;\eset,\{1,2\}\big), \qquad
\si_1=\big(2;\{1\},\{2\}\big), \quad\hbox{and}\quad
\si_2=\big(2;\{2\},\{1\}\big).$$
While $\ov\M_{1,\si_1}(\Pn,d)$ and $\ov\M_{1,\si_2}(\Pn,d)$ do intersect
in $\ov\M_{1,2}(\Pn,d)$, their proper transforms are disjoint
after the blowup along $\ov\M_{1,\si_{\eset}}(\Pn,d)$.
The second statement of Theorem~\ref{main_thm} follows immediately from 
the description of the blowup construction in this and the next subsections,
as each step of the construction commutes with the embeddings of the moduli spaces
induced by the embedding $\P^m\!\lra\!\Pn$.\\

\noindent
The main claim of this paper is that $\wt\M_{1,k}^0(\Pn,d)$ is a smooth variety.
The structure of the space of maps from 
curves with no contracted  genus-one components 
$$\M_{1,k}^{\eff}(\Pn,d) \equiv  \ov\M_{1,k}(\Pn,d)
- \bigcup_{\si\in\A_1(d,k)}\!\!\!\!\! \ov\M_{1,\si}(\Pn,d)$$
is well understood; see \cite[Chapter~27]{H} for example.
In particular, $\M_{1,k}^{\eff}(\Pn,d)$ is smooth.
Below we describe the structure of the complement 
$\part\wt\M_{1,k}^0(\Pn,d)$ of $\M_{1,k}^{\eff}(\Pn,d)$ in~$\wt\M_{1,k}^0(\Pn,d)$.\\ 

\noindent
If $J$ is a finite set and $g$ is a nonnegative integer, 
we denote by $\ov\cM_{g,J}$ the moduli space of stable genus-$g$ curves
with $|J|$ marked points, which are indexed by the set~$J$. 
Similarly, we denote by $\ov\M_{g,J}(\Pn,d)$ 
the moduli space of stable maps from genus-$g$ curves 
with marked points indexed by~$J$ to~$\Pn$.
If $j\!\in\!J$, let
$$\ev_j\!: \ov\M_{g,J}(\Pn,d)\lra\Pn$$
be the evaluation map at the marked point labeled by~$j$.\\

\noindent
If $\si\!=\!(m;J_P,J_B)$ is an element of~$\A_1(d,k)$, we define
\begin{alignat*}{1}
\ov\cM_{\si;P} &= \ov\cM_{1,[m]\sqcup J_P} \qquad\hbox{and}\\
\ov\M_{\si;B}(\Pn,d)&= 
\Big\{(b_1,\ldots,b_m)\in\prod_{i=1}^m\ov\M_{0,\{ 0 \}\sqcup J_i}(\Pn,d_i)\!:
d_i\!>\!0,~\sum_{i=1}^md_i\!=\!d;~ \bigsqcup_{i=1}^m J_i\!=\!J_B;\\
&\qquad\qquad\qquad\qquad\qquad\qquad\qquad\qquad\qquad\qquad~~
\ev_0(b_{i_1})\!=\!\ev_0(b_{i_2})~\forall\, i_1,i_2\!\in\![m]\Big\}.
\end{alignat*}
There is a natural node-identifying surjective immersion
$$\io_{\si}\!: \ov\cM_{\si;P}\!\times\!\ov\M_{\si;B}(\Pn,d)
\lra \ov\M_{1,\si}(\Pn,d) \subset \ov\M_{1,k}(\Pn,d).$$
As before, $P$ denotes ``principal component'', and $B$ denotes ``bubble components''.
This immersion descends to the quotient:
$$\bar\io_{\si}\!: \big(\ov\cM_{\si;P}\!\times\!\ov\M_{\si;B}(\Pn,d)\big) \big/G_{\si}
\lra \ov\M_{1,\si}(\Pn,d),$$
where $G_{\si}\!\equiv\!S_m$ is the symmetric group on $m$ elements.
If $m\!\ge\!3$, $\bar\io_{\si}$ is not an isomorphism as 
some subvarieties of the left side are identified.
An example of a point on the right which is the image of two points 
on the left is given in Figure~\ref{doublept_fig}.
In addition to the conventions used in Figure~\ref{m3_fig1},
in the first, symplectic-topology, diagram of Figure~\ref{doublept_fig}
we leave the components of the domain on which the map is constant unshaded.
The subvarieties identified by the map~$\bar\io_{\si}$ get ``unidentified'' 
after taking the proper transform of $\ov\M_{1,\si}(\Pn,d)$
in the blowup of $\ov\M_{1,k}(\Pn,d)$ at the step corresponding~to
$$\ti\si=\max\big\{\si'\!\in\!\A_1(d,k)\!: \si'\!\prec\!\si\big\}.$$
This is insured by the inductive assumption~($I13$) in Subsection~\ref{map1blconstr_subs}
and implies the third statement of Theorem~\ref{main_thm}.
For example, if $m\!=\!3$ and $k\!=\!0$ as in Figure~\ref{doublept_fig}, 
the ``identified'' subvarieties are ``unidentified'' after the blowup of 
the proper transform of $\ov\M_{1,(2;\eset,\eset)}(\Pn,d)$.\\

\begin{figure}
\begin{pspicture}(-1.1,-1.8)(10,1.25)
\psset{unit=.4cm}
\psellipse(5,-1.5)(1.5,2.5)
\psarc[linewidth=.05](6.8,-1.5){2}{150}{210}\psarc[linewidth=.05](3.2,-1.5){2}{330}{30}
\pscircle(2.5,-1.5){1}\pscircle*(3.5,-1.5){.2}
\pscircle(7.5,-1.5){1}\pscircle*(6.5,-1.5){.2}
\pscircle[fillstyle=solid,fillcolor=gray](1.09,-.09){1}\pscircle*(1.79,-.79){.2}
\pscircle[fillstyle=solid,fillcolor=gray](1.09,-2.91){1}\pscircle*(1.79,-2.21){.2}
\pscircle[fillstyle=solid,fillcolor=gray](8.91,-.09){1}\pscircle*(8.21,-.79){.2}
\pscircle[fillstyle=solid,fillcolor=gray](8.91,-2.91){1}\pscircle*(8.21,-2.21){.2}
\rput(-.4,-2.5){$d_1$}\rput(-.4,.4){$d_2$}
\rput(10.4,-2.5){$d_3$}\rput(10.4,.4){$d_4$}
\psarc(15,-1){3}{-60}{60}\rput(15.4,-4){\smsize{$(1,0)$}}
\psline(16.8,.4)(21,.4)\rput(15.6,.4){\smsize{$(0,0)$}}
\psline(16.8,-2.4)(21,-2.4)\rput(15.6,-2.4){\smsize{$(0,0)$}}
\psline(18.5,-.4)(20,2.5)\rput(21.3,2.4){\smsize{$(0,d_1)$}}
\psline(19.5,.6)(22.5,-.4)\rput(23.8,-.2){\smsize{$(0,d_2)$}}
\psline(19.5,-2.6)(22.5,-1.6)\rput(23.8,-1.8){\smsize{$(0,d_3)$}}
\psline(18.5,-1.6)(20,-4.5)\rput(21.3,-4.4){\smsize{$(0,d_4)$}}
\rput(32.5,-1){\smsize{$\begin{array}{l}
\si\!=\!(3;\eset,\eset)\\  \\ d_1\!+\!d_2\!+\!d_3\!+\!d_4\!=\!d \\ d_1,d_2,d_3,d_4\!>\!0
\end{array}$}}
\end{pspicture}
\caption{A Point in $\ov\M_{1,\si}(\Pn,d)\!\subset\!\ov\M_{1,0}(\Pn,d)$ 
with Two Preimages under $\bar\io_{\si}$}
\label{doublept_fig}
\end{figure}

\noindent
{\it Remark:} Throughout the paper, we use $\M$ (fraktur font) to denote moduli spaces
of stable {\it maps}, of genus-zero or one, into~$\Pn$.
We use $\cM$ (calligraphic font) to denote moduli spaces of stable {\it curves}.\\

\noindent
For each $i\!\in\![m]$, let 
$$\pi_i\!: \ov\M_{\si;B}(\Pn,d)\lra
\!\!\bigsqcup_{d_i>0,J_i\subset J_B}\!\!\!\!\!\!\ov\M_{0,\{0\}\sqcup J_i}(\Pn,d_i)$$
be the natural projection onto the $i$th component.
We put
$$F_{\si;B}=\bigoplus_{i=1}^m \pi_i^*L_0,$$
where $L_0\!\lra\!\ov\M_{0,\{0\}\sqcup J_i}(\Pn,d_i)$
is the universal tangent line bundle for the marked point $0$.
In Subsections~\ref{curve1bl_subs} and~\ref{map0blconstr_subs}, we construct blowups
\label{ravi1}
\begin{alignat*}{1}
\ti\pi_{\si;P}\!\equiv\!\pi_{1,([m],J_P)}
&:\wt\cM_{\si;P}\!\equiv\!\wt\cM_{1,([m],J_P)}
\lra \ov\cM_{\si;P}\!\equiv\!\ov\cM_{1,[m]\sqcup J_P}
\qquad\hbox{and}\\
\ti\pi_{\si;B}\!\equiv\!\pi_{0,([m],J_B)}
&:\wt\M_{\si;B}(\Pn,d)\!\equiv\!\wt\M_{0,([m],J_B)}(\Pn,d)
 \lra \P F_{\si;B}\!\equiv\!\P F_{([m],J_B)}.
\end{alignat*}
We also construct a section 
\begin{equation}\label{spacesec_e1}
\wt\cD_{\si;B}\!\equiv\!\wt\cD_{([m],J_B)}
\in  \Ga\big(\wt\M_{\si;B}(\Pn,d);
\E_{\si;B}^*\!\otimes\!\ti\pi_{\si;B}^*\pi_{\P F_{\si;B}}^*\ev_0^*T\Pn\big),
\end{equation}
where
$$\ev_0\!: \ov\M_{\si;B}(\Pn,d)\lra\Pn  \qquad\hbox{and}\qquad
\pi_{\P F_{\si;B}}\!: \P F_{\si;B}\lra\ov\M_{\si;B}(\Pn,d)$$
are the natural evaluation map and the bundle projection map, respectively, and
$$\E_{\si;B}\!\equiv\!\ti\E \lra
  \wt\M_{\si;B}(\Pn,d)\!\equiv\!\wt\M_{0,([m],J_B)}(\Pn,d)$$
is a line bundle.
This line bundle is the sum of the tautological line bundle
$$\ga_{\si;B} \lra \P F_{\si;B}$$
and all exceptional divisors. 
The section $\wt\cD_{\si;B}$ is transverse to the zero section.
Thus, its zero~set, 
\begin{equation}\label{spacedfn_e1}
\wt\cZ_{\si;B}(\Pn,d) \equiv \wt\cD_{\si;B}^{-1}(0) 
\subset \wt\M_{\si;B}(\Pn,d),
\end{equation}
is a smooth subvariety.
The boundary $\part\wt\M_{1,k}^0(\Pn,d)$ of $\wt\M_{1,k}^0(\Pn,d)$
is a union of smooth divisors:
$$\part\wt\M_{1,k}^0(\Pn,d)=\bigcup_{\si\in\A_1(d,k)}\!\!\!\!\! 
\wt\cZ_{\si}(\Pn,d)\big/G_{\si},
\qquad\hbox{where}\quad
\wt\cZ_{\si}(\Pn,d)=\wt\cM_{\si;P}\!\times\!\wt\cZ_{\si;B}(\Pn,d);$$
see the inductive assumptions~($I7$) and ($I8$) in Subsection~\ref{map1blconstr_subs}
and Figure~\ref{bdlift_fig}.
By the inductive assumption~($I6$) and~($I7$),
the normal bundle of $\wt\cZ_{\si}(\Pn,d)$ in $\wt\M_{1,k}(\Pn,d)$ is the quotient
of the line bundle
$$\L_{\si;P}\otimes \ti\pi_{\si;B}^*\ga_{\si;B}\lra
\wt\cM_{\si;P}\!\times\!\wt\cZ_{\si;B}(\Pn,d)$$
by the $G_{\si}$-action, where
$$\L_{\si;P}\!\equiv\!\L \lra  \wt\cM_{\si;P}\!\equiv\!\wt\cM_{1,([m],J_P)}$$
is the universal tangent line bundle constructed in Subsection~\ref{curve1bl_subs}.
Thus we conclude that 
$\wt\M_{1,k}^0(\Pn,d)$ is smooth, as the open subset 
$\M_{1,k}^{\eff}(\Pn,d)$ is smooth, and its complement is
a union of smooth divisors whose normal sheaves are line bundles
(i.e.\ with their reduced induced scheme structure, they are 
Cartier divisors).\\

\begin{figure}
\begin{pspicture}(-1.1,-1.8)(10,1.25)
\psset{unit=.4cm}
\rput(10,0){$\wt\cM_{\si;P}~\times~\wt\cZ_{\si;B}(\Pn,d)$}
\rput(25,0){$\wt\M_{1,k}^0(\Pn,d)$}
\rput(10,-4){$\ov\cM_{\si;P}~\times~\ov\M_{\si;B}(\Pn,d)$}
\rput(25,-4){$\ov\M_{1,k}(\Pn,d)$}
\psline{->}(16,0)(21,0)\psline{->}(16,-4)(21,-4)
\psline{->}(6.5,-1)(6.5,-3)\psline{->}(12.5,-1)(12.5,-3)\psline{->}(25,-1)(25,-3)
\rput(18.5,.7){\small{$\ti\io_{\si}$}}\rput(18.5,-3.4){\small{$\io_{\si}$}}
\rput(7.6,-2){\small{$\ti\pi_{\si;P}$}}\rput(13.6,-2){\small{$\ti\pi_{\si;B}$}}
\rput(25.5,-2){\small{$\ti\pi$}}
\end{pspicture}
\caption{Changes in the Boundary Structure of $\ov\M_{1,k}^0(\Pn,d)$
under the Desingularization}
\label{bdlift_fig}
\end{figure}

\noindent
{\it Remark 1:} In the Gromov-Witten theory, the symbol $\E$ is commonly used 
to denote the Hodge vector bundle of holomorphic differentials.
It is the zero vector bundle in the genus-zero case.
The line bundles over moduli spaces of genus-zero curves and maps
we denote by $\E$, with various decorations, play roles analogous
to that of the Hodge line bundle over moduli spaces of genus-one curves.
The most overt parallel is described at the end of Subsection~\ref{curvebldata_subs}.
There are deeper, more subtle, connections as well;
compare the structural descriptions of Lemmas~\ref{deriv0str_lmm} 
and~\ref{map1bl_lmm2}, for example.\\
 
\noindent
{\it Remark 2:} Throughout this paper, the symbols $\cD$ and $\D$,
with various decorations, denote vector bundle sections related 
to derivatives of holomorphic maps into $\Pn$ and of holomorphic bundle sections.
In most cases, such bundle sections are viewed as vector bundle homomorphisms.\\

\noindent
The final claim of Theorem~\ref{main_thm} follows from the 
fact that  $\ov\M_{1,k}^0(\Pn,d)$ is unibranch (locally irreducible).
If $\pi\!:Y\!\lra\!X$  is a surjective birational map of irreducible
varieties, and $\pi^{-1}(x)$ is not connected for some $x\!\in\!X$, 
then  $X$ is not unibranch at $x$.\\

\noindent
We next describe a desingularization of the sheaf $\pi_*\ev^*\O_{\Pn}(a)$ and 
of the corresponding cone $\V_{1,k}^d$ over $\ov\M_{1,k}^0(\Pn,d)$.
Let $\wt\U\!=\!\ti\pi^*\U$ be the pullback of $\U$ by~$\ti\pi$:
$$\xymatrix{\wt\U \ar[d]^{\pi} \ar[r]^{\ti\pi}& \U \ar[d]^{\pi} \ar[r]^{\ev} & \Pn \\
\wt\M_{1,k}^0(\Pn,d) \ar[r]^{\ti\pi}& \ov\M_{1,k}^0(\Pn,d).}$$
For each $\si\!\in\!\A_1(d,k)$, let
$$\V_{\si;B}\lra \ov\M_{\si;B}(\Pn,d)$$
be the cone induced by the sheaf $\O_{\Pn}(a)$, similarly to $\V_{g,k}^d$;
see Subsection~\ref{conehomomor_subs} for details.
It is a vector bundle of rank $da\!+\!1$.
We note that 
$$\ti\pi^*\V_{1,k}^d\big|_{\wt\cZ_{\si}(\Pn,d)} 
=\pi_B^*\big\{\ti\pi_{\si;B}^*
\pi_{\P F_{\si;B}}^*\V_{\si;B}|_{\wt\cZ_{\si;B}(\Pn,d)}\big\}\big/G_{\si},$$
where 
$$\pi_B\!:\wt\cM_{\si;P}\!\times\!\wt\cZ_{\si;B}(\Pn,d)\lra
\wt\cZ_{\si;B}(\Pn,d)$$
is the projection map. Let $\cL\!=\!\ga^{*\otimes a}$, where
$\ga\!\lra\!\Pn$ is the tautological line bundle.

\begin{thm}
\label{cone_thm}
Suppose $d,n,a\!\in\!\Z^+$ and $k\!\in\!\bar\Z^+$.\\
(1) The sheaf $\pi_*\ti\pi^*\ev^*\O_{\Pn}(a)$ over $\wt\M_{1,k}^0(\Pn,d)$
is locally free and of the expected rank, i.e.~$da$.\\
(2) If $\wt\V_{1,k}^d\!\subset\!\ti\pi^*\V_{1,k}^d$ is the corresponding vector bundle
and $\si\!\in\!\A_1(d,k)$, then there exists a surjective bundle homomorphism
$$\wt\D_{\si;B}\!: \ti\pi_{\si;B}^*\pi_{\P F_{\si;B}}^*\V_{\si;B}|_{\wt\cZ_{\si;B}(\Pn,d)}
\lra \E_{\si;B}^*\!\otimes\!\ti\pi_{\si;B}^*\pi_{\P F_{\si;B}}^*\ev_0^*\cL$$
over $\wt\cZ_{\si;B}(\Pn,d)$ such that 
$$\wt\V_{1,k}^d\big|_{\wt\cZ_{\si}(\Pn,d)}
= \big(\pi_B^*\, \ker \wt\D_{\si;B}\big) \big/G_{\si}.$$
(3) $\ti\pi_*\pi_*\ti\pi^*\ev^*\O_{\Pn}(a)=\pi_*\ev^*\O_{\Pn}(a)$
 over $\ov\M_{1,k}^0(\Pn,d)$.
\end{thm}

\noindent
We prove the first two statements of this theorem by working with 
the cone
$$p\!:\ov\M_{1,k}(\cL,d)\lra \ov\M_{1,k}(\Pn,d).$$
The sheaves $\pi_*\ev^*\O_{\Pn}(a)$ and $\pi_*\ti\pi^*\ev^*\O_{\Pn}(a)$  are the sheaves
of (holomorphic) sections of 
$$\V_{1,k}^d\equiv \ov\M_{1,k}(\cL,d)\big|_{\ov\M_{1,k}^0(\Pn,d)}
\lra \ov\M_{1,k}^0(\Pn,d)$$
and $\ti\pi^*\V_{1,k}^d$, respectively; see Lemma~\ref{conesheaf_lmm}.
In Subsection~\ref{conebl_subs}, we lift the blowup construction of 
Subsection~\ref{map0blconstr_subs} to $\ov\M_{1,k}(\cL,d)$.
In particular, we blow up $\ov\M_{1,k}(\cL,d)$ along the subvarieties
$$\ov\M_{1,\si}(\cL,d) = p^{-1}\big(\ov\M_{1,\si}(\Pn,d)\big),
\qquad\si\in\A_1(d,k),$$
and their proper transforms.
The end result of this construction, which we denote by $\wt\M_{1,k}^0(\cL,d)$,
is smooth for essentially the same reasons that $\ov\M_{1,k}^0(\Pn,d)$~is.
The only additional input we need is Lemma~\ref{cone1bl_lmm2}, 
which is a restatement of the key result concerning the structure
of the cone~$\V_{1,k}^d$ obtained in~\cite{g1cone}.
The bundle 
$$\ti{p}\!:\wt\M_{1,k}^0(\cL,d)\lra\wt\M_{1,k}^0(\Pn,d)$$
of vector spaces of the same rank contains $\wt\M_{1,k}^0(\Pn,d)$ as the zero section.
Thus, $\ti{p}$ is a vector bundle.
There is a natural inclusion
$$\wt\M_{1,k}^0(\cL,d) \lra \ti\pi^*\ov\M_{1,k}(\cL,d).$$
All sections of $\ti\pi^*\ov\M_{1,k}(\cL,d)$ must in fact be sections of
$\wt\M_{1,k}^0(\cL,d)$ and thus the sheaf $\pi_*\ti\pi^*\ev^*\O_{\Pn}(a)$ is indeed 
locally free.
The bundle map 
$$\wt\D_{\si;B}\equiv\wt\D_{([m],J_B)}$$ 
of the second statement of Theorem~\ref{cone_thm} is described in
Subsection~\ref{conehomomor_subs}.
It is the ``vertical'' part of the natural extension of the bundle map $\wt\cD_{\si;B}$
from stable maps into $\Pn$ to stable maps into~$\cL$.
Finally, the last statement of Theorem~\ref{cone_thm} is a consequence 
of the last statement of Theorem~\ref{main_thm}; see Lemma~\ref{pushfor_lmm3}.
At this point, this observation does not appear to have any applications though.\\

\noindent
{\it Remark:} By applying the methods of Section~\ref{cone_sec} and of \cite{g1cone},
it should be possible to show that the standard scheme structure on $\ov\M_{1,k}(\Pn,d)$
is in fact reduced.

\subsection{Outline of the Main Desingularization Construction}
\label{outline_subs}

\noindent
The main blowup construction of this paper is contained in
Subsections~\ref{map1prelim_subs} and~\ref{map1blconstr_subs}.
It is a sequence of {\it idealized} blowups along smooth subvarieties.
In other words, the blowup locus $\ov\M_{1,\si}^{\si-1}$ at each step comes with
an {\it idealized} normal bundle~$\N^{\ide}$.
It is a vector {\it bundle} (of the smallest possible rank) 
containing the normal cone $\N$ for~$\ov\M_{1,\si}^{\si-1}$.  
After taking the usual blowup of the ambient space along $\ov\M_{1,\si}^{\si-1}$,
we attach the {\it idealized exceptional divisor}
$$\cE^{\ide}\equiv\P\N^{\ide}$$
along the usual exceptional divisor 
$$\cE\equiv\P\N\subset \cE^{\ide}.$$
The blowup construction summarized in Theorem~\ref{main_thm} is contained in
the idealized blowup construction of Section~\ref{map1bl_sec}.
The latter turns out to be more convenient for describing the proper transforms
of $\ov\M_{1,k}^0(\Pn,d)$, including at the final stage, i.e.~$\wt\M_{1,k}^0(\Pn,d)$.\\

\noindent
The ambient space $\ov\M_{1,k}^{\si}$ at each step 
$\si\!\in\!\{0\}\!\sqcup\!\A_1(d,k)$ of the blowup construction 
contains a subvariety $\ov\M_{1,\si^*}^{\si}$ for each $\si^*\!\in\!\A_1(d,k)$.
We take $\ov\M_{1,\si}^{\si}$ to be the idealized exceptional divisor 
for the idealized blowup just constructed, i.e.~along~$\ov\M_{1,\si}^{\si-1}$.
If $\si^*\!<\!\si$ or $\si^*\!>\!\si$, $\ov\M_{1,\si^*}^{\si}$ is the proper 
transform of $\ov\M_{1,\si^*}^{\si^*}$ or $\ov\M_{1,\si^*}(\Pn,d)$, respectively.\\

\noindent
Every immersion $\io_{\si^*}$ of Subsection~\ref{descr_subs}  
comes with an {\it idealized} normal bundle~$\N_{\io_{\si^*}}^{\ide}$.
It is a vector bundle of the smallest possible rank containing the normal cone
to the immersion~$\io_{\si}$ (see Definition~\ref{virtvar_e}).
It is given~by
$$\N_{\io_{\si^*}}^{\ide} = \bigoplus_{i\in[m^*]}\!\!
\pi_P^*L_i\!\otimes\!\pi_B^*\pi_i^*L_0
\qquad\hbox{if}\quad \si^*=(m^*;J_P^*,J_B^*),$$
where 
$$\pi_P,\pi_B\!:\ov\cM_{\si^*;P}\!\times\!\ov\M_{\si^*;B}(\Pn,d)
\lra \ov\cM_{\si^*;P},\ov\M_{\si^*;B}(\Pn,d)$$
are the component projection maps.
In the case of Figure~\ref{m3_fig1}, $\N_{\io_{\si^*}}^{\ide}$ is a rank-three vector bundle
encoding the potential smoothings of the three nodes.
At each step~$\si$ of the blowup construction, $\io_{\si^*}$ induces an immersion
$\io_{\si,\si^*}$ onto~$\ov\M_{1,\si^*}^{\si}$.
Like the domain of~$\io_{\si^*}$, the domain of $\io_{\si,\si^*}$ splits
as a Cartesian product.
If $\si^*\!>\!\si$, the second component of the domain does not change from 
the previous step, while the first is modified by blowing up along a collection of 
disjoint subvarieties, as specified by the inductive assumption~($I9$)
in Subsection~\ref{map1blconstr_subs}.
The idealized normal bundle $\N_{\io_{\si,\si^*}}^{\ide}$ is obtained from 
$\N_{\io_{\si-1,\si^*}}^{\ide}$ by twisting the first factor in each summand by
a subset of the exceptional divisors,  as specified by the inductive assumption~($I11$).
These blowup and twisting procedures correspond to several interchangeable steps 
in the blowup construction of Subsection~\ref{curve1bl_subs}.
For $\si^*\!=\!\si$, the first component in the domain of $\io_{\si-1,\si}$
has already been blown up all the way to $\wt\cM_{\si;P}$ and the first component
of every summand of $\N_{\io_{\si-1,\si}}^{\ide}$ has already twisted to 
the universal tangent line bundle~$\L$,~i.e.
$$\N_{\io_{\si-1,\si}}^{\ide}= 
\bigoplus_{i\in[m]} \pi_P^*\L \!\otimes\!
\pi_B^*\pi_i^*L_0
=\pi_P^*\L\otimes \pi_B^*F_{\si;B}
\lra \wt\cM_{\si;P}\!\times\!\ov\M_{\si;B}(\Pn,d),$$
if $\si\!=\!(m;J_P,J_B)$.
In particular, the domain for $\io_{\si,\si}$,
$$\P\N_{\io_{\si-1,\si}}^{\ide}=\wt\cM_{\si;P}\!\times\!\P F_{\si;B},$$
still splits as a Cartesian product!
The idealized normal bundle for $\io_{\si,\si}$ is the tautological line for
$\P\N_{\io_{\si-1,\si}}^{\ide}$:
$$\N_{\io_{\si,\si}}^{\ide}=\ga_{\N_{\io_{\si-1,\si}}^{\ide}}
=\pi_P^*\L\!\otimes\!\pi_B^*\ga_{F_{\si;B}} 
\equiv \pi_P^*\L\!\otimes\!\pi_B^*\ga_{\si;B}.$$
On the other hand, if $\si^*\!<\!\si$, the domain of $\io_{\si,\si^*}$ is obtained
from the domain of $\io_{\si-1,\si^*}$ by blowing up the second component along
a collection of disjoint subvarieties, as specified by the inductive assumption~($I4$)
in Subsection~\ref{map1blconstr_subs}.
This corresponds to several interchangeable steps of the blowup construction
in Subsection~\ref{map0blconstr_subs}.
By the time we are done with the last step of the blowup construction in 
Subsection~\ref{map1blconstr_subs}, $\P F_{\si;B}$ has been blown up all the way
to $\wt\M_{\si;B}(\Pn,d)$.
In the $\si^*\!<\!\si$ case,
$$\N_{\io_{\si,\si^*}}^{\ide} = \N_{\io_{\si-1,\si^*}}^{\ide},$$
since $\ov\M_{1,\si^*}^{\si-1}$ is transverse to~$\ov\M_{1,\si}^{\si-1}$.\\

\noindent
We study the proper transform $\ov\M_{1,(0)}^{\si}$ of $\ov\M_{1,k}^0(\Pn,d)$
in $\ov\M_{1,k}^{\si}$ by looking at the structure~of
$$\bar\cZ_{\si^*}^{\si}=\io_{\si,\si^*}^{\,-1}\big(\ov\M_{1,(0)}^{\si}\big).$$
Given a finite set $J$, there are natural bundle sections
$$s_j\in\Ga(\ov\cM_{1,J};L_j^*\!\otimes\!\E^*),~~~j\!\in\!J, 
\qquad\hbox{and}\qquad 
\cD_0\in\Ga\big(\ov\M_{0,\{0\}\sqcup J}(\Pn,d);L_0^*\!\otimes\!\ev_0^*T\Pn\big);$$
see Subsection~\ref{curvebldata_subs} and~\ref{map0str_subs}, respectively.
By Lemma~\ref{map1bl_lmm2}, the intersection of 
$$\bar\cZ_{\si^*}^0\equiv\io_{\si^*}^{\,-1}\big(\ov\M_{1,k}^0(\Pn,d)\big)$$
with the main stratum $\cM_{\si^*;P}\!\times\!\M_{\si^*;B}(\Pn,d)$ of
$\ov\cM_{\si^*;P}\!\times\!\ov\M_{\si^*;B}(\Pn,d)$ is
$$\cZ_{\si^*}^0=\big\{b\!\in\!\cM_{\si^*;P}\!\times\!\M_{\si^*;B}(\Pn,d)\!:
\ker\cD_{\si^*}|_b\!\neq\!\{0\}\big\}$$
where
\begin{gather*}
\cD_{\si^*} \in \Ga\big(\ov\cM_{\si^*;P} \!\times\!\ov\M_{\si^*;B}(\Pn,d);
\Hom(\N_{\io_{\si^*}}^{\ide},\pi_P^*\E^*\!\otimes\!\pi_B^*\ev_0^*T\Pn)\big),\\
\cD_{\si^*}\big|_{\pi_P^*L_i\otimes\pi_B^*\pi_i^*L_0}
=\pi_P^*s_i\!\otimes\!\pi_B^*\pi_i^*\cD_0, 
\qquad\forall\,i\!\in\![m^*].
\end{gather*}
In addition, if $\N\bar\cZ_{\si^*}^{\si}\!\subset\!\N_{\io_{\si,\si^*}}^{\ide}$ 
is the normal cone for the immersion $\io_{\si,\si^*}|_{\bar\cZ_{\si^*}^{\si}}$ 
into $\ov\M_{1,(0)}^{\si}$, then
$$\N\bar\cZ_{\si^*}^0\big|_{\cZ_{\si^*}^0}=\ker \cD_{\si^*}\big|_{\cZ_{\si^*}^0}$$
and $\N\bar\cZ_{\si^*}^0$ is the closure of $\N\bar\cZ_{\si^*}^0\big|_{\cZ_{\si^*}^0}$
in~$\N_{\io_{\si^*}}^{\ide}$.
By Lemma~\ref{virimmer_lmm3}, $\N\bar\cZ_{\si^*}^{\si}$ is still the closure of 
$\N\bar\cZ_{\si^*}^0\big|_{\cZ_{\si^*}^0}$, but now in~$\N_{\io_{\si,\si^*}}^{\ide}$,
for all $\si\!<\!\si^*$.
In Subsection~\ref{curve1bl_subs}, we construct a non-vanishing section
$$\ti{s}_i\in \Ga(\ov\cM_{1,J};\L^*\!\otimes\!\E^*) \approx \Ga(\ov\cM_{1,J};\C)$$
obtained by twisting $s_i$ by some exceptional divisors.
Since $\ti{s}_i$ agrees with $s_i$ on $\cM_{\si^*;P}$, we can replace $s_i$ with
$\ti{s}_i$ in the descriptions of~$\cD_{\si^*}$, $\cZ_{\si^*}^0$, and 
$\N\bar\cZ_{\si^*}^0\big|_{\cZ_{\si^*}^0}$ above.
In particular, $\N\bar\cZ_{\si^*}^{\si^*-1}$ is the closure~of
\begin{gather*}
\N\bar\cZ_{\si^*}^0\big|_{\cZ_{\si^*}^0}=
\pi_P^*\L\!\otimes\!\pi_B^*\ker\cD_{\si^*;B}\big|_{\cZ_{\si^*}^0}
\subset \pi_P^*\L\!\otimes\!\pi_B^*F_{\si^*;B}, \qquad\hbox{where}\\
\cD_{\si^*;B} \in \Ga\big(\ov\M_{\si^*;B}(\Pn,d);
\Hom(F_{\si^*;B},\ev_0^*T\Pn)\big), \qquad
\cD_{\si^*;B}\big|_{\pi_i^*L_0}=\pi_i^*\cD_0, \quad\forall\,i\!\in\![m^*].
\end{gather*}\\

\noindent
The bundle homomorphism $\cD_{\si^*;B}$ induces a section
$$\wt\cD_0\in\Ga\big(\P F_{\si^*;B};
\ga_{\si;B}^*\!\otimes\!\pi_{\P F_{\si^*;B}}^*\ev_0^*T\Pn\big).$$
By the previous paragraph and Lemma~\ref{virimmer_lmm3},
$\bar\cZ_{\si^*}^{\si^*}$ is the closure~of 
\begin{gather*}
\wt\cM_{\si^*;P} \times \wt\cD_0^{-1}(0)\!\cap\!\P F_{\si^*;B}\big|_{\M_{\si^*;B}(\Pn,d)}
\subset \wt\cM_{\si^*;P} \times\P F_{\si^*;B} \qquad\hbox{and}\\
\N\bar\cZ_{\si^*}^{\si^*}=
\pi_P^*\L\!\otimes\!\pi_B^*\ga_{\si;B}\big|_{\bar\cZ_{\si^*}^{\si^*}}.
\end{gather*}
Since $\wt\M_{1,k}^0(\Pn,d)\!\equiv\!\wt\M_{1,k}^{\si_{\max}}(\Pn,d)$
is the proper transform of $\ov\M_{1,(0)}^{\si^*}$ in $\ov\M_{1,k}^{\si_{\max}}$,
$$\wt\cZ_{\si^*}\equiv \io_{\si_{\max},\si^*}^{\,-1}\big(\wt\M_{1,k}^0(\Pn,d)\big)$$
is still the closure of 
$$\wt\cM_{\si^*;P} \times \wt\cD_0^{-1}(0)\!\cap\!\P F_{\si^*;B}\big|_{\M_{\si^*;B}(\Pn,d)}
\subset \wt\cM_{\si^*;P} \times\wt\M_{\si^*;B}(\Pn,d).$$
On the other hand, in the process of constructing the blowup $\wt\M_{\si;B}(\Pn,d)$
of $\P F_{\si^*;B}$ in Subsection~\ref{map0blconstr_subs}, we also define 
a bundle section
$$\wt\cD_{\si^*;B}\in
\Ga\big(\wt\M_{\si^*;B}(\Pn,d);
\E_{\si;B}^*\!\otimes\!\ti\pi_{\si;B}^*\pi_{\P F_{\si^*;B}}^*\ev_0^*T\Pn\big)$$
by twisting $\wt\cD_0$ by the exceptional divisors.
In particular,
$$\wt\cD_{\si^*;B}^{-1}(0)\!\cap\!\P F_{\si^*;B}\big|_{\M_{\si^*;B}(\Pn,d)}
 =\wt\cD_0^{-1}(0)\!\cap\!\P F_{\si^*;B}\big|_{\M_{\si^*;B}(\Pn,d)}.$$
Since $\wt\cD_{\si^*;B}$ is transverse to the zero set, we conclude that
$$\wt\cZ_{\si^*} = \wt\cM_{\si^*;P} \times \wt\cD_{\si^*;B}^{-1}(0),$$
as stated in Subsection~\ref{descr_subs}.\\

\noindent
Finally,
the role played by the blowup construction of Subsection~\ref{curve0bl_subs} 
in the blowup construction of Section~\ref{map0bl_sec} is similar to the role played
by the construction of Subsection~\ref{curve1bl_subs}  
in the construction of Section~\ref{map1bl_sec}.
In the case of Section~\ref{map0bl_sec}, we blow up a moduli space of genus-zero stable
maps, $\P F_{(\aleph,J)}$, along certain subvarieties $\wt\M_{0,\vr}^0$
and their proper transforms.
These subvarieties are images of natural node-identifying immersions~$\io_{0,\vr}$.
The domain of~$\io_{0,\vr}$ splits as the Cartesian product of a moduli space
of genus-zero curves and a moduli space of genus-zero maps, defined in
Subsections~\ref{curve0bl_subs} and~\ref{map0prelim_subs}, respectively.
As we modify $\wt\M_{0,\vr}^0$ by taking its proper transforms in the blowups
of $\P F_{(\aleph,J)}$ constructed in Subsection~\ref{map0blconstr_subs},
the first factor in the domain of the corresponding immersion changes 
by blowups along collections of smooth disjoint subvarieties, as specified
by the inductive assumption~($I6)$.
This change corresponds to several interchangeable steps in
the blowup construction of Subsection~\ref{curve0bl_subs}.
By the time we are ready to blow up the proper transform of $\wt\M_{0,\vr}^0$,
the first component of the domain of the corresponding immersion has been blown up
all the way to $\wt\cM_{0,\rho_P}$, the end result in the blowup construction
of Subsection~\ref{curve0bl_subs}.\\

\noindent
In the blowup construction of Subsection~\ref{map0blconstr_subs},
we twist a natural bundle section
$$\wt\cD_0\in\Ga\big(\P F_{(\ale,J)};
\ga_{(\aleph,J)}^*\!\otimes\!\pi_{\P F_{(\ale,J)}}^*\ev_0^*T\Pn\big)$$
by the exceptional divisors to a bundle section
$$\wt\cD_{(\aleph,J)}\in
\Ga\big(\wt\M_{(\ale,J)}(\Pn,d);
\wt\E^*\!\otimes\!\pi_{0,(\ale,J)}^*\pi_{\P F_{(\ale,J)}}^*\ev_0^*T\Pn\big).$$
The two sections enter in an essential way in the main blowup construction 
of this paper.
It is also essential that $\wt\cD_{(\aleph,J)}$ is transverse to the zero set.
The section~$\wt\cD_0$ is transverse to the zero set outside of the subvarieties
$\wt\M_{0,\vr}^0$ and vanishes identically along~$\wt\M_{0,\vr}^0$.
Its derivative in the normal direction to~$\io_{0,\vr}$ is described  
by Lemma~\ref{map0bl_lmm2}, using Lemma~\ref{deriv0str_lmm}. 
The bundle sections~$s_i$ over a moduli space of genus-zero curves defined in 
Subsection~\ref{curvebldata_subs} and modified in Subsection~\ref{curve0bl_subs} 
enter into the expression of Lemma~\ref{deriv0str_lmm}.
In fact, this expression is identical to the expression for $\cD_{\si^*}$ above,
i.e.~in the genus-one case.
We use Lemma~\ref{map0bl_lmm2} to show that with each newly twisted version of $\wt\cD_0$
is transverse to the zero set outside of the proper transforms of the remaining 
subvarieties $\ov\M_{0,\vr}^0$, i.e.~the ones that have not been blown up yet;
see the inductive assumption~($I4$) in Subsection~\ref{map0blconstr_subs}.
In particular, at the end of the blowup of Subsection~\ref{map0blconstr_subs},
we end up with a twisted version of $\wt\cD_0$, which we call $\wt\cD_{(\aleph,J)}$,
which is transverse to the zero set.

\subsection{Localization Data}
\label{local_subs}

\noindent
Suppose the group $G\!=\!(S^1)^{n+1}$ or $G\!=\!(\C^*)^{n+1}$ acts in
a natural way on the projective space~$\Pn$.
In particular, the fixed locus
consists of $n\!+\!1$ points, which we denote by $p_0,\ldots,p_n$, and
the only curves preserved by $G$ are the lines passing through pairs of fixed points.  
The $G$-action on $\Pn$ lifts to an action on $\ov\M_{1,k}(\Pn,d)$ and on 
$\wt\M_{1,k}^0(\Pn,d)$.  
The fixed loci of these two actions that are contained in $\M_{1,k}^{\eff}(\Pn,d)$ and
their normal bundles are the same and are described in \cite[Sects.~27.3 and~27.4]{H}.  
We note that the four-term exact sequence \cite[(27.6)]{H} applies to such loci.\\

\noindent
In this subsection, we describe the fixed loci of the $G$-action on 
$\wt\M_{1,k}^0(\Pn,d)$ that are contained in $\part\wt\M_{1,k}^0(\Pn,d)$
and their normal bundles.  
To simplify the discussion, we ignore all automorphism groups
until the very end of this subsection.\\

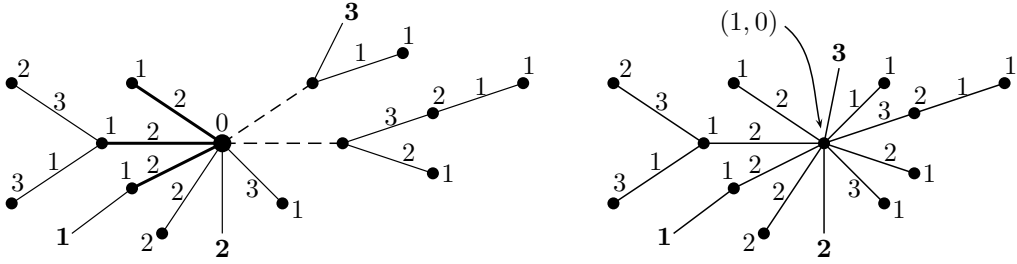
\begin{figure}
\begin{pspicture}(0,-2)(10,2)
\psset{unit=.4cm}
\pscircle*(10,0){.3}\rput(10,.7){\smsize{$0$}}
\psline[linewidth=.1](10,0)(7,2)\pscircle*(7,2){.2}
\rput(8.6,1.5){\smsize{$2$}}\rput(7.3,2.5){\smsize{$1$}}
\psline[linewidth=.1](10,0)(6,0)\pscircle*(6,0){.2}
\rput(7.7,.4){\smsize{$2$}}\rput(6.3,.6){\smsize{$1$}}
\psline[linewidth=.04](6,0)(3,2)\pscircle*(3,2){.2}
\rput(4.6,1.4){\smsize{$3$}}\rput(3.4,2.5){\smsize{$2$}}
\psline[linewidth=.04](6,0)(3,-2)\pscircle*(3,-2){.2}
\rput(4.4,-.6){\smsize{$1$}}\rput(3.2,-1.3){\smsize{$3$}}
\psline[linewidth=.1](10,0)(7,-1.5)\pscircle*(7,-1.5){.2}
\rput(7.7,-.7){\smsize{$2$}}\rput(6.8,-.9){\smsize{$1$}}
\psline[linewidth=.04](7,-1.5)(5,-3)\rput(4.7,-3.1){\smsize{$\bf 1$}}
\psline[linewidth=.04](10,0)(8,-3)\pscircle*(8,-3){.2}
\rput(8.5,-1.6){\smsize{$2$}}\rput(7.5,-3.3){\smsize{$2$}}
\psline[linewidth=.04](10,0)(10,-3)\rput(10,-3.5){\smsize{$\bf 2$}}
\psline[linewidth=.04](10,0)(12,-2)\pscircle*(12,-2){.2}
\rput(11,-1.5){\smsize{$3$}}\rput(12.5,-2.2){\smsize{$1$}}
\psline[linewidth=.05,linestyle=dashed](10,0)(14,0)\pscircle*(14,0){.2}
\psline[linewidth=.04](14,0)(17,-1)\pscircle*(17,-1){.2}
\rput(16.2,-.3){\smsize{$2$}}\rput(17.5,-1){\smsize{$1$}}
\psline[linewidth=.04](14,0)(17,1)\pscircle*(17,1){.2}
\rput(15.6,.9){\smsize{$3$}}\rput(17.2,1.6){\smsize{$2$}}
\psline[linewidth=.04](17,1)(20,2)\pscircle*(20,2){.2}
\rput(18.6,1.9){\smsize{$1$}}\rput(20.2,2.6){\smsize{$1$}}
\psline[linewidth=.05,linestyle=dashed](10,0)(13,2)\pscircle*(13,2){.2}
\psline[linewidth=.04](13,2)(16,3)\pscircle*(16,3){.2}
\rput(14.6,3){\smsize{$1$}}\rput(16.2,3.6){\smsize{$1$}}
\psline[linewidth=.04](13,2)(14,4)\rput(14.3,4.4){\smsize{$\bf 3$}}
\pscircle*(30,0){.2}
\psline[linewidth=.05](30,0)(27,2)\pscircle*(27,2){.2}
\rput(28.6,1.5){\smsize{$2$}}\rput(27.3,2.5){\smsize{$1$}}
\psline[linewidth=.05](30,0)(26,0)\pscircle*(26,0){.2}
\rput(27.7,.4){\smsize{$2$}}\rput(26.3,.6){\smsize{$1$}}
\psline[linewidth=.05](26,0)(23,2)\pscircle*(23,2){.2}
\rput(24.6,1.4){\smsize{$3$}}\rput(23.4,2.5){\smsize{$2$}}
\psline[linewidth=.05](26,0)(23,-2)\pscircle*(23,-2){.2}
\rput(24.4,-.6){\smsize{$1$}}\rput(23.2,-1.3){\smsize{$3$}}
\psline[linewidth=.05](30,0)(27,-1.5)\pscircle*(27,-1.5){.2}
\rput(27.7,-.7){\smsize{$2$}}\rput(26.8,-.9){\smsize{$1$}}
\psline[linewidth=.05](27,-1.5)(25,-3)\rput(24.7,-3.1){\smsize{$\bf 1$}}
\psline[linewidth=.05](30,0)(28,-3)\pscircle*(28,-3){.2}
\rput(28.5,-1.6){\smsize{$2$}}\rput(27.5,-3.3){\smsize{$2$}}
\psline[linewidth=.05](30,0)(30,-3)\rput(30,-3.5){\smsize{$\bf 2$}}
\psline[linewidth=.05](30,0)(32,-2)\pscircle*(32,-2){.2}
\rput(31,-1.5){\smsize{$3$}}\rput(32.5,-2.2){\smsize{$1$}}
\psline[linewidth=.05](30,0)(33,-1)\pscircle*(33,-1){.2}
\rput(32.2,-.3){\smsize{$2$}}\rput(33.5,-1){\smsize{$1$}}
\psline[linewidth=.05](30,0)(33,1)\pscircle*(33,1){.2}
\rput(32,1.1){\smsize{$3$}}\rput(33.2,1.6){\smsize{$2$}}
\psline[linewidth=.05](33,1)(36,2)\pscircle*(36,2){.2}
\rput(34.6,1.9){\smsize{$1$}}\rput(36.2,2.6){\smsize{$1$}}
\psline[linewidth=.05](30,0)(32,2)\pscircle*(32,2){.2}
\rput(31,1.5){\smsize{$1$}}\rput(32.2,2.6){\smsize{$1$}}
\psline[linewidth=.05](30,0)(30.5,2.5)\rput(30.5,3){\smsize{$\bf 3$}}
\rput(27.5,4){\smsize{$(1,0)$}}
\pnode(28.3,4.2){A}\pnode(30,0){B}
\ncarc[linewidth=.05,nodesep=.5,arcangleA=30,arcangleB=10,ncurv=.7]{->}{A}{B}
\end{pspicture}
\caption{A Refined Decorated Rooted Tree and its Decorated Graph}
\label{graph_fig}
\end{figure}

\noindent
The boundary fixed loci $\wt\cZ_{\ti{\Ga}}$ 
will be indexed by {\it refined decorated rooted trees}~$\ti\Ga$.
Figure~\ref{graph_fig} shows such a tree $\ti\Ga$ and 
the corresponding decorated graph $\Ga\!=\!\pi(\ti\Ga)$.
In \cite[Section~27.3]{H} the fixed loci $\cZ_{\Ga}$ of
the $G$-action on $\ov\M_{g,k}(\Pn,d)$ are indexed by decorated graphs~$\Ga$.
If $\Ga$ is a decorated graph such that $\cZ_{\Ga}$
is a $G$-fixed locus contained in $\part\ov\M_{1,k}(\Pn,d)$, we~will have
$$\cZ_{\Ga}\cap\ov\M_{1,k}^0(\Pn,d)
=\ti\pi\Big(\bigsqcup_{\pi(\ti\Ga)=\Ga}\!\! \wt\cZ_{\ti{\Ga}}\Big),$$
where $\ti\Ga$ denotes a refined decorated rooted tree.\\

\noindent
We now formally describe what we mean by a refined decorated rooted tree
and its corresponding decorated graph.
A {\it graph} consists of a set~$\Ver$ of {\it vertices} and 
a collection $\Edg$ of {\it edges}, i.e.~of two-element subsets of~$\Ver$.
In Figure~\ref{graph_fig}, the vertices are represented by dots,
while each edge $\{v_1,v_2\}$ is shown as the line segment between $v_1$ and~$v_2$.
A graph is a {\it tree} if it contains no~{\it loops}, i.e.~the set $\Edg$ contains no subset
of the~form
$$\big\{\{v_1,v_2\},\{v_2,v_3\},\ldots,\{v_N,v_1\}\big\},
\qquad v_1,\ldots,v_N\!\in\!\Ver,~ N\!\ge\!1.$$
A tree is {\it rooted} if $\Ver$ contains a distinguished element $v_0$.
It is represented by the large dot in the first diagram of Figure~\ref{graph_fig}.
A rooted tree is {\it refined} if $\Ver\!-\!\{v_0\}$ contains two, possibly empty,
distinguished subsets $\Ver_+$ and $\Ver_0$ such~that
$$\Ver_+\!\cap\!\Ver_0=\eset \quad\hbox{and}\quad
\{v_0,v\}\!\in\!\Edg ~~\forall\, v\!\in\!\Ver_+\!\cup\!\Ver_0.$$
We put
$$\Edg_+=\big\{\{v_0,v\}\!: v\!\in\!\Ver_+\big\} \quad\hbox{and}\quad
\Edg_0=\big\{\{v_0,v\}\!: v\!\in\!\Ver_0\big\}.$$
The elements of $\Edg_+$ and $\Edg_0$ are shown in the first diagram of Figure~\ref{graph_fig}
as the thick solid lines and the thin dashed lines, respectively.
Finally, a {\it refined decorated rooted tree} is a tuple
\begin{equation}\label{treedfn_e}
\ti\Ga = \big(\Ver,\Edg;v_0;\Ver_+,\Ver_0;\mu,\d,\eta \big),
\end{equation}
where $(\Ver,\Edg;v_0;\Ver_+,\Ver_0)$ is refined rooted tree
and
$$\mu\!:\Ver\!-\!\Ver_0\lra\big\{0,\ldots,n\}, \qquad
\d\!: \Edg\!-\!\Edg_0\lra\Z^+, \quad\hbox{and}\quad
\eta\!: \{1,\ldots,k\}\lra\Ver$$
are maps such that\\
${}\quad$ (i) $\mu(v_1)\!=\!\mu(v_2)$ and $\d(\{v_0,v_1\})\!=\!\d(\{v_0,v_2\})$
for all $v_1,v_2\!\in\!\Ver_+$;\\
${}\quad$ (ii) if $v_1\!\in\!\Ver_+$, $v_2\!\in\!\Ver\!-\!\Ver_0\!-\!\Ver_+$,
and $\{v_0,v_2\}\!\in\!\Edg$, then 
$$\mu(v_1)\neq\mu(v_2) \qquad\hbox{or}\qquad
\d(\{v_0,v_1\})\!\neq\!\d(\{v_0,v_2\});$$
${}\quad$ (iii) if $\{v_1,v_2\}\!\in\!\Edg$ and $v_2\!\not\in\!\Ver_0\!\cup\!\{v_0\}$,
then 
$$\mu(v_2)\!\neq\!\mu(v_1)  \quad\hbox{if}~~~ v_1\!\not\in\!\Ver_0
\qquad\hbox{and}\qquad 
\mu(v_2)\!\neq\!\mu(v_0) \quad\hbox{if}~~~ v_1\!\in\!\Ver_0;$$
${}\quad$ (iv) if $v_1\!\in\!\Ver_0$, then $\{v_1,v_2\}\!\in\!\Edg$ for some
$v_2\!\in\!\Ver\!-\!\{v_0\}$ and
$$\val(v_1)\equiv \big|\{v_2\!\in\!\Ver\!: \{v_1,v_2\}\!\in\!\Edg\}\big|
+ \big|\{l\!\in\![k]\!: \eta(l)\!=\!v_1\big|\ge 3;$$
${}\quad$ (v) $\sum_{e\in\Edg_+}\!\!\d(e)\ge2$.\\
In Figure~\ref{graph_fig}, the value of the map $\mu$ on each vertex,
not in $\Ver_0$, is indicated by the number next to the vertex.
Similarly, the value of the map $\d$ on each edge,
not in $\Edg_0$, is indicated by the number next to the edge.
The elements of the~set $[k]\!=\![3]$ are shown in bold face.
Each of them is linked by a line segment to its image under~$\eta$.
The first condition above implies that all of the thick edges have the same labels,
and so do their vertices, other than the root~$v_0$.
By the second condition, the set of thick edges is 
a maximal set of edges leaving~$v_0$ which satisfies the first condition.
By the third condition, no two consecutive vertex labels are the same.
By the fourth condition, there are at least two solid lines, at least one of which
is an edge, leaving from every vertex which is connected to the root by a dashed line.
The final condition implies that either the set $\Edg_+$ contains at least two elements
or its only element is marked by at least~2.\\

\noindent
A {\it decorated graph} is a tuple
$$\Ga=\big(\Ver,\Edg;g,\mu,\d,\eta\big),$$
where $(\Ver,\Edg)$ is a graph and 
$$g\!:\Ver\lra\bar\Z^+, \quad \mu\!:\Ver\lra\big\{0,\ldots,n\},  \quad 
\d\!: \Edg\lra\Z^+, \quad\hbox{and}\quad \eta\!: \{1,\ldots,k\}\lra\Ver$$
are maps such that
$$\mu(v_1)\neq\mu(v_2) \qquad\hbox{if}\quad \{v_1,v_2\}\!\in\!\Edg.$$
The domain $[k]$ of the map $\eta$ can be replaced by any finite set. 
A decorated graph can be represented graphically as in the second diagram
of Figure~\ref{graph_fig}.
In this case, every vertex $v$ should be labeled by the pair $(g(v),\mu(v))$.
However, we drop the first entry if it is zero.
If $\ti\Ga$ is a refined decorated rooted tree as in~\e_ref{treedfn_e}, 
the corresponding decorated graph $\Ga$ is obtained by identifying all
elements of $\Ver_0$ with~$v_0$, dropping $\Edg_0$ from~$\Edg$,
and setting
$$ g(v) =\begin{cases} 
1,& \hbox{if}~v\!=\!v_0;\\
0,& \hbox{otherwise}. 
\end{cases}$$
In terms of the first diagram in Figure~\ref{graph_fig},
this procedure corresponds to contracting the dashed edges
and adding $1$ to the label for~$v_0$.\\

\noindent
The fixed locus $\cZ_{\Ga}$ of $\ov\M_{1,k}(\Pn,d)$
consists of the stable maps~$u$ from a genus-one nodal curve $\Si_u$
with $k$ marked points into~$\Pn$ that satisfy the following conditions.
The components of $\Si_u$ on which the map $u$ is not constant
are rational and correspond to the edges of~$\Ga$.
Furthermore, if $e\!=\!\{v_1,v_2\}$ is an edge,
the restriction of $u$ to the component $\Si_{u,e}$ corresponding to~$e$
is a degree-$\d(e)$ cover of the line 
$$\P^1_{p_{\mu(v_1)},p_{\mu(v_2)}}\subset\Pn$$
passing through the fixed points $p_{\mu(v_1)}$ and $p_{\mu(v_2)}$.
The map $u|_{\Si_{u,e}}$ is ramified only over $p_{\mu(v_1)}$ and~$p_{\mu(v_2)}$.
In particular, $u|_{\Si_{u,e}}$ is unique up to isomorphism.
The remaining, contracted, components of $\Si_u$ correspond to the vertices 
$v\!\in\!\Ver$ such that 
$$\val(v)+g(v)\ge 3.$$
For such a vertex $v$, $g(v)$ specifies the genus of the component corresponding to~$v$.
The map $u$ takes this component to the fixed point~$\mu(v)$.
Thus,
$$\cZ_{\Ga}\approx \ov\cM_{\Ga}\!\equiv\!\prod_{v\in\Ver}\!\!\ov\cM_{g(v),\val(v)};$$
see \cite[Section 27.3]{H}.
For the purposes of this definition, $\ov\cM_{0,1}$ and $\ov\cM_{0,2}$
denote one-point spaces.
For example, in the case of the second diagram in Figure~\ref{graph_fig},
$$\cZ_{\Ga}\approx
\ov\cM_{\Ga} \!\equiv \ov\cM_{1,10} \!\times\! \ov\cM_{0,3}
\!\times\! \ov\cM_{0,2}^2 \!\times\! \ov\cM_{0,1}^8 \approx \ov\cM_{1,10}.$$
In this case, $\cZ_{\Ga}$ is a locus in $\ov\M_{1,3}(\Pn,22)$,
with $n\!\ge\!3$.\\

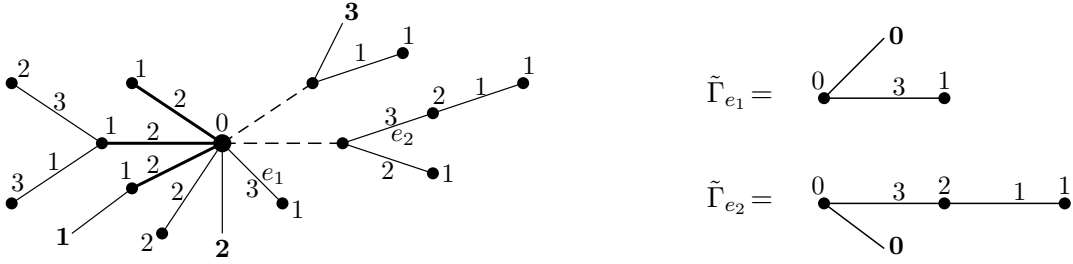
\begin{figure}
\begin{pspicture}(0,-2)(10,2)
\psset{unit=.4cm}
\pscircle*(10,0){.3}\rput(10,.7){\smsize{$0$}}
\psline[linewidth=.1](10,0)(7,2)\pscircle*(7,2){.2}
\rput(8.6,1.5){\smsize{$2$}}\rput(7.3,2.5){\smsize{$1$}}
\psline[linewidth=.1](10,0)(6,0)\pscircle*(6,0){.2}
\rput(7.7,.4){\smsize{$2$}}\rput(6.3,.6){\smsize{$1$}}
\psline[linewidth=.04](6,0)(3,2)\pscircle*(3,2){.2}
\rput(4.6,1.4){\smsize{$3$}}\rput(3.4,2.5){\smsize{$2$}}
\psline[linewidth=.04](6,0)(3,-2)\pscircle*(3,-2){.2}
\rput(4.4,-.6){\smsize{$1$}}\rput(3.2,-1.3){\smsize{$3$}}
\psline[linewidth=.1](10,0)(7,-1.5)\pscircle*(7,-1.5){.2}
\rput(7.7,-.7){\smsize{$2$}}\rput(6.8,-.9){\smsize{$1$}}
\psline[linewidth=.04](7,-1.5)(5,-3)\rput(4.7,-3.1){\smsize{$\bf 1$}}
\psline[linewidth=.04](10,0)(8,-3)\pscircle*(8,-3){.2}
\rput(8.5,-1.6){\smsize{$2$}}\rput(7.5,-3.3){\smsize{$2$}}
\psline[linewidth=.04](10,0)(10,-3)\rput(10,-3.5){\smsize{$\bf 2$}}
\psline[linewidth=.04](10,0)(12,-2)\pscircle*(12,-2){.2}
\rput(11,-1.5){\smsize{$3$}}\rput(12.5,-2.2){\smsize{$1$}}\rput(11.7,-1.1){\smsize{$e_1$}}
\psline[linewidth=.05,linestyle=dashed](10,0)(14,0)\pscircle*(14,0){.2}
\psline[linewidth=.04](14,0)(17,-1)\pscircle*(17,-1){.2}
\rput(15.5,-.9){\smsize{$2$}}\rput(17.5,-1){\smsize{$1$}}
\psline[linewidth=.04](14,0)(17,1)\pscircle*(17,1){.2}
\rput(15.6,.9){\smsize{$3$}}\rput(17.2,1.6){\smsize{$2$}}\rput(16,.2){\smsize{$e_2$}}
\psline[linewidth=.04](17,1)(20,2)\pscircle*(20,2){.2}
\rput(18.6,1.9){\smsize{$1$}}\rput(20.2,2.6){\smsize{$1$}}
\psline[linewidth=.05,linestyle=dashed](10,0)(13,2)\pscircle*(13,2){.2}
\psline[linewidth=.04](13,2)(16,3)\pscircle*(16,3){.2}
\rput(14.6,3){\smsize{$1$}}\rput(16.2,3.6){\smsize{$1$}}
\psline[linewidth=.04](13,2)(14,4)\rput(14.3,4.4){\smsize{$\bf 3$}}
\pscircle*(30,1.5){.2}\rput(29.8,2.1){\smsize{$0$}}
\psline[linewidth=.05](30,1.5)(34,1.5)\pscircle*(34,1.5){.2}
\rput(32.5,1.9){\smsize{$3$}}\rput(34,2.1){\smsize{$1$}}
\psline[linewidth=.05](30,1.5)(32,3.5)\rput(32.4,3.6){\smsize{$\bf 0$}}
\rput(27.2,1.7){$\ti\Ga_{e_1}\!=$}
\pscircle*(30,-2){.2}\rput(29.8,-1.4){\smsize{$0$}}
\psline[linewidth=.05](30,-2)(34,-2)\pscircle*(34,-2){.2}
\rput(32.5,-1.6){\smsize{$3$}}\rput(34,-1.4){\smsize{$2$}}
\psline[linewidth=.05](34,-2)(38,-2)\pscircle*(38,-2){.2}
\rput(36.5,-1.6){\smsize{$1$}}\rput(38,-1.4){\smsize{$1$}}
\psline[linewidth=.05](30,-2)(32,-3.5)\rput(32.4,-3.4){\smsize{$\bf 0$}}
\rput(27.2,-1.8){$\ti\Ga_{e_2}\!=$}
\end{pspicture}
\caption{A Refined Decorated Rooted Tree and Some of its Components Graphs}
\label{graph_fig2}
\end{figure}

\noindent
If $\ti\Ga$ is a refined decorated rooted tree as in~\e_ref{treedfn_e}, we put
$$\Edg(v_0)=\big\{\{v_0,v_1\}\!\in\!\Edg\!: v_1\!\in\!\Ver\big\}
\qquad\hbox{and}\qquad
J_{v_0}=\big\{l\!\in\![k]\!: \mu(l)\!=\!v_0\big\}.$$
Similarly, for each $v\!\in\!\Ver_0$, we set
$$\Edg(v)=\big\{\{v,v_1\}\!\in\!\Edg\!: v_1\!\in\!\Ver\!-\!\{v_0\}\big\}
\qquad\hbox{and}\qquad
J_v=\big\{l\!\in\![k]\!: \mu(l)\!=\!v\big\}.$$
If $e\!=\!\{v,v_1\}$ is an element of $\Edg(v)$ for some $v\!\in\!\Ver_0$
or of $\Edg(v_0)\!-\!\Edg_0$ with $v\!=\!v_0$,
let $(\Ver_e,\Edg_e)$ be the branch of the tree $(\Ver,\Edg)$ 
beginning at $v$ with the edge~$e$.
We~put 
$$J_e=\big\{l\!\in\![k]\!: \mu(l)\!\in\!\Ver_e\!-\!\{v\}\big\}
\qquad\hbox{and}\qquad
d_e=\sum_{e'\in\Edg_e}\!\!\!\d(e').$$
Let $\ti\Ga_e$ be the decorated graph defined by
\begin{gather*}
\ti\Ga_e=\big(\Ver_e,\Edg_e;g_e\!\equiv\!0,\mu_e,
\d_e\!\equiv\!\d|_{\Edg_e},\eta_e),
\qquad\hbox{where}\\
\mu_e(v')=\begin{cases}
\mu(v'),&\hbox{if}~v'\!\neq\!v;\\
\mu(v_0),&\hbox{if}~v'\!=\!v;
\end{cases} \qquad
\eta_e\!:\{0\}\!\sqcup\!J_e\lra\Ver_e, \quad
\eta_e(l)=\begin{cases}
\eta(l),& \hbox{if}~l\!\in\!J_e;\\
v,& \hbox{if}~l\!=\!0;
\end{cases}
\end{gather*}
see Figure~\ref{graph_fig2} for two examples.\\

\noindent
If $e$ is an element of $\Edg(v_0)\!-\!\Edg_0$ or of $\Edg(v)$
for some $v\!\in\!\Ver_0$, let
$$\cZ_{\ti\Ga_e} \subset \ov\M_{0,\{0\}\sqcup J_e}(\Pn,d_e)$$
be the fixed locus corresponding to the decorated graph $\ti\Ga_e$.
We~put 
\begin{gather*}
\si(\ti\Ga)=\big(|\Edg(v_0)|;J_{v_0},[k]\!-\!J_{v_0}\big)\in\A_1(d,k), 
\qquad
\wt\cM_{\ti\Ga;P}=\wt\cM_{\si(\ti\Ga);P};\\
\bar\cZ_{\ti\Ga;B}=\!\!
\prod_{e\in\Edg(v_0)-\Edg_0}\!\!\!\!\!\!\!\!\!\!\!\! \cZ_{\ti\Ga_e}
~\times  \prod_{v\in\Ver_0}\!\! \Big(\ov\cM_{0,\{0\}\sqcup\Edg(v)\sqcup J_v}
\times\! \prod_{e\in\Edg(v)}\!\!\!\!\!\! \cZ_{\ti\Ga_e}\Big)
\subset \ov\M_{\si(\ti\Ga);B}(\Pn,d);\\
F_{\ti\Ga;B}=\bigoplus_{e\in\Edg_+}\!\!\!L_{e;0}\subset F_{\si(\ti\Ga);B}
\lra \bar\cZ_{\ti\Ga;B},
\end{gather*}
where $L_{e;0}\!\lra\!\cZ_{\ti\Ga_e}$ is the tangent line 
bundle for the marked point~$0$.
If $e\!=\!\{v_0,v_1\}$ is an element of~$\Edg_+$, let
$$\mu_+(\ti\Ga)=\mu(v_1), \qquad \d_+(\ti\Ga)=\d(e),
\quad\hbox{and}\quad 
\dim_+(\ti\Ga)=\begin{cases}
|\Edg_+|\!-\!2,& \hbox{if}~\d_+(\ti\Ga)\!=\!1;\\
|\Edg_+|\!-\!1,& \hbox{if}~\d_+(\ti\Ga)\!\ge\!2.
\end{cases}$$
By the assumption~(i) above, the numbers $\mu_+(\ti\Ga)$ and $\d_+(\ti\Ga)$
are independent of the choice of $e\!\in\!\Edg_+$.
Furthermore, if $e,e'\!\in\!\Edg_+$, then the line bundles
$L_{e;0}$ and $L_{e';0}$ are $G$-equivariantly isomorphic.
Thus,
$$F_{\ti\Ga;B} \approx \Bbb{C}^{|\Edg_+|}\otimes L_{e;0}
\qquad\hbox{if}\quad e\!\in\!\Edg_+.$$
The group $G$ acts trivially on $\Bbb{C}^{|\Edg_+|}$.
Let
\begin{gather*}
F_{\ti\Ga;B}'= 
\begin{cases}
\big\{(w_e)_{e\in\Edg_+}\!\in\!\C^{\Edg_+}:\sum_{e\in\Edg_+}\!w_e\!=\!0\big\}, 
&\hbox{if}~ \d_+(\ti\Ga)\!=\!1;\\
\C^{\Edg_+}, &\hbox{if}~ \d_+(\ti\Ga)\!\ge\!2;
\end{cases}\\
\wt\cZ_{\ti\Ga;B}=\P\big(F_{\ti\Ga;B}'\!\otimes\!L_{e;0}\big)
\approx \bar\cZ_{\ti\Ga;B}\times\Bbb{P}^{\dim_+(\ti\Ga)}.
\end{gather*}
While the moduli space $\wt\M_{\si(\ti\Ga);B}(\Pn,d)$  is a blowup
of $\P F_{\si(\ti\Ga);B}$, none of the blowup loci intersects~$\wt\cZ_{\ti\Ga;B}$.
Thus,
$$\wt\cZ_{\ti\Ga;B} \subset \wt\M_{\si(\ti\Ga);B}(\Pn,d).$$
In fact,
$$\wt\cZ_{\ti\Ga;B} \subset \wt\cZ_{\si(\ti\Ga);B}(\Pn,d).$$
We put
$$\wt\cZ_{\ti\Ga}=\wt\cM_{\ti\Ga;P}\times \wt\cZ_{\ti\Ga;B}.$$
By the above, $\wt\cZ_{\ti\Ga}$ is a fixed point locus
in $\wt\M_{1,k}^0(\Pn,d)$.
For example, in the case of the first diagram in Figure~\ref{graph_fig2},
\begin{gather*}
\si(\ti\Ga)=\big(7;\{2\},\{1,3\}\big),  \qquad
\wt\cM_{\ti\Ga;P}=\wt\cM_{1,([7],\{2\})}, \\
\bar\cZ_{\ti\Ga;B}=\Big(\ov\cM_{0,2}^6\!\times\!\ov\cM_{0,1}^5\!\times\!\ov\cM_{0,3}\Big)
\times\Big(\ov\cM_{0,3}^2\!\times\!\ov\cM_{0,2}^4\!\times\!\ov\cM_{0,1}^3\Big)
\approx\{pt\};\\
\rk\, F_{\ti\Ga;B}=\rk\, F_{\ti\Ga;B}'=3, \qquad 
\wt\cZ_{\ti\Ga;B}\approx\P^2, \qquad 
\wt\cZ_{\ti\Ga}\approx \wt\cM_{1,([7],\{2\})}\!\times\!\P^2.
\end{gather*}
The weight of the $G$-action on the line $L_{e;0}$ is $1/2$
of the weight of the $G$-action on $T_{p_0}\P_{p_0,p_1}^1$;
see \cite[Sects~27.1 and~27.2]{H}.\\

\noindent
We next describe the equivariant normal bundle ${\cal N}\wt\cZ_{\ti\Ga}$
of $\wt\cZ_{\ti\Ga}$ in~$\wt\M_{1,k}^0(\Pn,d)$. 
Let 
$$\N_{\ov\M_{\si(\ti\Ga);B}(\Pn,d)}\bar\cZ_{\ti\Ga;B} \lra \bar\cZ_{\ti\Ga;B}$$
be the normal bundle of $\bar\cZ_{\ti\Ga;B}$ in $\ov\M_{\si(\ti\Ga);B}(\Pn,d)$.
This normal bundle can easily be described using \cite[Section~27.4]{H}.
Let
$$F_{\ti\Ga;B}^- = F_{\si(\ti\Ga);B} \big/ (F_{\ti\Ga;B}'\!\otimes\!L_{e;0})
\approx \!\!
\bigoplus_{e'\in\Edg(v_0)-\Edg^+}\!\!\!\!\!\!\!\!\!\!\!\!\!\! L_{e';0}~~
\oplus
\begin{cases}
L_{e;0},& \hbox{if}~\d_+(\ti\Ga)\!=\!1;\\
\{0\},& \hbox{if}~\d_+(\ti\Ga)\!\ge\!2,
\end{cases}$$
where $e$ is an element of $\Edg_+$.
The normal bundle of $\wt\cZ_{\ti\Ga;B}$ in $\wt\M_{\si(\ti\Ga);B}(\Pn,d)$
is given~by
$$\N_{\wt\M_{\si(\ti\Ga);B}(\Pn,d)}\wt\cZ_{\ti\Ga;B}= 
\N_{\ov\M_{\si(\ti\Ga);B}(\Pn,d)}\bar\cZ_{\ti\Ga;B}
\oplus \ga_{\dim_+\ti\Ga}^*\!\otimes\!L_{e;0}^*\!\otimes\!F_{\ti\Ga;B}^-,$$
where $\ga_{\dim_+\ti\Ga}\!\lra\!\P^{\dim_+\ti\Ga}$ is the tautological line bundle.
Since none of the exceptional divisors intersects~$\wt\cZ_{\ti\Ga;B}$,
\begin{equation}\label{bunrestr_e1}
\E_{\si;B}\big|_{\wt\cZ_{\ti\Ga;B}} = \ga_{\dim_+\ti\Ga}\!\otimes\!L_{e;0}.
\end{equation}
Since the section $\wt\cD_{\si;B}$ is transverse to the zero set, 
the normal bundle of $\wt\cZ_{\ti\Ga;B}$  in $\wt\cZ_{\si(\ti\Ga);B}(\Pn,d)$~is
$$\N_{\wt\cZ_{\si(\ti\Ga);B}(\Pn,d)}\wt\cZ_{\ti\Ga;B} =
\N_{\wt\M_{\si(\ti\Ga);B}(\Pn,d)}\wt\cZ_{\ti\Ga;B}
\big/ \big(\ga_{\dim_+\ti\Ga}^*\!\otimes\!L_{e;0}^*\!\otimes\!T_{\mu(v_0)}\Pn\big)$$
by \e_ref{spacesec_e1} and \e_ref{spacedfn_e1}.
Finally,
$${\cal N}\wt\cZ_{\ti\Ga}=
\N_{\wt\cZ_{\si(\ti\Ga);B}(\Pn,d)}\wt\cZ_{\ti\Ga;B}
\oplus \L_{\si(\ti\Ga);P}\!\otimes\!\ga_{\dim_+\ti\Ga}\!\otimes\!L_{e;0},$$
since the normal bundle of $\cZ_{\si(\ti\Ga)}(\Pn,d)$ in $\ov\M_{1,k}^0(\Pn,d)$
is $\L_{\si(\ti\Ga);P}\!\otimes\!\ga_{\si(\ti\Ga);B}$.\\

\noindent
In order to compute the last number in~\e_ref{euler_e}, we also 
need to determine the restriction of the vector bundle $\wt\V_{1,k}^d$
to~$\wt\cZ_{\ti\Ga}$.
By Theorem~\ref{cone_thm} and~\e_ref{bunrestr_e1},
there is a short exact sequence of vector bundles:
$$0\lra \wt\V_{1,k}^d\big|_{\wt\cZ_{\ti\Ga}} \lra
\V_{\si(\ti\Ga);B}^d\big|_{\wt\cZ_{\ti\Ga;B}} \lra
\ga_{\dim_+\ti\Ga}^*\!\otimes\!L_{e;0}^*\!\otimes\!\cL_{\mu(v_0)} \lra 0$$
over $\wt\cZ_{\ti\Ga}$. 
This exact sequence describes the euler class of 
the restriction of $\wt\V_{1,k}^d$ to~$\wt\cZ_{\ti\Ga}$.\\

\noindent
If $\si\!=\!(m;J_P,J_B)\!\in\!\A_1(d,k)$, 
$$ \blr{c_1^{|m|+|J_P|}(\L_{\si;P}^*),\wt\cM_{\si;P}}
=\frac{m^{|J_P|}\cdot(m\!-\!1)!}{24},$$
by \cite{g1desing2}.
This is the only intersection number on $\wt\cM_{\si;P}$ needed
for computing the last number in~\e_ref{euler_e} and 
the integrals of the cohomology classes on $\ov\M_{1,k}^0(\Pn,d)$
that count elliptic curves in $\Pn$ passing through specified constraints.
For more general enumerative problems, such as counting curves with tangency
conditions, as in~\cite{enumtang}, and 
with singularities, as in~\cite{genuss0pr}, one would need to compute
the intersection numbers of the form
$$\Blr{c_1^{\be_0}(\L_{\si;P}^*)\cdot\prod_{l\in J_P}\!\psi_l^{\be_l},\wt\cM_{\si;P}},
\qquad\hbox{where}\qquad  \be_0+\sum_{l\in J_P}\be_l=|m|+|J_P|.$$
The argument in~\cite{g1desing2} gives a recursive formula for a generalization
of such numbers. 
The recursion is on $|m|\!+\!|J_P|$, i.e.~the total number of marked points.
The starting data for the recursion are the numbers
$$\Blr{\prod_{l=1}^{l=k}\!\psi_l^{\be_l},\ov\cM_{1,k}},
\qquad\hbox{where}\qquad  \sum_{l=1}^{l=k}\be_l=k,$$
which in turn are computable from the genus-one string and dilaton equations;
see~\cite[Section~25.2]{H}.\\

\noindent
In the above discussion we ignored all automorphism groups.
As in~\cite[Chapter~27]{H}, the rational function for 
each refined decorated rooted tree~$\ti\Ga$
obtained following the above algorithm and applying the localization theorem of~\cite{ABo} 
should be divided by the order of the appropriate automorphism group~$\bA_{\ti\Ga}:$
$$\big|\bA_{\ti\Ga}\big|=\big|\Aut(\ti\Ga)\big|\cdot
\prod_{e\in\Edg-\Edg_0}\!\!\!\!\!\!\! \d(e).$$
For example, in the case of the first diagram in Figure~\ref{graph_fig2},
$$\big|\bA_{\ti\Ga}\big|= 1\cdot \big(1^3\cdot2^5\cdot3^3)=864.$$

\section{Blowups of Moduli Spaces of Curves}
\label{curvebl_sec}

\subsection{Blowups and Subvarieties}
\label{curveprelim_subs}

\noindent
In this section we construct blowups of certain moduli spaces of genus-one 
and genus-zero curves; see Subsections~\ref{curve1bl_subs} and~\ref{curve0bl_subs}.
The former appear in Subsection~\ref{map1blconstr_subs} as the first factor in 
the domain of the proper transforms of the immersion~$\io_{\si}$ of
Subsection~\ref{descr_subs}.
The latter play the analogous role in Subsection~\ref{map0blconstr_subs},
where we blow up certain moduli spaces of genus-zero maps.
In turn, these last blowups describe the second factor of the domain of maps
induced by~$\io_{\si}$ in Subsection~\ref{map1blconstr_subs};
see Subsection~\ref{outline_subs} for more details.\\

\noindent
We begin by introducing convenient terminology and reviewing standard facts 
from algebraic geometry.
If $\ov\cM$ is a smooth variety and $Z$ is a smooth subvariety of $\ov\cM$,
let
$$\N_{\ov\cM}Z\equiv T\ov\cM|_Z\big/TZ$$
be the normal bundle of $Z$ in $\ov\cM$.
We denote~by
$$\pi_Z^{\perp}\!: T\ov\cM|_Z \lra \N_{\ov\cM}Z$$
the quotient projection map.

\begin{dfn}
\label{ag_dfn1}
Let $\ov\cM$ be a smooth variety.\\
(1) Smooth subvarieties $X$ and $Y$ of $\ov\cM$ {\tt intersect properly} if 
$X\!\cap\!Y$ is a smooth subvariety of $\ov\cM$ and
$$T(X\!\cap\!Y)=TX|_{X\cap Y}\cap TY|_{X\cap Y}.
\footnote{In other words, the scheme-theoretic intersection of $X$ and $Y$ is smooth. 
If the set-theoretic intersection $X\!\cap\!Y$ is smooth, 
the second part of this condition is also equivalent to the injectivity of 
the natural homomorphism $$TX|_{X\cap Y}/T(X\!\cap\!Y)\lra T\ov\cM/TY.$$ }$$
(2) If $Z$ is a smooth subvariety of $\ov\cM$,
properly intersecting subvarieties $X$ and $Y$ of~$\ov\cM$
{\tt intersect properly relative to~$Z$} if
$$\pi_Z^{\perp}\big( T(X\!\cap\!Y)|_{X\cap Y\cap Z}\big)
=\pi_Z^{\perp}\big( TX|_{X\cap Y\cap Z}\big)\cap \pi_Z^{\perp}\big( TY|_{X\cap Y\cap Z}\big)
\subset \N_{\ov\cM}Z.$$\\
\end{dfn}

\noindent
For example, if $X$ and $Y$ are two smooth curves in a projective space 
that intersect without being tangent to each other, then $X$ and $Y$ intersect properly
(but not transversally, unless the dimension of the projective space is~$2$).
If $X$, $Y$, and $Z$ are three distinct lines that lie a plane,
then they intersect properly pairwise, but $X$ and $Y$ do not intersect properly 
relative to~$Z$.

\begin{dfn}
\label{subvarcoll_dfn}
If $\ov\cM$ is a smooth variety, a collection $\{\ov\cM_{\rho}\}_{\rho\in\A}$ 
of smooth subvarieties is {\tt properly intersecting} if 
$\ov\cM_{\rho_1}$ and $\ov\cM_{\rho_2}$ intersect properly relative to~$\ov\cM_{\rho_3}$
for all $\rho_1,\rho_2,\rho_3\!\in\!\A$. 
\end{dfn}

\noindent
If $Z$ is a smooth subvariety of $\ov\cM$, let
$$\pi\!:\Bl_Z\ov\cM\lra\ov\cM$$
be the blowup of $\ov\cM$ along $Z$.
If $X$ is a subvariety of $\ov\cM$, we denote by $\Pr_ZX$ the proper transform of
$X$ in $\Bl_Z\ov\cM$, i.e.~the closure of $\pi^{-1}(X\!-\!Z)$ in~$\Bl_Z\ov\cM$. 

\begin{lmm}
\label{ag_lmm1}
Let $\ov\cM$ be a smooth variety.\\
(1) If $X$ and $Z$ are properly intersecting subvarieties of $\ov\cM$, 
then $\Pr_ZX$ is a smooth subvariety of $\Bl_Z\ov\cM$ and
$$\Pr_ZX=\Bl_{X\cap Z}X.$$
(2) If $X$, $Y$, and $Z$ are pairwise properly intersecting subvarieties of $\ov\cM$
and $X$ and $Y$ intersect properly relative to~$Z$, then
$\Pr_ZX$ and $\Pr_ZY$ are properly intersecting subvarieties of $\Pr_Z\ov\cM$ and 
$$\Pr_ZX \cap \Pr_ZY = \Pr_Z(X\!\cap\!Y).$$
(3) If $X$, $Y$, $Z$, and $Z'$ are pairwise properly intersecting subvarieties of $\ov\cM$
and $X$ and $Y$ intersect properly relative to~$Z$ and~$Z'$, then
$\Pr_ZX$ and $\Pr_ZY$ intersect properly relative to~$\Pr_ZZ'$.
\end{lmm}

\begin{crl}
\label{subvarcoll_crl}
If $\ov\cM$ is a smooth variety, $\{\ov\cM_{\rho}\}_{\rho\in\A}$
is a properly intersecting collection of subvarieties of $\ov\cM$, and $\rho\!\in\!\A$, 
then $\{\Pr_{\ov\cM_{\rho}}\ov\cM_{\rho'}\}_{\rho'\in\A-\{\rho\}}$
is a properly intersecting collection of subvarieties of $\Bl_{\ov\cM_{\rho}}\ov\cM$.
\end{crl}

\subsection{Moduli Spaces of Genus-One and Zero Curves}
\label{curvebldata_subs}

\noindent
In this subsection, we describe natural subvarieties of moduli spaces of 
genus-one and -zero curves and natural bundle sections over these moduli spaces.
These bundle sections and their twisted versions introduced in the next two subsections 
are used in Subsections~\ref{map0blconstr_subs} and~\ref{map1blconstr_subs} to describe 
the structure of the proper transforms of $\ov\M_{1,k}^0(\Pn,d)$.
Below we also state the now-standard facts about these objects 
that are used in the next two subsections.\\

\noindent
If $I$ is a finite set, let
\begin{equation}\label{g0and1curvsubv_e}\begin{split}
\A_1(I) &=\big\{\big(I_P,\{I_k\!:k\!\in\!K\}\big)\!: 
K\!\neq\!\eset;~I\!=\!\bigsqcup_{k\in\{P\}\sqcup K}\!\!\!\!\!I_k;~
|I_k|\!\ge\!2 ~\forall\, k\!\in\!K\big\};\\
\A_0(I) &=\big\{\big(I_P,\{I_k\!:k\!\in\!K\}\big)\!: 
K\!\neq\!\eset;~I\!=\!\bigsqcup_{k\in\{P\}\sqcup K}\!\!\!\!\!I_k;~
|I_k|\!\ge\!2 ~\forall\, k\!\in\!K;~|K|\!+\!|I_P|\!\ge\!2\big\}.
\end{split}\end{equation}
If $\rho\!=\!(I_P,\{I_k\!:k\!\in\!K\})$ is an element of $\{(I,\eset)\}\!\sqcup\!\A_1(I)$, 
we denote by $\cM_{1,\rho}$ the subset of $\ov\cM_{1,I}$
consisting of the stable curves~$\cC$ such~that\\
${}\quad$ (i) $\cC$ is a union of a smooth torus and $|K|$ projective lines,
indexed by~$K$;\\
${}\quad$ (ii) each line is attached directly to the torus;\\
${}\quad$ (iii) for each $k\!\in\!K$,
the marked points on the line corresponding to $k$ are indexed by $I_k$.\\
Let $\ov\cM_{1,\rho}$ be the closure of $\cM_{1,\rho}$ in~$\ov\cM_{1,I}$.
Figure~\ref{g1curv_fig} illustrates this definition,
from the points of view of symplectic topology and of algebraic geometry.
In the first diagram, each circle represents a sphere, or~$\P^1$.
In the second diagram, the irreducible components of $\cC$ are represented by curves,
and the integer next to each component shows its genus. 
Similarly, if 
$$\rho\!=\!(I_P,\{I_k\!:k\!\in\!K\})\in
\big\{(I,\eset)\big\}\!\sqcup\!\A_0(I),$$
let $\cM_{0,\rho}$ be the subset of $\ov\cM_{0,\{0\}\sqcup I}$
consisting of the stable curves~$\cC$ such~that\\
${}\quad$ (i) the components of $\cC$ are indexed by the set $\{P\}\!\sqcup\!K$;\\
${}\quad$ (ii) for each $k\!\in\!K$, the component $\cC_k$ of $\cC$
is attached directly to~$\cC_P$;\\
${}\quad$ (iii) for each $k\!\in\!K$,
the marked points on $\cC_k$ are indexed by $I_k$.\\
We denote by $\ov\cM_{0,\rho}$ the closure of $\cM_{0,\rho}$ in $\ov\cM_{0,\{0\}\sqcup I}$.
This definition is illustrated in Figure~\ref{g0curv_fig}.
In this case, we do not indicate the genus of the irreducible components
in the second diagram, as all of the curves are rational.

\begin{figure}
\begin{pspicture}(-1.1,-1.8)(10,1.3)
\psset{unit=.4cm}
\rput{45}(0,-4){\psellipse(5,-1.5)(2.5,1.5)\pscircle*(7.5,-1.5){.2}\pscircle*(2.5,-1.5){.2}
\psarc[linewidth=.05](5,-3.3){2}{60}{120}\psarc[linewidth=.05](5,0.3){2}{240}{300}
\pscircle(5,-4){1}\pscircle*(5,-3){.15}
\pscircle*(5,-5){.2}\pscircle*(4,-4){.2}\pscircle*(6,-4){.2}
\pscircle(6.83,.65){1}\pscircle*(6.44,-.28){.15}
\pscircle*(5.9,1.04){.2}\pscircle*(7.76,.26){.2}
\pscircle(3.17,.65){1}\pscircle*(3.56,-.28){.15}
\pscircle*(4.1,1.04){.2}\pscircle*(2.24,.26){.2}}
\rput(7,.4){\smsize{$i_1$}}\rput(2.4,-3.7){\smsize{$i_2$}}
\rput(1.1,-2.7){\smsize{$i_3$}}\rput(2,.3){\smsize{$i_4$}}
\rput(3.1,.5){\smsize{$i_5$}}\rput(6,1.8){\smsize{$i_6$}}
\rput(7.7,-2.5){\smsize{$i_7$}}\rput(7.7,-4.2){\smsize{$i_8$}}
\rput(5.3,-4.5){\smsize{$i_9$}}
\psarc(15,-1){3}{-60}{60}\pscircle*(16.93,1.3){.2}\pscircle*(16.93,-3.3){.2}
\rput(17.6,1.4){\smsize{$i_1$}}
\rput(17.6,-3.4){\smsize{$i_2$}}
\psline(16.8,0)(22.05,1.25)\pscircle*(18.9,.5){.2}\pscircle*(21,1){.2}
\rput(18.9,1.2){\smsize{$i_3$}}\rput(21,1.7){\smsize{$i_4$}}
\psline(17,-1)(22,-1)\pscircle*(19.5,-1){.2}\pscircle*(21,-1){.2}
\rput(19.5,-.3){\smsize{$i_5$}}\rput(21,-.3){\smsize{$i_6$}}
\psline(16.8,-2)(22.05,-3.25)\pscircle*(18.9,-2.5){.2}
\pscircle*(19.95,-2.75){.2}\pscircle*(21,-3){.2}
\rput(18.9,-1.8){\smsize{$i_7$}}\rput(19.95,-2.05){\smsize{$i_8$}}
\rput(21,-2.3){\smsize{$i_9$}}
\rput(16.1,-3.5){$1$}\rput(22.5,1.2){$0$}\rput(22.4,-1){$0$}\rput(22.5,-3.2){$0$}
\rput(32,-1){\begin{small}
\begin{tabular}{l}
$I_P\!=\!\{i_1,i_2\}$\\
$K\!=\!\{1,2,3\}$\\
$I_1\!=\!\{i_3,i_4\}$\\
$I_2\!=\!\{i_5,i_6\}$\\
$I_3\!=\!\{i_7,i_8,i_9\}$\\
\end{tabular}\end{small}}
\end{pspicture}
\caption{A Typical Element of $\ov\cM_{1,\rho}$}
\label{g1curv_fig}
\end{figure}

\begin{figure}
\begin{pspicture}(-1.1,-1.8)(10,1.3)
\psset{unit=.4cm}
\pscircle(5,-1.5){1.5}
\pscircle*(5,-3){.2}\rput(4.9,-3.6){\smsize{$0$}}
\pscircle*(6.06,-2.56){.2}\rput(6.6,-2.9){\smsize{$i_1$}}
\pscircle(2.5,-1.5){1}\pscircle*(3.5,-1.5){.15}
\pscircle*(2.5,-.5){.2}\rput(2.9,0){\smsize{$i_2$}}
\pscircle*(2.5,-2.5){.2}\rput(2.3,-3){\smsize{$i_3$}}
\pscircle(5,1){1}\pscircle*(5,0){.15}
\pscircle*(4,1){.2}\rput(3.6,.7){\smsize{$i_4$}}
\pscircle*(6,1){.2}\rput(6.6,1){\smsize{$i_5$}}
\pscircle(7.5,-1.5){1}\pscircle*(6.5,-1.5){.15}
\pscircle*(7.5,-.5){.2}\rput(7.9,0){\smsize{$i_6$}}
\pscircle*(8.5,-1.5){.2}\rput(9.1,-1.5){\smsize{$i_7$}}
\pscircle*(7.5,-2.5){.2}\rput(8,-2.9){\smsize{$i_8$}}
\psline(18,2)(18,-4)\pscircle*(18,1.5){.2}\rput(17.5,1.3){\smsize{$i_1$}}
\pscircle*(18,-3.5){.2}\rput(17.5,-3.6){\smsize{$0$}}
\psline(16.8,0)(22.05,1.25)\pscircle*(18.9,.5){.2}\pscircle*(21,1){.2}
\rput(18.9,1.2){\smsize{$i_2$}}\rput(21,1.7){\smsize{$i_3$}}
\psline(17,-1)(22,-1)\pscircle*(19.5,-1){.2}\pscircle*(21,-1){.2}
\rput(19.5,-.3){\smsize{$i_4$}}\rput(21,-.3){\smsize{$i_5$}}
\psline(16.8,-2)(22.05,-3.25)\pscircle*(18.9,-2.5){.2}
\pscircle*(19.95,-2.75){.2}\pscircle*(21,-3){.2}
\rput(18.9,-1.8){\smsize{$i_6$}}\rput(19.95,-2.05){\smsize{$i_7$}}
\rput(21,-2.3){\smsize{$i_8$}}
\rput(32,-1){\smsize{\begin{tabular}{l}
$I_P\!=\!\{i_1\}$\\
$K\!=\!\{1,2,3\}$\\
$I_1\!=\!\{i_2,i_3\}$\\
$I_2\!=\!\{i_4,i_5\}$\\
$I_3\!=\!\{i_6,i_7,i_8\}$\\
\end{tabular}}}
\end{pspicture}
\caption{A Typical Element of $\ov\cM_{0,\rho}$}
\label{g0curv_fig}
\end{figure}

\begin{lmm}
\label{curvstr_lmm}
If $g\!=\!0,1$ and $I$ is a finite set, the collection
$\{\ov\cM_{g,\rho}\}_{\rho\in\A_g(I)}$ is properly intersecting.
\end{lmm}

\noindent
We define a partial ordering on the sets $\A_g(I)$ for $g\!=\!0,1$ by setting
\begin{equation}\label{partorder_e}
\rho'\!\equiv\!\big(I_P',\{I_k'\!: k\!\in\!K'\}\big)
\prec \rho\!\equiv\!\big(I_P,\{I_k\!: k\!\in\!K\}\big)
\end{equation}
if $\rho'\!\neq\!\rho$ and there exists a map $\vph\!:K\!\lra\!K'$
such that $I_k\!\subset\!I_{\vph(k)}'$ for all $k\!\in\!K$.
This condition means that the elements of $\cM_{\rho'}$ can be obtained
from the elements of $\cM_{\rho}$ by moving more points onto the bubble components
or combining the bubble components; see Figure~\ref{partorder_fig}.
In the $g\!=\!0$ case, we define the bubble components to be the components not
containing the marked point~$0$.\\

\begin{figure}
\begin{pspicture}(-1.1,-1.8)(10,1.3)
\psset{unit=.4cm}
\psarc(-2,-1){3}{-60}{60}\pscircle*(-.07,1.3){.2}\pscircle*(-.07,-3.3){.2}
\rput(.6,1.4){\begin{small}$i_1$\end{small}}
\rput(.6,-3.4){\begin{small}$i_2$\end{small}}
\psline(-.2,0)(6.05,1.25)\pscircle*(1.8,.4){.2}\pscircle*(3.05,.65){.2}
\pscircle*(4.3,.9){.2}\pscircle*(5.55,1.15){.2}
\rput(1.8,1.1){\begin{small}$i_3$\end{small}}
\rput(2.9,0){\begin{small}$i_4$\end{small}}
\rput(4.2,1.6){\begin{small}$i_5$\end{small}}
\rput(5.5,.5){\begin{small}$i_6$\end{small}}
\psline(-.2,-2)(5.05,-3.25)\pscircle*(1.9,-2.5){.2}
\pscircle*(2.95,-2.75){.2}\pscircle*(4,-3){.2}
\rput(1.9,-1.8){\begin{small}$i_7$\end{small}}
\rput(2.95,-2.05){\begin{small}$i_8$\end{small}}
\rput(4,-2.3){\begin{small}$i_9$\end{small}}
\rput(-.9,-3.5){$1$}\rput(6.45,1.3){$0$}\rput(5.5,-3.2){$0$}
\rput(8,-1){\begin{Huge}$\prec$\end{Huge}}
\psarc(9,-1){3}{-60}{60}\pscircle*(10.93,1.3){.2}\pscircle*(10.93,-3.3){.2}
\rput(11.6,1.4){\begin{small}$i_1$\end{small}}
\rput(11.6,-3.4){\begin{small}$i_2$\end{small}}
\psline(10.8,0)(17.05,1.25)\pscircle*(14.05,.65){.2}
\pscircle*(15.3,.9){.2}\pscircle*(16.55,1.15){.2}
\pscircle*(12,-1){.2}\rput(12.7,-.8){\begin{small}$i_3$\end{small}}
\rput(13.9,0){\begin{small}$i_4$\end{small}}
\rput(15.2,1.6){\begin{small}$i_5$\end{small}}
\rput(16.5,.5){\begin{small}$i_6$\end{small}}
\psline(10.8,-2)(16.05,-3.25)\pscircle*(12.9,-2.5){.2}
\pscircle*(13.95,-2.75){.2}\pscircle*(15,-3){.2}
\rput(12.9,-1.8){\begin{small}$i_7$\end{small}}
\rput(13.95,-2.05){\begin{small}$i_8$\end{small}}
\rput(15,-2.3){\begin{small}$i_9$\end{small}}
\rput(10.1,-3.5){$1$}\rput(17.45,1.3){$0$}\rput(16.5,-3.2){$0$}
\psarc(19,-1){3}{-60}{60}\pscircle*(20.93,1.3){.2}\pscircle*(20.93,-3.3){.2}
\rput(21.6,1.4){\begin{small}$i_1$\end{small}}
\rput(21.6,-3.4){\begin{small}$i_2$\end{small}}
\psline(20.8,0)(27.05,1.25)\pscircle*(22.8,.4){.2}\pscircle*(24.05,.65){.2}
\pscircle*(25.3,.9){.2}\pscircle*(26.55,1.15){.2}
\rput(22.8,1.1){\begin{small}$i_3$\end{small}}
\rput(23.9,0){\begin{small}$i_4$\end{small}}
\rput(25.2,1.6){\begin{small}$i_5$\end{small}}
\rput(26.5,.5){\begin{small}$i_6$\end{small}}
\pscircle*(21.95,-.48){.2}\rput(21.4,-.7){\begin{small}$i_7$\end{small}}
\pscircle*(21.95,-1.52){.2}\rput(22.6,-1.6){\begin{small}$i_8$\end{small}}
\pscircle*(21.6,-2.5){.2}\rput(21,-2.4){\begin{small}$i_9$\end{small}}
\rput(20.1,-3.5){$1$}\rput(27.45,1.3){$0$}
\psarc(29,-1){3}{-60}{60}\pscircle*(30.93,1.3){.2}\pscircle*(30.93,-3.3){.2}
\rput(31.6,1.4){\begin{small}$i_1$\end{small}}
\rput(31.6,-3.4){\begin{small}$i_2$\end{small}}
\psline(30.8,0)(36.05,1.25)\pscircle*(32.9,.5){.2}\pscircle*(35,1){.2}
\rput(32.9,1.2){\begin{small}$i_3$\end{small}}
\rput(35,1.7){\begin{small}$i_4$\end{small}}
\psline(31,-1)(36,-1)\pscircle*(33.5,-1){.2}\pscircle*(35,-1){.2}
\rput(33.5,-.3){\begin{small}$i_5$\end{small}}
\rput(35,-.3){\begin{small}$i_6$\end{small}}
\psline(30.8,-2)(36.05,-3.25)\pscircle*(32.9,-2.5){.2}
\pscircle*(33.95,-2.75){.2}\pscircle*(35,-3){.2}
\rput(32.9,-1.8){\begin{small}$i_7$\end{small}}
\rput(33.95,-2.05){\begin{small}$i_8$\end{small}}
\rput(35,-2.3){\begin{small}$i_9$\end{small}}
\rput(30.1,-3.5){$1$}\rput(36.5,1.2){$0$}\rput(36.4,-1){$0$}\rput(36.5,-3.2){$0$}
\end{pspicture}
\caption{Examples of Partial Ordering~\e_ref{partorder_e}}
\label{partorder_fig}
\end{figure}

\noindent
In the blowup constructions of the next two subsections we will twist certain
line bundles over moduli spaces of curves and homomorphisms between them.
In the rest of this subsection we describe the relevant starting data.\\

\noindent
For each $i\!\in\!I$, let $L_i\!\lra\!\ov\cM_{1,I}$ be the universal tangent line bundle
at the marked point labeled~$i$.
Let $\E\!\lra\!\ov\cM_{1,I}$ be the Hodge line bundle of holomorphic differentials.
The natural pairing of tangent vectors with cotangent vectors induces a section
$$s_i\in\Ga\big(\ov\cM_{1,I};\Hom(L_i,\E^*)\big).$$
Explicitly,
\begin{gather*}
\big\{s_i([\cC;w])\big\}([\cC,\psi])=\psi_{x_i(\cC)}w \qquad\hbox{if}\\
[\cC]\!\in\!\ov\cM_{1,I}, \quad [\cC,w]\!\in\!L_i|_{\cC}\!=\!T_{x_i(\cC)}\cC, 
\quad [\cC,\psi]\!\in\!\E|_{\cC}\!=\!H^0(\cC;T^*\cC),
\end{gather*}
and $x_i(\cC)\!\in\!\cC$ is the marked point on $\cC$ labeled by~$i$.\\

\noindent
In the genus-zero case, the line bundle $L_0\!\lra\!\ov\cM_{0,\{0\}\sqcup I}$
will be one of the substitutes for~$\E$.
We note that for every $p\!\in\!\P^1$, there is a natural isomorphism 
between the tangent space $T_p\P^1$ of $\P^1$ at $p$ and 
the space of holomorphic differentials $H^0(\P^1;T^*\P^1\!\otimes\!\O(2p))$
on $\P^1$ that have a pole of order two at~$p$.
More precisely, let $w$ be a meromorphic function on $\P^1$
such that $p$ is the only zero of $w$ and this zero is a simple one.
We can then view $w$ as a coordinate around $p$ in~$\P^1$.
Every tangent vector $v\!\in\!T_p\P^1$ can be written~as
$$v=c_w(v)\frac{\part}{\part w}, \qquad c_w(v)\in\C.$$
We define the isomorphism
$$\psi\!: T_p\P^1 \lra H^0(\P^1;T^*\P^1\!\otimes\!\O(2p)) \qquad\hbox{by}\qquad
v \lra \psi_v=\frac{c_w(v)\, dw}{w^2}.$$ 
If $w'$ is another  meromorphic function on $\P^1$
such that $p$ is the only zero of $w'$ and this zero is a simple one, then
$$w'=\frac{w}{\al w\!+\!\be} \quad\Lra\quad
dw'=\frac{\be\, dw}{(\al w\!+\!\be)^2}  , ~~
c_{w'}(v)= \frac{c_w(v)}{\be} 
\quad\Lra\quad
\frac{c_{w'}(v)\, dw'}{{w'}^2} = \frac{c_w(v)\, dw}{w^2}. $$
Thus, the isomorphism $\psi$ is well-defined.
If $i\!\in\!I$, we define the section
\begin{gather*}
s_i \in \Ga\big(\ov\cM_{0,\{0\}\sqcup I};\Hom(L_i,L_0^*)\big)
\quad\hbox{by}\quad  
\big\{s_i([\cC;w])\big\}([\cC,v])=\psi_v\big|_{x_i(\cC)}w\\
\hbox{if}\qquad
[\cC]\!\in\!\ov\cM_{0,\{0\}\sqcup I}, \quad [\cC,w]\!\in\!L_i|_{\cC}\!=\!T_{x_i(\cC)}\cC, 
\quad [\cC,v]\!\in\!L_0|_{\cC}\!=\!T_{x_0(\cC)}\cC.
\end{gather*}\\

\noindent
We note that in both cases the section $s_i$ vanishes precisely on the curves for
which the point $i$ lies on a bubble component.
In fact,
\begin{equation}\label{curvezero_e}
s_i^{-1}(0)=\!\sum_{\rho\in\B_g(I;i)}\!\!\!\!\!\ov\cM_{1,\rho},
\qquad\hbox{where}\qquad
\B_g(I;i)=\big\{ \big(I_P,\{I_B\}\big)\!\in\!\A_g(I) \!: i\!\in\!I_B\big\}.
\end{equation}

\subsection{A Blowup of a Moduli Space of Genus-One Curves}
\label{curve1bl_subs}

\noindent
Let $I$ and $J$ be finite sets such that $I$ is nonempty.
In this subsection, we construct a blowup 
$$\pi_{1,(I,J)}\!:\wt\cM_{1,(I,J)} \lra \ov\cM_{1,I\sqcup J}$$
of the moduli space $\ov\cM_{1,I\sqcup J}$,
$|I|\!+\!1$ line bundles 
$$\ti\E,\, \ti{L}_i \lra \wt\cM_{1,(I,J)}, \qquad i\!\in\!I,$$
and $|I|$ nowhere vanishing sections 
$$\ti{s}_i\in\Ga\big(\wt\cM_{1,(I,J)};\Hom(\ti{L}_i,\ti\E^*)\big),
\qquad i\!\in\!I.$$
Since the sections $\ti{s}_i$ do not vanish, all $|I|\!+\!1$ bundles 
$\ti{L}_i$ and $\ti\E^*$ are  explicitly isomorphic.
They will be denoted by~$\L$ and called the universal tangent line bundle.\\

\noindent
The smooth variety $\wt\cM_{1,(I,J)}$ is obtained by blowing up some of
the subvarieties $\ov\cM_{1,\rho}$, defined in the previous subsection, and 
their proper transforms in an order consistent with the partial ordering~$\prec$.
The line bundle $\ti\E$ is the sum of the Hodge line bundle $\E$
and all exceptional divisors.
For each given $i\!\in\!I$, $\ti{L}_i$ is the tangent line bundle $L_i$ 
for the marked point~$i$ minus some of these divisors.
The section~$\ti{s}_i$ is induced from the pairing~$s_i$ of the previous subsection.\\

\noindent
With $I$ and $J$ as above and $\A_g(I\!\sqcup\!J)$ as in~\e_ref{g0and1curvsubv_e}, let
\begin{equation}\label{smallcolldfn_e}
\A_g(I,J) =\big\{\big((I_P\!\sqcup\!J_P),\{I_k\!\sqcup\!J_k\!: k\!\in\!K\}\big)\!\in\!\A_g(I\!\sqcup\!J)\!:
~I_k\!\neq\!\eset~ \forall\, k\!\in\!K \big\}.
\end{equation}
We note that if $\rho\!\in\!\A_g(I\!\sqcup\!J)$, then $\rho\!\in\!\A_g(I,J)$ 
if and only if every bubble component of an element of $\cM_{\rho}$ 
carries at least one element of~$I$.
Furthermore,
\begin{equation}\label{curve1zero_e}
\B_g(I\!\sqcup\!J;i)\subset \A_g(I,J) \qquad\forall\,i\!\in\!I.
\end{equation}
If $|I|\!+\!|J|\!\ge\!2$, with respect to the partial ordering $\prec$
the set $\A_1(I,J)$ has a unique minimal element:
$$\rho_{\min}\equiv\big(\eset,\{I\!\sqcup\!J\}\big).$$
Let $<$ be an ordering on $\A_1(I,J)$ extending the partial ordering $\prec$.
We denote the corresponding maximal element by~$\rho_{\max}$.
If $\rho\!\in\!\A_1(I,J)$, we~put
\begin{equation}\label{minusdfn_e}
\rho\!-\!1=
\begin{cases}
\max\{\rho'\!\in\!\A_1(I,J)\!: \rho'\!<\!\rho\},&
\hbox{if}~ \rho\!\neq\!\rho_{\min};\\
0,& \hbox{if}~\rho\!=\!\rho_{\min},
\end{cases}
\end{equation}
where the maximum is taken with respect to the ordering $<$.\\

\noindent
We now describe the starting data for the inductive blowup procedure
involved in constructing the space $\wt\cM_{1,(I,J)}$ and the line bundle
$\L$ over~$\wt\cM_{1,(I,J)}$.
Let
$$\ov\cM_{1,(I,J)}^0=\ov\cM_{1,I\sqcup J}, \qquad
\E_0=\E\lra\ov\cM_{1,(I,J)}^0, \quad\hbox{and}\quad
\ov\cM_{1,\rho}^0=\ov\cM_{1,\rho} ~~~\forall \,\rho\!\in\!\A_1(I,J).$$
For each $i\!\in\!I$, let
$$L_{0,i}=L_i\lra\ov\cM_{1,(I,J)}^0 \qquad\hbox{and}\qquad
s_{0,i}=s_i\in\Ga\big(\ov\cM_{1,(I,J)}^0;\Hom(L_{0,i},\E_0^*)\big).$$
By~\e_ref{curvezero_e},
$$s_{0,i}^{\,-1}(0)=\!\sum_{\rho^*\in\B_1(I\sqcup J;i)}\!\!\!\!\!\!\!\! \ov\cM_{1,\rho^*}^0.$$
\\

\noindent
Suppose $\rho\!\in\!\A_1(I,J)$ and we have constructed\\
${}\quad$ ($I1$) a blowup $\pi_{\rho-1}\!:\ov\cM_{1,(I,J)}^{\rho-1}\!\lra\!\ov\cM_{1,(I,J)}^0$
of $\ov\cM_{1,(I,J)}^0$ such that $\pi_{\rho-1}$ is an isomorphism 
outside of the preimages of 
the spaces $\ov\cM_{1,\rho'}^0$ with $\rho'\!\le\!\rho\!-\!1$;\\
${}\quad$ ($I2$) line bundles $L_{\rho-1,i}\!\lra\!\ov\cM_{1,(I,J)}^{\rho-1}$ 
for $i\!\in\!I$ and $\E_{\rho-1}\!\lra\!\ov\cM_{1,(I,J)}^{\rho-1}$;\\
${}\quad$ ($I3$) sections $s_{\rho-1,i}\!\in\!
\Ga(\ov\cM_{1,(I,J)}^{\rho-1};\Hom(L_{\rho-1,i},\E_{\rho-1}^*))$ for $i\!\in\!I$.\\
For each $\rho^*\!>\!\rho\!-\!1$, let $\ov\cM_{1,\rho^*}^{\rho-1}$ 
be the proper transform of~$\ov\cM_{1,\rho^*}^0$ in $\ov\cM_{1,(I,J)}^{\rho-1}$.
We assume that\\
${}\quad$ ($I4$) the collection 
$\{\ov\cM_{1,\rho^*}^{\rho-1}\}_{\rho^*\in\A_1(I,J),\rho^*>\rho-1}$
is properly intersecting;\\
${}\quad$ ($I5$) for all $i\!\in\!I$,
$$s_{\rho-1,i}^{\,-1}(0)=\sum_{\rho^*\in\B_1(I\sqcup J;i),\rho^*>\rho-1}
  \!\!\!\!\!\!\!\!\!\!\!\!\!\!\! \ov\cM_{1,\rho^*}^{\rho-1}.$$
The assumption ($I5$) means that we will gradually be killing the zero locus of
the section~$s_i$.
We note that all five assumptions are satisfied 
if $\rho\!-\!1$ is replaced by~$0$.\\

\noindent
If $\rho$ is as above, let 
$$\ti\pi_{\rho}\!:\ov\cM_{1,(I,J)}^{\rho}\lra\ov\cM_{1,(I,J)}^{\rho-1}$$ 
be the blowup of $\ov\cM_{1,(I,J)}^{\rho-1}$ along $\ov\cM_{1,\rho}^{\rho-1}$.
We denote by $\ov\cM_{1,\rho}^{\rho}$ the corresponding exceptional divisor.
If $\rho^*\!>\!\rho$, let $\ov\cM_{1,\rho^*}^{\rho}\!\subset\!\ov\cM_{1,(I,J)}^{\rho}$ be 
the proper transform of~$\ov\cM_{1,\rho^*}^{\rho-1}$.
If 
\begin{equation}\label{rhodfn_e}
\rho=\big(I_P\!\sqcup\!J_P,\{I_k\!\sqcup\!J_k\!: k\!\in\!K\}\big)
\end{equation}
and $i\!\in\!I$, we~put
\begin{equation}\label{bundtwist_e}
L_{\rho,i}=\begin{cases}
\ti\pi_{\rho}^{\,*}L_{\rho-1,i},& \hbox{if}~ i\!\not\in\!I_P;\\
\ti\pi_{\rho}^{\,*}L_{\rho-1,i}\otimes\O(-\ov\cM_{1,\rho}^{\rho}),& \hbox{if}~ i\!\in\!I_P;
\end{cases}  \qquad
\E_{\rho}=\ti\pi_{\rho}^{\,*}\,\E_{\rho-1}\otimes\O(\ov\cM_{1,\rho}^{\rho}).
\end{equation}
The section $\ti\pi_{\rho}^{\,*}s_{\rho-1,i}$ induces a section
$$\ti{s}_{\rho,i}\in 
\Ga\big(\ov\cM_{1,(I,J)}^{\rho};\Hom(L_{\rho,i},\ti\pi_{\rho}^{\,*}\E_{\rho-1}^*)\big).$$
This section vanishes along $\ov\cM_{1,\rho}^{\rho}$,
by the inductive assumption~($I5$) if $i\!\not\in\!I_P$.
Thus, $\ti{s}_{\rho,i}$ induces a section
$$s_{\rho,i}\in \Ga\big(\ov\cM_{1,(I,J)}^{\rho};\Hom(L_{\rho,i},\E_{\rho}^*)\big).$$
We have now described the inductive step of the procedure.
It is immediate that the requirements ($I1$)-($I3$) and ($I5$) are satisfied,
with $\rho\!-\!1$ replaced by~$\rho$, are satisfied.
Corollary~\ref{subvarcoll_crl} and the assumption ($I4$) imply that 
the assumption ($I4$) with $\rho\!-\!1$ replaced by~$\rho$ is also satisfied.\\

\noindent
We conclude the blowup construction after the $\rho_{\max}$ step.
Let
$$\wt\cM_{1,(I,J)}=\ov\cM_{1,(I,J)}^{\rho_{\max}}; \qquad
\ti\E=\E_{\rho_{\max}}; \qquad
\ti{L}_i=L_{\rho_{\max},i},\quad
\ti{s}_i=s_{\rho_{\max},i} \qquad\forall\, i\!\in\!I.$$
By ($I5$), with $\rho\!-\!1$ replaced by $\rho_{\max}$, and \e_ref{curve1zero_e},
the section $\ti{s}_i$ does not vanish.
We note that by~($I1$), the stratum 
$$\cM_{1,(I,J)}\subset\ov\cM_{1,(I,J)}$$ 
consisting of the smooth curves
is a Zariski open subset of $\ov\cM_{1,(I,J)}^{\rho}$ for all
$\rho\!\in\!\{0\}\!\sqcup\!\A_1(I,J)$.\\

\noindent
By the next lemma, different extensions of the partial order $\prec$ to an order $<$
on $\A_1(I,J)$ correspond to blowing up along disjoint subvarieties in different orders.
Thus, the end result of the above blowup construction is well-defined, 
i.e.~independent of the choice of the ordering $<$ extending the partial ordering~$\prec$.

\begin{lmm}
\label{curve1bl_lmm}
Suppose $\rho,\rho'\!\in\!\A_1(I,J)$
are such that $\rho\!\not\prec\!\rho'$ and $\rho'\!\not\prec\!\rho$.
If $\rho\!\neq\!\rho'$, then 
the spaces $\ov\cM_{1,\rho}^{\ti\rho}$ and $\ov\cM_{1,\rho'}^{\ti\rho}$
are disjoint for some $\ti\rho\!\prec\!\rho,\rho'$.
\end{lmm}

\noindent
{\it Proof:} (1) Suppose
$$\rho=\big(I_P\!\sqcup\!J_P,\{I_k\!\sqcup\!J_k\!: k\!\in\!K\}\big)
\qquad\hbox{and}\qquad
\rho'=\big(I_P'\!\sqcup\!J_P',\{I_k'\!\sqcup\!J_k'\!: k\!\in\!K'\}\big).$$
For each $k\!\in\!K$ and $k'\!\in\!K'$, let
\begin{gather*}
\rho_k\!=\!\big((I\!-\!I_k)\!\sqcup\!(J\!-\!J_k),\{I_k\!\sqcup\!J_k\}\big)\in \A_1(I,J)
\qquad\hbox{and}\\
\rho_{k'}'\!=\!\big((I\!-\!I_{k'}')\!\sqcup\!(J\!-\!J_{k'}'),\{I_{k'}'\!\sqcup\!J_{k'}'\}\big)
\in \A_1(I,J).
\end{gather*}
By definition, $\ov\cM_{1,\rho_k}^0$ and $\ov\cM_{1,\rho_{k'}'}^0$
are divisors in $\ov\cM_{1,(I,J)}^0\!=\!\ov\cM_{1,I\sqcup J}$, 
$$\ov\cM_{1,\rho}^0=\bigcap_{k\in K}\ov\cM_{1,\rho_k}^0,
\qquad\hbox{and}\qquad
\ov\cM_{1,\rho'}^0=\bigcap_{k'\in K'}\ov\cM_{1,\rho_{k'}'}^0.$$
Furthermore, if $\ov\cM_{1,\rho_k}^0\!\cap\!\ov\cM_{1,\rho_{k'}'}^0\!\neq\!\eset$,
then either
$$I_k\!\sqcup\!J_k \subset I_{k'}'\!\sqcup\!J_{k'}',
\qquad\hbox{or}\qquad
I_k\!\sqcup\!J_k \supset I_{k'}'\!\sqcup\!J_{k'}',
\qquad\hbox{or}\qquad
(I_k\!\sqcup\!J_k) \cap (I_{k'}'\!\sqcup\!J_{k'}')=\eset.$$
(2) Suppose $\ov\cM_{1,\rho}^0\!\cap\!\ov\cM_{1,\rho'}^0\!\neq\!\eset$.
By the above, there exist decompositions
$$K=K_+\sqcup K_0\sqcup \bigsqcup_{l'\in K_+'}\!\!\!K_{l'}
\qquad\hbox{and}\qquad
K'=K_+'\sqcup K_0'\sqcup \bigsqcup_{l\in K_+}\!\!\!K_l'$$
and a bijection $\vph\!:K_0\!\lra\!K_0'$ such that 
\begin{gather*}
I_k\!\sqcup\!J_k \subsetneq I_{l'}'\!\sqcup\!J_{l'}'
\quad\forall\, k\!\in\!K_{l'},\, l'\!\in\!K_+',  \qquad
I_l\!\sqcup\!J_l\supsetneq I_{k'}'\!\sqcup\!J_{k'}'
\quad\forall\, k'\!\in\!K_l',\, l\!\in\!K_+,\\
\hbox{and}\qquad
I_k\!\sqcup\!J_k = I_{\vph(k)}'\!\sqcup\!J_{\vph(k)}' \quad\forall\, k\!\in\!K_0.
\end{gather*}
We note that the subsets $K_+$ and $K_+'$ of $K$ and $K'$ are nonempty.
For example, if $K_+$ were empty, then we would have 
$\rho'\!\prec\!\rho$, contrary to our assumptions.
Let
$$\ti\rho= \big( \ti{I}_P\!\sqcup\!\ti{J}_P, 
\{\ti{I}_k\!\sqcup\!\ti{J}_k\!: k\!\in\!K_0\!\sqcup\!K_+\!\sqcup\!K_+'\}\big)
\in \A_1(I,J)$$
be given by
$$\ti{I}_P \sqcup \ti{J}_P = (I_P\!\cap\!I_P') \sqcup (J_P\!\cap\!J_P'), 
\qquad
\ti{I}_k\!\sqcup\!\ti{J}_k =\begin{cases}
I_k\!\sqcup\!J_k,& \hbox{if}~ k\!\in\!K_0\!\sqcup\!K_+;\\
I_k'\!\sqcup\!J_k',& \hbox{if}~ k\!\in\!K_+'.\\
\end{cases}$$
For example, if $\rho$ corresponds to the second diagram on the right side
of Figure~\ref{partorder_fig} and $\rho'$ corresponds to either 
the first or the third diagram on the right side,
then $\ti\rho$ corresponds to the diagram on the left side of Figure~\ref{partorder_fig}.
By definition, $\ti\rho\!\prec\!\rho,\rho'$.
Furthermore,
$$\ov\cM_{1,\rho}^0\!\cap\!\ov\cM_{1,\rho'}^0 \subset \ov\cM_{1,\ti\rho}^0.$$
Thus, by Lemma~\ref{curvstr_lmm}, Corollary~\ref{subvarcoll_crl}, and~(2) 
of Lemma~\ref{ag_lmm1}, 
$$\ov\cM_{1,\rho}^{\ti\rho}\cap\ov\cM_{1,\rho'}^{\ti\rho}
\subset \ov\cM_{1,(I,J)}^{\ti\rho}$$ 
is the closure of the empty set.

\subsection{A Blowup of a Moduli Space of Genus-Zero Curves}
\label{curve0bl_subs}

\noindent
Suppose $\ale$ is a nonempty finite set and $\vr\!=\!(I_l,J_l)_{l\in\ale}$
is a tuple of finite sets such that $I_l\!\neq\!\eset$ and 
$|I_l|\!+\!|J_l|\!\ge\!2$ for all $l\!\in\!\ale$.
Let 
$$\ov\cM_{0,\vr}= \prod_{l\in\ale}\! \ov\cM_{0,\{0\}\sqcup I_l\sqcup J_l}
\qquad\hbox{and}\qquad
F_{\vr}=\bigoplus_{l\in\ale} \pi_l^*L_0 \lra \ov\cM_{0,\vr},$$
where $L_0\!\lra\!\ov\cM_{0,\{0\}\sqcup I_l\sqcup J_l}$ is 
the universal tangent line bundle for the marked point~$0$ and
$$\pi_l\!: \ov\cM_{0,\vr} \lra \ov\cM_{0,\{0\}\sqcup I_l\sqcup J_l}$$
is the projection map.
In this subsection, we construct a blowup 
$$\pi_{0,\vr}\!:\wt\cM_{0,\vr}\lra \P F_{\vr}$$
of the projective bundle $\P F_{\vr}$ over $\ov\cM_{0,\vr}$.
We also construct line bundles 
$$\ti\E,\, \ti{L}_{(l,i)}\lra \wt\cM_{0,\vr}, 
\qquad i\!\in\!I_l,\, l\!\in\!\ale,$$
and nowhere vanishing sections 
$$\ti{s}_{(l,i)}\in\Ga\big(\wt\cM_{0,\vr};\Hom(\ti{L}_{(l,i)},\ti\E^*)\big),
\qquad i\!\in\!I_l,\, l\!\in\!\ale.$$
In particular, all line bundles $\ti{L}_{(l,i)}$ and $\ti\E^*$ 
are explicitly isomorphic.
They will be denoted by~$\L$ and called the universal tangent line bundle.\\

\noindent
Similarly to the previous subsection, the smooth variety $\wt\cM_{0,\vr}$
is obtained by blowing up the subvarieties $\ov\cM_{0,\rho}$ defined below and 
their proper transforms in an order consistent with a natural partial ordering~$\prec$.
The line bundle $\ti\E$ is the sum of the tautological line bundle
$$\ga_{\vr}\lra\P F_{\vr}$$
and all exceptional divisors.
For every $l\!\in\!\aleph$ and $i\!\in\!I_l$, $\ti{L}_{(l,i)}$ is 
$\pi_l^*L_i$ minus some of these divisors.
The section~$\ti{s}_{(l,i)}$ is induced from 
the pairings~$s_i$ of Subsection~\ref{curvebldata_subs}.\\

\noindent
With $\vr$ as above and $\A_0(I_l,J_l)$ as in~\e_ref{smallcolldfn_e}, let
\begin{equation}\label{vr0setdfn_e}\begin{split}
\A_0(\vr) =\big\{ \big(\ale_+,(\rho_l)_{l\in\ale}\big)\!:~ 
&\ale_+\!\subset\!\ale,~\ale_+\!\neq\!\eset;~
\rho_l\!\in\!\{(I_l\!\sqcup\!J_l,\eset)\}\!\sqcup\!\A_0(I_l,J_l)~\forall\, l\!\in\!\ale;\\
&\rho_l\!=\!(I_l\!\sqcup\!J_l,\eset)~\forall\, l\!\in\!\ale\!-\!\ale_+; 
\big(\ale_+,(\rho_l)_{l\in\ale}\big)\!\neq\!
\big(\ale,(I_l\!\sqcup\!J_l,\eset)_{l\in\ale}\big)\big\}.
\end{split}\end{equation}
We define a partial ordering on $\A_0(\vr)$ by setting
\begin{equation}\label{rho0dfn_e}
\rho'\!\equiv\!\big(\ale_+',(\rho_l')_{l\in\ale}\big)
\prec  \rho\!\equiv\!\big(\ale_+,(\rho_l)_{l\in\ale}\big)
\end{equation}
if $\rho'\!\neq\!\rho$, $\ale_+'\!\subset\!\ale_+$,
and for every $l\!\in\!\ale$ either $\rho_l'\!=\!\rho_l$, $\rho_l'\!\prec\!\rho_l$,
or $\rho_l'\!=\!(I_l\!\sqcup\!J_l,\eset)$.
Let $<$ be an ordering on $\A_0(\vr)$ extending the partial ordering $\prec$.
We denote the corresponding minimal and maximal elements of $\A_0(\vr)$ 
by~$\rho_{\min}$ and~$\rho_{\max}$, respectively. 
If $\rho\!\in\!\A_0(\vr)$, we define
$$\rho\!-\!1 \in \{0\} \!\sqcup\! \A_0(\vr)$$
as in~\e_ref{minusdfn_e}.\\

\noindent
If $\rho\!\in\!\A_0(\vr)$ is as in \e_ref{rho0dfn_e},
let 
\begin{gather*}
\ov\cM_{0,\rho}=\prod_{l\in\ale}\ov\cM_{0,\rho_l},
\qquad 
F_{\rho}=\bigoplus_{l\in\ale^+} \pi_l^*L_0\big|_{\ov\cM_{0,\rho}}
\subset F_{\vr},\\
\hbox{and}\qquad
\wt\cM_{0,\rho}^0\!=\!\P F_{\rho} \subset \wt\cM_{0,\vr}^0\!\equiv\!\P F_{\vr}.
\end{gather*}
The spaces $\wt\cM_{0,\vr}^0$ and $\wt\cM_{0,\rho}^0$ can be represented by diagrams 
as in Figure~\ref{g0curv_fig2}.
The trees of circles attached to the vertical lines correspond to 
the tuples $\rho_l$, with conventions as 
in the first, symplectic-topology, diagram in Figure~\ref{g0curv_fig}.
For each such tree, the marked point~$0$ is the point on the line.
We indicate the elements of $\ale_+\!\subset\!\ale$ with plus signs
next to these points.
Note that by~\e_ref{vr0setdfn_e}, every dot on a vertical line
for which the corresponding tree of circles contains more than one circle
must be labeled with a plus sign.
From Lemma~\ref{curvstr_lmm}, we immediately obtain

\begin{figure}
\begin{pspicture}(-2,-2.4)(10,1.3)
\psset{unit=.4cm}
\psline[linewidth=.1](0,2.5)(0,-5.5)
\pscircle(1.2,1.5){1.2}\pscircle*(0,1.5){.25}\rput(-.8,1.5){\smsize{$+$}}
\pscircle*(2.05,2.35){.2}\pscircle*(2.05,.65){.2}
\pscircle(1.2,-1.5){1.2}\pscircle*(0,-1.5){.25}\rput(-.8,-1.5){\smsize{$+$}}
\pscircle*(1.2,-.3){.2}\pscircle*(2.4,-1.5){.2}\pscircle*(1.2,-2.7){.2}
\pscircle(1.2,-4.5){1.2}\pscircle*(0,-4.5){.25}\rput(-.8,-4.5){\smsize{$+$}}
\pscircle*(1.2,-3.3){.2}\pscircle*(2.4,-4.5){.2}\pscircle*(1.2,-5.7){.2}
\rput(7.5,-1.5){\begin{small}\begin{tabular}{l}
$\aleph\!=\!\{1,2,3\}$\\
$|I_1\!\sqcup\!J_1|\!=\!2$\\
$|I_2\!\sqcup\!J_2|\!=\!|I_3\!\sqcup\!J_3|\!=\!3$\\ 
~\\
$\aleph_+\!=\!\{1,2,3\}$\\
$\rho_1\!=\!(I_1\!\sqcup\!J_1,\eset)$\\
$\rho_2\!=\!(I_2\!\sqcup\!J_2,\eset)$\\
$\rho_3\!=\!(I_3\!\sqcup\!J_3,\eset)$\\
\end{tabular}\end{small}}
\psline[linewidth=.1](20,2.5)(20,-5.5)
\pscircle(21.2,1.5){1.2}\pscircle*(20,1.5){.25}\rput(19.2,1.5){\smsize{$-$}}
\pscircle*(22.05,2.35){.2}\pscircle*(22.05,.65){.2}
\pscircle(21.2,-1.5){1.2}\pscircle*(20,-1.5){.25}\rput(19.2,-1.5){\smsize{$+$}}
\pscircle*(21.2,-.3){.2}
\pscircle(23.4,-1.5){1}\pscircle*(22.4,-1.5){.15}
\pscircle*(24.11,-.79){.2}\pscircle*(24.11,-2.21){.2}
\pscircle(21.2,-4.5){1.2}\pscircle*(20,-4.5){.25}\rput(19.2,-4.5){\smsize{$+$}}
\pscircle*(21.2,-3.3){.2}\pscircle*(22.4,-4.5){.2}\pscircle*(21.2,-5.7){.2}
\rput(29,-1.5){\begin{small}\begin{tabular}{l}
$\aleph\!=\!\{1,2,3\}$\\
$|I_1\!\sqcup\!J_1|\!=\!2$\\
$|I_2\!\sqcup\!J_2|\!=\!|I_3\!\sqcup\!J_3|\!=\!3$\\ 
~\\
$\aleph_+\!=\!\{2,3\}$\\
$\rho_1\!=\!(I_1\!\sqcup\!J_1,\eset)$\\
$\rho_2\!\neq\!(I_2\!\sqcup\!J_2,\eset)$\\
$\rho_3\!=\!(I_3\!\sqcup\!J_3,\eset)$\\
\end{tabular}\end{small}}
\end{pspicture}
\caption{Typical Elements of $\wt\cM_{0,\vr}^0$  and $\wt\cM_{0,\rho}$}
\label{g0curv_fig2}
\end{figure}

\begin{lmm}
\label{curv0str_lmm}
Suppose $\ale$ is a nonempty finite set and $\vr\!=\!(I_l,J_l)_{l\in\ale}$
is a tuple of finite sets such that $I_l\!\neq\!\eset$ and 
$|I_l|\!+\!|J_l|\!\ge\!2$ for all $l\!\in\!\ale$.
If $\A_0(\vr)$ is as above, the collection
$\{\wt\cM_{0,\rho}\}_{\rho\in\A_0(\vr)}$ is properly intersecting.
\end{lmm}

\noindent
We now describe the starting data for the sequential blowup construction of 
this subsection.
Let
$$\E_0\!=\!\ga_{\vr} \lra \wt\cM_{0,\vr}^0\!=\!\P F_{\vr}
\qquad\hbox{and}\qquad
L_{0,(l,i)}\!=\!\pi_{0,\vr}^*\pi_l^*L_i \lra \wt\cM_{0,\vr}^0
~~~\forall\,i\!\in\!I_l,\,l\!\in\!\ale.$$
We take
$$s_{0,(l,i)}\in\Ga\big(\wt\cM_{0,\vr}^0;\Hom(L_{0,(l,i)},\E_0^*)\big)$$
to be the section induced by $\pi_{0,\vr}^{\,*}\pi_l^*s_i$, where
$s_i$ is the natural homomorphism described in Subsection~\ref{curvebldata_subs}.
It follows immediately from~\e_ref{curvezero_e} that
$$s_{0,(l,i)}^{~-1}(0)=\!\sum_{\rho^*\in\B_0(\vr;l,i)}\!\!\!\!\! \wt\cM_{0,\rho^*}^0,
\qquad\hbox{where}$$
\begin{equation*}\begin{split}
\B_0(\vr;l,i)=\Big\{\big(\ale_+,(\rho_{l'})_{l'\in\ale}\big)\!\in\!\A_0(\vr)\!: 
\ale_+\!=\!\ale\!-\!\{l\} ~\hbox{and}~
\rho_{l'}\!=\!(I_{l'}\!\sqcup\!J_{l'},\eset) ~\forall l'\!\in\!\ale, ~~~\hbox{or}&\\
\ale_+\!=\!\ale,~\rho_l\!\in\!\B_0(I_l\!\sqcup\!J_l;i),~ 
\rho_{l'}\!=\!(I_{l'}\!\sqcup\!J_{l'},\eset) ~\forall l'\!\in\!\ale\!-\!\{l\}&\Big\}.
\end{split}\end{equation*}\\

\noindent
The rest of the construction proceeds as in Subsection~\ref{curve1bl_subs}.
The analogue of~\e_ref{bundtwist_e} now~is
\begin{gather}\label{bundletwist_e2a}
L_{\rho;(l,i)}=\begin{cases}
\ti\pi_{\rho}^*L_{\rho-1,(l,i)},& \hbox{if}~ 
l\!\not\in\!\ale_+ ~\hbox{or}~
\rho_l\!\neq\!(I_l\!\sqcup\!J_l,\eset),\,i\!\not\in\!I_{l,P};\\
\ti\pi_{\rho}^*L_{\rho-1,(l,i)}\otimes
\O(-\wt\cM_{0,\rho}^{\rho}),& \hbox{otherwise};
\end{cases}  \\
\label{bundletwist_e2b}
\E_{\rho}=\ti\pi_{\rho}^*\,\E_{\rho-1}\!\otimes\!\O(\wt\cM_{0,\rho}^{\rho}).
\end{gather}
As before, we~take
\begin{gather*}
\wt\cM_{0,\vr}=\wt\cM_{0,\vr}^{\rho_{\max}}; \qquad
\ti\E=\E_{\rho_{\max}};\\
\ti{L}_{(l,i)}=L_{\rho_{\max},(l,i)} \quad\hbox{and}\quad
\ti{s}_{(l,i)}=s_{\rho_{\max},(l,i)} \qquad\forall~ i\!\in\!I_l,~l\!\in\!\ale.
\end{gather*}
The analogue of the inductive assumption ($I5$) insures that
each section $\ti{s}_{(l,i)}$ does not vanish.
The statement and the proof of Lemma~\ref{curve1bl_lmm}
remain valid in the present setting, with only minor changes.
Thus, the end result of the above blowup construction
is again well-defined, i.e.~independent of the choice of the ordering~$<$ 
extending the partial ordering~$\prec$.

\section{A Blowup of a Moduli Space of Genus-Zero Maps}
\label{map0bl_sec}

\subsection{Blowups and Immersions}
\label{map0prelim_subs1}

\noindent
In this section we construct blowups of certain moduli spaces of genus-zero maps; 
see Subsections~\ref{map0prelim_subs} and~\ref{map0blconstr_subs}.
As outlined in Subsection~\ref{outline_subs},
these blowups appear in Subsection~\ref{map1blconstr_subs} as 
the second factor in the domain of the immersions induced by
the immersions~$\io_{\si}$ of Subsection~\ref{descr_subs}.\\

\noindent
As in Section~\ref{curvebl_sec}, we begin by introducing convenient terminology 
and reviewing standard facts from algebraic geometry.
If $\ov\M$ is a variety, we denote its Zariski tangent space and its tangent cone
by $T\ov\M$ and $TC\ov\M$, respectively. 
If $X$ is a smooth variety (but not necessarily equidimensional),
we recall that a morphism $\io_X\!:X\!\lra\!\ov\M$ is an {\tt immersion} 
if the differential of~$\io_X$,
$$d\io_X\!: TX\lra \io_X^*TC\ov\M,$$
is injective at every point of $X$.
Let
$$\Im^s\,\io_X\equiv \big\{p\!\in\!\ov\M\!: |\io_X^{-1}(p)|\!\ge\!2\big\}
\qquad\hbox{and}\qquad \N_{\io_X}\equiv \io_X^*TC\ov\M\big/ \Im\, d\io_X$$
be {\tt the singular locus of $\io_X$} and
{\tt the normal cone of $\io_X$ in $\ov\M$}, respectively.
We denote~by
$$\pi_{\io_X}^{\perp}\!: \io_X^*TC\ov\M \lra \N_{\io_X}$$
the projection map.
If $Z$ is a subvariety of $\ov\M$, let
$$\io_Z\!: Z\lra\ov\M$$
the inclusion map.

\begin{dfn}
\label{ag_dfn2}
Let $\ov\M$ be a variety.\\
(1) An immersion $\io_X\!:X\!\lra\!\ov\M$ is {\tt properly self-intersecting} if
for all $x_1,x_2\!\in\!X$ such that $\io_X(x_1)\!=\!\io_X(x_2)$
and sufficiently small neighborhoods $U_1$ of $x_1$ and 
$U_2$ of $x_2$ in $X$
$$TC_{\io_X(x_1)}\big(\io_X(U_1)\!\cap\!\io_X(U_2)\big) = 
\Im\,d\io_X|_{x_1} \cap \Im\,d\io_X|_{x_2} \subset TC_{\io_X(x_1)}\ov\M.
~\footnote{We emphasize that intersections are taken to be set-theoretic
intersections unless otherwise noted.}$$
(2) If $\io_X\!:X\!\lra\!\ov\M$ and $\io_Y\!:Y\!\lra\!\ov\M$ are immersions
such that $\io_X$ is properly self-intersecting,
$\io_X$ is {\tt properly self-intersecting relative to $\io_Y$}
if for all $x_1,x_2\!\in\!X$ and $y\!\in\!Y$ such that 
$$\io_X(x_1) =\io_X(x_2)=\io_Y(y)$$
and for all sufficiently small neighborhoods $U_1$ of $x_1$ and $U_2$ of $x_2$ in $X$, 
$$\pi_{\io_Y}^{\perp}\big|_y\big(TC_{\io_Y(y)} (\io_X(U_1)\!\cap\!\io_X(U_2))\big)
=  \pi_{\io_Y}^{\perp}\big|_y\Im\,d\io_X|_{x_1} \cap 
\pi_{\io_Y}^{\perp}\big|_y\Im\,d\io_X|_{x_2} 
\subset \N_{\io_Y}\big|_y.$$
\end{dfn}

\noindent
This definition generalizes Definition~\ref{ag_dfn1};
see the paragraph following the latter for some examples.

\begin{dfn}
\label{immercoll_dfn}
If $\ov\M$ is a variety, a collection 
$\{\io_{\vr}\!:X_{\vr}\!\lra\!\ov\M\}_{\vr\in\A}$ 
of immersions is {\tt properly self-intersecting} if 
for all $\rho_1,\rho_2,\rho_3\!\in\!\A$ the immersion 
$\io_{\rho_1}\!\sqcup\!\io_{\rho_2}$ is properly 
self-intersecting relative to~$\io_{\rho_3}$. 
\end{dfn}

\begin{lmm}
\label{ag_lmm2a}
Suppose $\ov\M$ is a variety and $Z$ is a smooth subvariety of $\ov\M$.\\
(1) If $\io_X\!:X\!\lra\!\ov\M$ is an immersion such that 
the immersion $\io_X\!\sqcup\!\io_Z\!:X\!\sqcup\!Z\!\lra\!\ov\M$ is properly 
self-intersecting, then $\io_X$ lifts to an immersion
$$\Pr_Z\io_X\!: \Bl_{\io_X^{-1}(Z)}X \lra \ov\M 
\qquad\st\qquad \Im\,\Pr_Z\io_X=\Pr_Z\,\Im\,\io_X.$$
(2) If in addition $\io_X$ is properly self-intersecting relative to~$\io_Z$, 
then $\Pr_Z\io_X$ is properly self-intersecting and
$$\Im^s\,\Pr_Z\io_X=\Pr_Z\,\Im^s\,\io_X.$$
(3) If in addition $\io_Y\!:Y\!\lra\!\ov\M$ is an immersion such that 
$\io_X\!\sqcup\io_Y\!\sqcup\!\io_Z$ is properly self-intersecting and 
$\io_X$ is properly self-intersecting relative to~$\io_Y$, then $\Pr_Z\io_X$ 
is properly self-intersecting relative to~$\Pr_Z\io_Y$. Furthermore, $$\big\{\Pr_Z\io_X\big\}^{-1}\big(\Pr_Z\Im\,\io_Y\big)=
\Pr_{\io_X^{-1}(Z)}\io_X^{-1}(\Im\,\io_Y).$$\\
\end{lmm}

\noindent
{\it Remark:} Since we always require that the blowup locus be smooth,
an implicit conclusion of~(1) of Lemma~\ref{ag_lmm2a} is that $\io_X^{-1}(Z)$
is a smooth subvariety of~$X$; this is immediate from the local situation.
Note that $X$ itself is smooth, as it is the domain of the immersion~$\io_X$.

\begin{crl}
\label{immercoll_crl}
If $\ov\M$ is a variety, $\{\io_{\vr}\!:X_{\vr}\!\lra\!\ov\M\}_{\vr\in\A}$
is a properly self-intersecting collection of immersions, and
$\vr\!\in\!\A$ is such that $\io_{\vr}$ is an embedding, 
then $\{\Pr_{\Im\,\io_{\vr}}\io_{\vr'}\}_{\vr'\in\A-\{\vr\}}$
is a properly self-intersecting collection of immersions into $\Bl_{\Im\,\io_{\vr}}\ov\M$.
\end{crl}

\begin{lmm}
\label{ag_lmm2b}
Suppose $\ov\M$ is a smooth variety, $Z$ is a smooth subvariety of $\ov\M$,
$\io_X\!:X\!\lra\!\ov\M$ is an immersion such that 
the immersion $\io_X\!\sqcup\!\io_Z$ is properly self-intersecting.
Let
$$\io_X^{-1}(Z) = \bigsqcup_{\vr\in\A} Z_{\vr}$$
be the decomposition of $\io_X^{-1}(Z)$ into path components.
If there exist a splitting
$$\N_{\io_X}=\bigoplus_{i\in I}L_i \lra X$$
and a subset $I_{\vr}$ of $I$ for each $\vr\!\in\!\A$ such that
\begin{equation}\label{ag_lmm2b_e}
\io_X|_{Z_{\vr}}^{\,*}TZ\big/TZ_{\vr}
=\bigoplus_{i\in I-I_{\vr}}\!\!\!L_i|_{Z_{\vr}} \qquad\forall\,\vr\!\in\!A,
\end{equation}
then 
$$\N_{\Pr_Z\io_X}=\bigoplus_{i\in I}
\Big(\pi^*L_i\otimes \bigotimes_{i\in I_{\vr}}\O(-E_{\vr})\Big),$$
where $E_{\vr}$ is the component of the exceptional divisor for the blowup
$\pi\!: \Bl_{\io_X^{-1}(Z)}X\!\lra\!X$ that projects onto~$Z_{\vr}$.
\end{lmm}

\noindent
We note that by (1) of Definition~\ref{ag_dfn2}, the homomorphism
$$\io_X|_{Z_{\vr}}^{\,*}TZ\big/TZ_{\vr} \lra 
\N_{\io_X}\!\equiv\!\io_X^*T\ov\M\big/ \Im\, d\io_X$$
induced by the inclusions is injective. 
Thus, we can identify $\io_X|_{Z_{\vr}}^{\,*}TZ\big/TZ_{\vr}$
with a subbundle of~$\N_{\io_X}$, as we have done in Lemma~\ref{ag_lmm2b}.

\subsection{Moduli Spaces of Genus-Zero Maps}
\label{map0str_subs}

\noindent
In this subsection, we describe natural subvarieties of the moduli space of
genus-zero maps and a natural bundle section over them.
This bundle section induces other bundle sections, introduced in the next two subsections,
that are used in the blowup construction of Subsection~\ref{map1blconstr_subs} to describe 
the structure of the proper transforms of~$\ov\M_{1,k}^0(\Pn,d)$;
see Subsection~\ref{outline_subs} for more details.
Below we also state two well-known facts in the Gromov-Witten theory,
Lemmas~\ref{g0mapstr_lmm1} and~\ref{g0mapstr_lmm2}, and 
a more recent result, Lemma~\ref{deriv0str_lmm}.\\ 

\noindent
If $d\!\in\!\Z^+$ and $J$ is a finite set, let
\begin{gather}\label{g0map_e}
\A_0(d,J)=\big\{(m;J_P,J_B)\!: m\!\in\!\Z^+,\,m\!\le\!d;~
J_P\!\subset\!J; m\!+\!|J_P|\!\ge\!2\big\};\\
\ov\M_{0,(0,J)}(\Pn,d)=\ov\M_{0,\{0\}\sqcup J}(\Pn,d).\notag
\end{gather}
If $\si\!=\!(m;J_P,J_B)$ is an element of $\A_0(d,J)$, 
let $\M_{0,\si}(\Pn,d)$ be the subset of $\ov\M_{0,\{0\}\sqcup J}(\Pn,d)$ 
consisting of the stable maps $[\Si,u]$ such~that\\
${}\quad$ (i) the components of $\Si$ are $\Si_i\!=\!\Bbb{P}^1$ 
with $i\!\in\!\{P\}\!\sqcup\![k]$;\\
${}\quad$ (ii) $u|_{\Si_P}$ is constant and 
the marked points on $\Si_P$ are indexed by the set $\{0\}\!\sqcup\!J_P$;\\
${}\quad$ (iii) for each $i\!\in\![m]$, 
$\Si_i$ is attached to $\Si_P$ and $u|_{\Si_i}$ is not constant.\\
We denote by $\ov\M_{0,\si}(\Pn,d)$ the closure of 
$\M_{0,\si}(\Pn,d)$ in $\ov\M_{0,\{0\}\sqcup J}(\Pn,d)$.
Figure~\ref{g0map_fig} illustrates this definition,
from the points of view of symplectic topology and of algebraic geometry.
In the first diagram, each disk represents a sphere,
and we shade the components on which the map~$u$ is non-constant.
In the second diagram, the irreducible components of $\Si$ are represented by lines,
and the integer next to each component shows the degree of $u$ on that component.
In both cases, we indicate the marked points lying on the component $\Si_P$ only.\\

\begin{figure}
\begin{pspicture}(-1.1,-2.2)(10,1.5)
\psset{unit=.4cm}
\pscircle(5,-1.5){1.5}
\pscircle[fillstyle=solid,fillcolor=gray](2.5,-1.5){1}\pscircle*(3.5,-1.5){.15}
\pscircle[fillstyle=solid,fillcolor=gray](5,1){1}\pscircle*(5,0){.15}
\pscircle[fillstyle=solid,fillcolor=gray](7.5,-1.5){1}\pscircle*(6.5,-1.5){.15}
\pscircle*(5,-3){.2}\rput(4.9,-3.6){\smsize{$0$}}
\pscircle*(6.06,-2.56){.2}\rput(6.6,-2.9){\smsize{$j_1$}}
\psline(18,2)(18,-4)\pscircle*(18,1.5){.2}\rput(18.6,1.6){\smsize{$j_1$}}
\pscircle*(18,-3.2){.2}\rput(18.5,-3.3){\smsize{$0$}}
\psline(16.8,0)(22.05,1.25)\rput(22.6,1.2){\smsize{$d_1$}}
\psline(17,-1)(22,-1)\rput(22.5,-1.1){\smsize{$d_2$}}
\psline(16.8,-2)(22.05,-3.25)\rput(22.6,-3.3){\smsize{$d_3$}}
\rput(32,-1){\smsize{\begin{tabular}{l}
$m\!=\!3,\,J_P\!=\!\{j_1\}$\\ 
$d_1,d_2,d_3\!>\!0$\\ $d_1\!+\!d_2\!+\!d_3\!=\!d$\\
\end{tabular}}}
\end{pspicture}
\caption{A Typical Element of $\ov\M_{0,\si}(\Pn,d)$}
\label{g0map_fig}
\end{figure}

\noindent 
We define a partial ordering on the set $\A_0(d,J)$ by setting
\begin{equation}\label{partorder_e2}
\si'\!\equiv\!(m';J_P',J_B') \prec \si\!\equiv\!(m;J_P,J_B)
\qquad\hbox{if}\quad \si'\!\neq\!\si,~m'\!\le\!m,~J_P'\!\subset\!J_P.
\end{equation}
Similarly to Subsection~\ref{curvebldata_subs}, 
this condition means that the elements of $\M_{0,\si'}(\Pn,d)$ can be obtained
from the elements of $\M_{0,\si}(\Pn,d)$ by moving more points onto the bubble components
or combining the bubble components; see Figure~\ref{partorder_fig2}.
As in the $g\!=\!0$ case of Subsection~\ref{curvebldata_subs}, 
the bubble components are the components not containing the marked point~$0$.

\begin{figure}
\begin{pspicture}(-0.3,-2.2)(10,1.5)
\psset{unit=.4cm}
\pscircle(10,-1.5){1.5}
\pscircle[fillstyle=solid,fillcolor=gray](8.23,.27){1}\pscircle*(8.94,-.44){.15}
\pscircle[fillstyle=solid,fillcolor=gray](11.77,.27){1}\pscircle*(11.06,-.44){.15}
\pscircle*(10,-3){.2}\rput(9.9,-3.6){\smsize{$0$}}
\rput(15,-1){\begin{Huge}$\prec$\end{Huge}}
\pscircle(22,-1.5){1.5}
\pscircle[fillstyle=solid,fillcolor=gray](19.5,-1.5){1}\pscircle*(20.5,-1.5){.15}
\pscircle[fillstyle=solid,fillcolor=gray](22,1){1}\pscircle*(22,0){.15}
\pscircle*(23.5,-1.5){.2}\rput(24.2,-1.5){\smsize{$j_1$}}
\pscircle*(22,-3){.2}\rput(21.9,-3.6){\smsize{$0$}}
\pscircle(30,-1.5){1.5}
\pscircle[fillstyle=solid,fillcolor=gray](27.5,-1.5){1}\pscircle*(28.5,-1.5){.15}
\pscircle[fillstyle=solid,fillcolor=gray](30,1){1}\pscircle*(30,0){.15}
\pscircle[fillstyle=solid,fillcolor=gray](32.5,-1.5){1}\pscircle*(31.5,-1.5){.15}
\pscircle*(30,-3){.2}\rput(29.9,-3.6){\smsize{$0$}}
\end{pspicture}
\caption{Examples of Partial Ordering~\e_ref{partorder_e2}}
\label{partorder_fig2}
\end{figure}

\begin{lmm}
\label{g0mapstr_lmm1}
If $\si_1,\si_2\!\in\!\A_0(d,J)$, $\si_1\!\neq\!\si_2$
$\si_1\!\not\prec\!\si_2$, and $\si_2\!\not\prec\!\si_1$, then
\begin{gather*}
\ov\M_{0,\si_1}(\Pn,d) \cap \ov\M_{0,\si_2}(\Pn,d)
\subset \ov\M_{0,\ti\si(\si_1,\si_2)}(\Pn,d),\\
\hbox{where}\qquad
\ti\si(\si_1,\si_2)=\max\big\{\si'\!\in\!\A_0(d,J)\!: \si'\!\prec\!\si_1,\si_2\big\}.
\end{gather*}
If $\ti\si(\si_1,\si_2)$ is not defined, $\ov\M_{0,\si_1}(\Pn,d)$ and 
$\ov\M_{0,\si_2}(\Pn,d)$ are disjoint.
\end{lmm}

\noindent
For example, if $\si_1$ and $\si_2$ correspond to the two diagrams on the right of 
Figure~\ref{partorder_fig2}, then $\ti\si(\si_1,\si_2)$ corresponds to 
the diagram on the left of Figure~\ref{partorder_fig2}.\\

\noindent
If $\si\!=\!(m;J_P,J_B)$ is an element of $\A_0(d,J)$, let 
\begin{gather*}
\ov\M_{\si;B}(\Pn,d) \subset \prod_{i\in[m]}~\bigsqcup_{d_i>0,J_i\subset J_B}
\!\!\!\!\!\!\!\!\ov\M_{0,\{0\}\sqcup J_i}(\Pn,d_i)    \qquad\hbox{and}\\
\pi_i\!: \ov\M_{\si;B}(\Pn,d)\lra
\!\!\bigsqcup_{d_i>0,J_i\subset J_B}\!\!\!\!\!\!\!\!\ov\M_{0,\{0\}\sqcup J_i}(\Pn,d_i),
\quad i\!\in\![m],
\end{gather*}
be as in Subsection~\ref{descr_subs}.
Since each of the spaces $\ov\M_{0,\{0\}\sqcup J_i}(\Pn,d_i)$ is smooth and
each of the evaluation maps
$$\ev_0\!: \ov\M_{0,\{0\}\sqcup J_i}(\Pn,d_i)\lra \Pn$$
is a submersion, the space $\ov\M_{\si;B}(\Pn,d)$ is smooth.
We denote~by
\begin{equation}\label{g0nodeiden_e}
\io_{\si}\!: \ov\cM_{0,\{0\}\sqcup[m]\sqcup J_P}\!\times\!\ov\M_{\si;B}(\Pn,d)
\lra \ov\M_{0,\si}(\Pn,d) \subset \ov\M_{0,\{0\}\sqcup J}(\Pn,d)
\end{equation}
the natural node-identifying map.
It descends to an immersion
$$\bar\io_{\si}\!: 
\big(\ov\cM_{0,\{0\}\sqcup[m]\sqcup J_P}\!\times\!\ov\M_{\si;B}(\Pn,d)\big)
\big/S_m \lra \ov\M_{0,\{0\}\sqcup J}(\Pn,d).$$
Let
$$\pi_P,\pi_B\!: \ov\cM_{0,\{0\}\sqcup [m]\sqcup J_P}\!\times\!\ov\M_{\si;B}(\Pn,d)
\lra \ov\cM_{0,\{0\}\sqcup [m]\sqcup J_P},\ov\M_{\si;B}(\Pn,d)$$
be the natural projection maps.
The following lemma can be easily deduced from~\cite{P}.

\begin{lmm}
\label{g0mapstr_lmm2}
If $d\!\in\!\Z^+$ and $J$ is a finite set,  the collections 
$\{\io_{\si}\}_{\si\in\A_0(d,J)}$  and $\{\bar\io_{\si}\}_{\si\in\A_0(d,J)}$
of immersions are properly self-intersecting.
If $\si\!\in\!\A_0(d,J)$ is as in~\e_ref{partorder_e2},
$$\Im^s\,\bar\io_{\si} \subset \bigcup_{\si'\prec\si}\! \ov\M_{0,\si'}(\Pn,d)
\qquad\hbox{and}\qquad 
\N_{\io_{\si}}=\bigoplus_{i\in[m]} \pi_P^*L_i\!\otimes\!\pi_B^*\pi_i^*L_0.$$
If in addition $\si'\!\in\!\A_0(d,J)$, $\si'\!\prec\!\si$, 
and $\si'$ is as in~\e_ref{partorder_e2}, then
\begin{gather*}
\io_{\si}^{-1} \big( \ov\M_{0,\si'}(\Pn,d) \big)
=\Big(\bigcup_{\rho\in\A_0(\si;\si')} \!\!\!\!\!\!\!\ov\cM_{0,\rho}\Big)
\times \ov\M_{\si;B}(\Pn,d),
\qquad\hbox{where}\\
\A_0(\si;\si')=\big\{ \rho\!=\!\big(I_P\!\sqcup\!J_P',\{I_k\!\!\sqcup\!J_k\!: k\!\in\!K\})
\!\in\!\A_0([m],J_P)\!: |K|\!+\!|I_P|\!=\!m'\big\}
\end{gather*}
and $\A_0([m],J_P)$ and $\ov\cM_{0,\rho}$ are as in Subsection~\ref{curvebldata_subs}.
Finally, if $\rho\!\in\!\A_0(\si;\si')$ is as above,
$$\io_{\si}  \big|_{\ov\cM_{0,\rho}\times\ov\M_{\si;B}(\Pn,d)}^* 
T\ov\M_{0,\si'}(\Pn,d)
\big/T\big(\ov\cM_{0,\rho}\!\times\!\ov\M_{\si;B}(\Pn,d)\big)
 =\bigoplus_{i\in[m]-I_P}\!\!\!\!\!\!\pi_P^*L_i\!\otimes\!\pi_B^*\pi_i^*L_0.$$\\
\end{lmm}

\noindent
We finish this subsection by describing a natural bundle section 
$$\cD_0\in\Ga\big(\ov\M_{0,\{0\}\sqcup J}(\Pn,d),\Hom(L_0;\ev_0^*T\Pn)\big)$$
which plays a central role in the rest of the paper.
An element $[b]\!\in\!\ov\M_{0,\{0\}\sqcup J}(\Pn,d)$ consists of a prestable nodal 
curve $\Si$ with marked points and a map $u\!:\Si\!\lra\!\Pn$.
One of the marked points is labeled by~$0$. We denote it by~$x_0(b)$.
We define $\cD_0$ by 
$$\cD_0\big|_b = du|_{x_0(b)}\!: T_{x_0(b)}\Si\lra T_{\ev_0(b)}\Pn.$$
If $\U\!\lra\!\ov\M_{0,\{0\}\sqcup J}(\Pn,d)$ is the universal curve and 
$\ev\!:\U\!\lra\!\Pn$ is the natural evaluation map,
then $\cD_0|_b$ is simply the restriction of $d\ev|_{x_0(b)}$ 
to the vertical tangent bundle of~$\U$.
The bundle section~$\cD_0$ vanishes identically along the subvarieties 
$\ov\M_{0,\si}(\Pn,d)$ with  $\si\!\in\!\A_0(d,J)$.

\begin{lmm}
\label{deriv0str_lmm}
If $d\!\in\!\Z^+$ and $J$ is a finite set,
the section $\cD_0$ is transverse to the zero set
on the complement of the subvarieties $\ov\M_{0,\si}(\Pn,d)$ 
with  $\si\!\in\!\A_0(d,J)$.
Furthermore, for every 
$$\si\!\equiv\!(m;J_P,J_B) \in \A_0(d,J),$$
the differential of~$\cD_0$,
$$\na\cD_0\!: \N_{\io_{\si}} \lra \io_{\si}^*\,\Hom(L_0,\ev_0^*T\Pn)
=\pi_P^*L_0^*\!\otimes\!\pi_B^*\ev_0^*T\Pn,$$
in the normal direction to the immersion $\io_{\si}$ is given~by
$$ \na\cD_0\big|_{\pi_P^*L_i\otimes\pi_B^*\pi_i^*L_0}
=  \pi_P^*s_i\!\otimes\!\pi_B^*\pi_i^*\cD_0
\qquad\forall\, i\!\in\![m],$$
where $s_i$ is the homomorphism defined in Subsection~\ref{curvebldata_subs}.
\end{lmm}

\noindent
The first claim of the lemma is an immediate consequence of the fact that
$$H^1\big(\Si;u^*T\Pn\!\otimes\!\O(-2z)\big)=\{0\}$$
for every genus-zero stable map $(\Si,u)$ and a smooth point $z\!\in\!\Si$
such that the restriction of $u$ to the irreducible component of $\Si$
containing $z$ is not constant.
The second statement of the lemma follows from Theorem~2.8 in~\cite{g2n2and3}.

\subsection{Initial Data}
\label{map0prelim_subs}

\noindent
If $\ale$ and $J$ are finite sets and $d$ is positive integer, let
\begin{equation*}\begin{split}
\ov\M_{0,(\aleph,J)}(\Pn,d)&= 
\Big\{(b_l)_{l\in\ale}\in\prod_{l\in\ale}\ov\M_{0,\{0\}\sqcup J_l}(\Pn,d_l)\!:
d_l\!\in\!\Z^+,~\sum_{l\in\ale}d_l\!=\!d;~ \bigsqcup_{l\in\ale} J_l\!=\!J;\\
&\qquad\qquad\qquad\qquad\qquad\qquad\qquad\qquad\qquad\qquad~~
\ev_0(b_l)\!=\!\ev_0(b_{l'})~\forall\, l,l'\!\in\!\ale\Big\};
\end{split}\end{equation*}
\begin{equation*}\begin{split}
\M_{0,(\aleph,J)}(\Pn,d)&= 
\Big\{(b_l)_{l\in\ale}\in\prod_{l\in\ale}\M_{0,\{0\}\sqcup J_l}(\Pn,d_l)\!:
d_l\!\in\!\Z^+,~\sum_{l\in\ale}d_l\!=\!d;~ \bigsqcup_{l\in\ale} J_l\!=\!J;\\
&\qquad\qquad\qquad\qquad\qquad\qquad\qquad\qquad\qquad\qquad~~
\ev_0(b_l)\!=\!\ev_0(b_{l'})~\forall\, l,l'\!\in\!\ale\Big\},
\end{split}\end{equation*}
where $\M_{0,\{0\}\sqcup J_l}(\Pn,d_l)$ is the subset of 
$\ov\M_{0,\{0\}\sqcup J_l}(\Pn,d_l)$ consisting of stable maps
with smooth domains.
For each $l\!\in\!\ale$, let 
$$\pi_l\!: \ov\M_{0,(\ale,J)}(\Pn,d)\lra
\!\!\bigsqcup_{d_l>0,J_l\subset J}\!\!\!\!\ov\M_{0,\{0\}\sqcup J_l}(\Pn,d_l)$$
be the projection map.
We put
$$F_{(\ale,J)}=\bigoplus_{l\in\ale} \pi_l^*L_0,$$
where $L_0\!\lra\!\ov\M_{0,\{0\}\sqcup J_l}(\Pn,d_l)$
is the universal tangent line bundle for the marked point $0$.
In the next subsection, we construct a blowup 
$$\pi_{0,(\ale,J)}\!: \wt\M_{0,(\ale,J)}(\Pn,d) \lra \P F_{(\ale,J)}$$
of the projective bundle $\P F_{(\ale,J)}$ over $\ov\M_{0,(\ale,J)}(\Pn,d)$
and a line bundle
$$ \ti\E \lra \wt\M_{0,(\ale,J)}(\Pn,d).$$
We also describe a natural bundle section
$$\wt\cD_{(\ale,J)} \in \Ga\big( \wt\M_{0,(\ale,J)}(\Pn,d);
 \ti\E^*\!\otimes\!\pi_{0,(\ale,J)}^*\pi_{\P F_{(\ale,J)}}^*\ev_0^*T\Pn\big),$$
where 
$$\pi_{\P F_{(\ale,J)}}\!: \P F_{(\ale,J)} \lra \ov\M_{0,(\ale,J)}(\Pn,d)$$
is the bundle projection map.
This section is transverse to the zero~set.\\

\noindent
Similarly to Subsection~\ref{curve0bl_subs}, the smooth variety 
$\wt\M_{0,(\ale,J)}(\Pn,d)$ is obtained by blowing up the subvarieties 
$\wt\M_{0,\vr}^0(\Pn,d)$ defined below and their proper transforms in 
an order consistent with a natural partial ordering~$\prec$.
The line bundle $\ti\E$ is the sum of the tautological line bundle
$$\ga_{(\ale,J)}\lra\P F_{(\ale,J)}$$
and all exceptional divisors.
The section~$\wt\cD_{(\ale,J)}$ is induced from the sections $\pi_l^*\cD_0$,
with $l\!\in\!\ale$, where $\cD_0$ is as in Subsection~\ref{map0str_subs}.\\

\noindent
If $\ale$, $J$, and $d$ are as above, let 
\begin{equation*}\begin{split}
\A_0(\ale;d,J) = \big\{ \big((\si_l)_{l\in\ale},J_B\big)\!:
\,&(\si_l,\eset)\!\in\!\{(0,\eset)\}\!\sqcup\!\A_0(d_l,J_{l,P}),~
(\si_l)_{l\in\ale}\!\neq\!(0)_{l\in\ale};\\
&\sum_{l\in\ale}d_l\!=\!d,~\bigsqcup_{l\in\ale}J_{l,P}\!=\!J\!-\!J_B\big\}.
\end{split}\end{equation*}
We define a partial ordering $\prec$ on $\A_0(\ale;d,J)$ by setting
\begin{equation}\label{vrdfn_e2}
\vr'\!\equiv\!\big((\si_l')_{l\in\ale},J_B'\big) 
\prec \vr\!\equiv\!\big((\si_l)_{l\in\ale},J_B)
\end{equation}
if $\vr'\!\neq\!\vr$ and for every $l\!\in\!\aleph$ either 
$\si_l'\!=\!\si_l$, $(\si_l',\eset)\!\prec\!(\si_l,\eset)$, or
$\si_l'\!=\!0$.
If $\vr\!\in\!\A_0(\ale;d,J)$ is as in~\e_ref{vrdfn_e2}, we~put
$$\ale_P(\vr)=\big\{l\!\in\!\aleph\!: \si_l\!\neq\!0\big\}
\qquad\hbox{and}\qquad
\ale_S(\vr)=\big\{l\!\in\!\aleph\!: \si_l\!=\!0\big\}.$$
Here $P$ and $S$ stand for the subsets of principal and secondary elements
of~$\aleph$, respectively; see the next paragraph.
Note that 
\begin{gather}
\vr'\prec \vr ~~~\Lra~~~ \ale_P(\vr')\subset\ale_P(\vr), \qquad
\ale_P(\vr)\neq\eset~~~\forall\,\vr\!\in\!\A_0(\ale;d,J),\qquad\hbox{and}\notag\\
\label{vrdfn_e3}
\vr=\big((m_l;J_{l,P})_{l\in\ale_P(\vr)},(0)_{l\in\ale_S(\vr)},J_B\big)
\end{gather}
for some $m_l$ and $J_{l,P}$.
Choose an ordering $<$ on $\A_0(\ale;d,J)$ extending the partial ordering~$\prec$.
We denote the corresponding minimal and maximal element by~$\vr_{\min}$
and~$\vr_{\max}$, respectively.
For every  $\vr\!\in\!\A_0(\ale;d,J)$, define
$$\vr\!-\!1 \in \{0\}\!\sqcup\!\A_0(\ale;d,J)$$
as in~\e_ref{minusdfn_e}.\\

\noindent
If $\vr\!\in\!\A_0(\ale;d,J)$ is as in \e_ref{vrdfn_e2}, let
\begin{equation*}\begin{split}
\ov\M_{0,\vr}(\Pn,d) = 
\Big\{(b_l)_{l\in\ale}\in\prod_{l\in\ale}\ov\M_{0,(\si_l,J_{l,B})}(\Pn,d_l)\!:
\sum_{l\in\ale}d_l\!=\!d;~\bigsqcup_{l\in\ale}J_{l,B}\!=\!J_B;\quad&\\
\ev_0(b_{l_1})\!=\!\ev_0(b_{l_2})~\forall\, l_1,l_2\!\in\!\ale&\Big\}
\subset\ov\M_{0,(\ale,J)}(\Pn,d).
\end{split}\end{equation*}
With $*\!=\!P,S$, we define
$$F_{\vr;*}=\bigoplus_{l\in\ale_*(\vr)}
\!\!\! \pi_l^*L_0\Big|_{\ov\M_{0,\vr}(\Pn,d)}
\subset F_{(\ale,J)}\big|_{\ov\M_{0,\vr}(\Pn,d)}.$$
Let 
$$\wt\M_{0,\vr}^0(\Pn,d)=\P F_{\vr;P}
\subset \wt\M_{0,(\ale,J)}^0(\Pn,d)\!\equiv\!\P F_{(\ale,J)}.$$
From Lemma~\ref{g0mapstr_lmm1}, we immediately obtain 

\begin{lmm}
\label{map0bl_lmm1a}
If $\vr_1,\vr_2\!\in\!\A_0(\ale;d,J)$, $\vr_1\!\neq\!\vr_2$
$\vr_1\!\not\prec\!\vr_2$, and $\vr_2\!\not\prec\!\vr_1$, then
\begin{gather*}
\wt\M_{0,\vr_1}^0(\Pn,d) \cap \wt\M_{0,\vr_2}^0(\Pn,d)
\subset \wt\M_{0,\ti\vr(\vr_1,\vr_2)}^0(\Pn,d),\\
\hbox{where}\qquad
\ti\vr(\vr_1,\vr_2)=\max\big\{\vr'\!\in\!\A_0(\ale;d,J)\!: \vr'\!\prec\!\vr_1,\vr_2\big\}.
\end{gather*}
If $\ti\vr(\vr_1,\vr_2)$ is not defined, 
$\wt\M_{0,\vr_1}^0(\Pn,d)$ and $\wt\M_{0,\vr_2}^0(\Pn,d)$ are disjoint.
\end{lmm}

\noindent
With $\vr$ as \e_ref{vrdfn_e3}, let
$$\vr_P=\big([m_l],J_{l,P}\big)_{l\in\ale_P(\vr)}, \quad
\ale_B(\vr)=\ale_S(\vr)\sqcup\!\bigsqcup_{l\in\ale_P(\vr)}\!\!\!\![m_l], \quad
J_B(\vr)=J_B, \quad\hbox{and}\quad G_{\vr}=\prod_{l\in\ale_P(\vr)}\!\!\!\!\!S_{m_l}.$$
With $\wt\cM_{0,\vr_P}^0$ as in Subsection~\ref{curve0bl_subs},
we denote by
$$\io_{0,\vr}\!: \wt\cM_{0,\vr_P}^0 \times 
\ov\M_{0,(\ale_B(\vr),J_B(\vr))}(\Pn,d) \lra \wt\M_{0,\vr}^0(\Pn,d)
\subset \wt\M_{0,(\ale,J)}^0(\Pn,d)$$
the natural node-identifying map induced by the immersions $\io_{(\si_l,J_{l,B})}$
in~\e_ref{g0nodeiden_e}.
It descends to an immersion
$$\bar\io_{0,\vr}\!: 
\big(\wt\cM_{0,\vr_P}^0 \!\times\! 
\ov\M_{0,(\ale_B(\vr),J_B(\vr))}(\Pn,d)\big)\big/G_{\vr}
\lra \wt\M_{0,(\ale,J)}^0(\Pn,d).$$
Let 
$$\pi_P,\pi_B\!: \wt\cM_{0,\vr_P}^0 \times  \ov\M_{0,(\ale_B(\vr),J_B(\vr))}(\Pn,d) 
\lra \wt\cM_{0,\vr_P}^0,\ov\M_{0,(\ale_B(\vr),J_B(\vr))}(\Pn,d)$$
be the projection maps.\\

\noindent
For the rest of this section, as well as for Section~\ref{map1bl_sec}, we take
\begin{gather*}
\ov\M_{0,(\ale,J)}=\ov\M_{0,(\ale,J)}(\Pn,d), ~~~
\M_{0,(\ale,J)}=\M_{0,(\ale,J)}(\Pn,d), ~~~
 \wt\M_{0,(\ale,J)}^0=\wt\M_{0,(\ale,J)}^0(\Pn,d) \quad\forall~(\ale,J);\\
\wt\M_{0,\vr}^0=\wt\M_{0,\vr}^0(\Pn,d)
\qquad\forall~\vr\!\in\!\A_0(\ale;d,J).
\end{gather*}

\begin{lmm}
\label{map0bl_lmm1}
If $\ale$ and $J$ are finite sets and $d\!\in\!\Z^+$, the collections 
$$\{\io_{0,\vr}\}_{\vr\in\A_0(\ale;d,J)}  \qquad\hbox{and}\qquad 
\{\bar\io_{0,\vr}\}_{\vr\in\A_0(\aleph;d,J)}$$
of  immersions are properly self-intersecting.
If 
$$\vr^*\!\equiv\!\big((m_l^*;J_{l,P}^*)_{l\in\ale_P(\vr^*)},
(0)_{l\in\ale_S(\vr^*)},J_B^*\big)  \in \A_0(\aleph;d,J),$$
then
\begin{gather*}
\Im^s\,\bar\io_{0,\vr^*} \subset \bigcup_{\vr'\prec\vr^*}\! \wt\M_{0,\vr'}^0
\qquad\hbox{and}\\
\N_{\io_{0,\vr^*}}= \!\bigoplus_{l\in\ale_P(\vr^*)} \bigoplus_{i\in[m_l^*]}
\!\!\pi_P^*L_{0,(l,i)}\!\otimes\!\pi_B^*\pi_{(l,i)}^*L_0
~\oplus~
\pi_P^*\E_0^*\!\otimes\pi_B^*\!\bigoplus_{l\in\ale_S(\vr^*)}\!\!\!\!\!\pi_l^*L_0,
\end{gather*}
where $\E_0$ is as in Subsections~\ref{curve0bl_subs}.
If $\vr,\vr^*\!\in\!\A_0(\ale;d,J)$,
$\vr$ is as \e_ref{vrdfn_e3}, and $\vr\!\prec\!\vr^*$, then
$$\io_{0,\vr^*}^{~-1} \big( \wt\M_{0,\vr}^0\big)
=\Big(\bigcup_{\rho\in\A_0(\vr^*;\vr)} \!\!\!\!\!\!\!  \wt\cM_{0,\rho}^0 \Big)
\times \ov\M_{0,(\ale_B(\vr^*),J_B(\vr^*))},$$
where
\begin{equation*}\begin{split}
\A_0(\vr^*;\vr)=\Big\{ \rho\!=\!\big(\ale_P(\vr),
\big(I_{l,P}^*\!\sqcup\!J_{l,P},\{I_{l,k}^*\!\sqcup\!J_{l,k}^*\!: k\!\in\!K_l^*\}
\big)_{l\in\ale_P(\vr^*)}\big)\!\in\!\A_0(\vr_P^*)\!: \qquad\qquad&\\
|K_l^*|\!+\!|I_{l,P}^*|\!=\!m_l~\forall\, l\!\in\!\ale_P(\vr^*)&\Big\}
\end{split}\end{equation*}
and $\A_0(\vr_P^*)$ and $\wt\cM_{0,\rho}^0$ are as in Subsection~\ref{curve0bl_subs}.
Finally, if $\rho\!\in\!\A_0(\vr^*;\vr)$ is as above,
\begin{equation*}\begin{split}
&\io_{0,\vr^*}  \big|_{\wt\cM_{0,\rho}^0\times\ov\M_{0,(\ale_B(\vr^*),J_B(\vr^*))}}^* 
T\ov\M_{0,\vr}
\big/T\big(\wt\cM_{0,\rho}^0\!\times\!\ov\M_{0,(\ale_B(\vr^*),J_B(\vr^*))}\big)\\
&\qquad 
 =\!\bigoplus_{l\in\ale_P(\vr^*)-\ale_P(\vr)} \bigoplus_{i\in[m_l^*]}
\!\!\pi_P^*L_{0,(l,i)}\!\otimes\!\pi_B^*\pi_{(l,i)}^*L_0
~\oplus~\!\bigoplus_{l\in\ale_P(\vr)} \bigoplus_{i\in[m_l^*]-I_{l,P}}
\!\!\!\!\!\!\!\!\!\pi_P^*L_{0,(l,i)}\!\otimes\!\pi_B^*\pi_{(l,i)}^*L_0.
\end{split}\end{equation*}\\
\end{lmm}

\noindent 
The normal bundle $\N_{\io_{0,\vr^*}}$ for the immersion $\io_{0,\vr^*}$ 
splits into horizontal and vertical bundles:
$$\N_{\io_{0,\vr^*}} = \N_{\io_{0,\vr^*}}^{\bot} \oplus \N_{\io_{0,\vr^*}}^{\top}.$$
It is immediate from the definitions that
$$\N_{\io_{0,\vr^*}}^{\top}=
\io_{0,\vr^*}^{\,*}\big(\ga_{(\aleph,J)}^{\,*}\!\otimes\!F_{\vr;S}\big)
=\pi_P^*\E_0^*\!\otimes\pi_B^*\!\!\bigoplus_{l\in\ale_S(\vr^*)}\!\!\!\!\!\pi_l^*L_0.$$ 
The horizontal normal bundle $\N_{\io_{0,\vr^*}}^{\bot}$ is the pullback of the normal bundle
for the node identifying immersion
$$\ov\cM_{0,\vr_P^*} \times 
\ov\M_{0,(\ale_B(\vr^*),J_B(\vr^*))} \lra \ov\M_{0,\vr^*}(\Pn,d)
\subset \ov\M_{0,(\ale,J)}$$
induced by the immersions $\io_{(\si_l^*,J_{l,B}^*)}$  in~\e_ref{g0nodeiden_e}
by the bundle projection map~$\pi_{\P F_{\vr;P}}$.
The normal bundle for this immersion is 
the sum of component-wise normal bundles given by Lemma~\ref{g0mapstr_lmm2}.
The remaining claims of Lemma~\ref{map0bl_lmm1} follow easily from the corresponding 
statements of Lemma~\ref{g0mapstr_lmm2} as~well.\\

\noindent
We note that for every $\vr^*\!\in\!\A_0(\ale;d,J)$,
$$\A_0(\vr_P^*)= \bigsqcup_{\vr\prec\vr^*}\!\A_0(\vr^*;\vr).$$
Furthermore, if $\vr_1,\vr_2\!\in\!\A_0(\ale;d,J)$ are such that
$\vr_1,\vr_2\!\prec\!\vr^*$, then
$$\rho_1^*\!\in\!\A_0(\vr^*;\vr_1), \quad \rho_2^*\!\in\!\A_0(\vr^*;\vr_2), \quad
\rho_1^*\!\prec\!\rho_2^*  \qquad\Lra\qquad \vr_1\!\prec\!\vr_2.$$
Thus, we can choose an ordering $<$ on $\A_0(\vr_P^*)$ extending the partial ordering~$\prec$
of Subsection~\ref{curve0bl_subs} such~that
$$\vr_1\!<\!\vr_2, \quad \rho_1^*\!\in\!\A_0(\vr^*;\vr_1), \quad 
\rho_2^*\!\in\!\A_0(\vr^*;\vr_2) \qquad\Lra\qquad \rho_1^*\!<\!\rho_2^*,$$
whenever $\vr_1,\vr_2\!\in\!\A_0(\ale;d,J)$ are such that
$\vr_1,\vr_2\!\prec\!\vr^*$.
In the next subsection, we will refer to the blowup construction of
Subsection~\ref{curve0bl_subs} corresponding to such an ordering.\\

\noindent
Via the projection maps $\pi_l$, the bundle sections $\cD_0$ of Subsection~\ref{map0str_subs}
induce a linear bundle~map
$$\cD_{(\ale,J)}\!: F_{(\ale,J)}\lra \ev_0^*T\Pn$$
over $\ov\M_{0,(\ale,J)}$. 
In turn, this homomorphism induces a bundle section
$$\wt\cD_0\in\Ga\big(\wt\M_{0,(\ale,J)}^0;
  \E_0^*\!\otimes\!\pi_{\P F_{(\ale,J)}}^*\ev_0^*T\Pn\big),
\quad\hbox{where}\quad
\E_0=\ga_{(\aleph,J)}\lra \wt\M_{0,(\ale,J)}^0.$$
This section vanishes identically on the subvarieties 
$\wt\M_{0,\vr}^0$ of $\wt\M_{0,(\ale,J)}^0$ 
with $\vr\!\in\!\A_0^*(\ale;d,J)$.

\begin{lmm}
\label{map0bl_lmm2}
The section $\wt\cD_0$ is transverse to the zero set
on the complement of the subvarieties $\wt\M_{0,\vr^*}^0$ 
with  $\vr^*\!\in\!\A_0(\ale;d,J)$.
Furthermore, for every $\vr^*\!\in\!\A_0(\ale;d,J)$ as in Lemma~\ref{map0bl_lmm1},
the differential of~$\wt\cD_0$,
$$\na\wt\cD_0\!: \N_{\io_{0,\vr^*}} \lra 
\io_{0,\vr^*}^{\,*}
\big( \E_0^*\!\otimes\!\pi_{\P F_{(\ale,J)}}^*\ev_0^*T\Pn\big)
=\pi_P^*\E_0^*\!\otimes\!\pi_B^*\ev_0^*T\Pn,$$
in the normal direction to the immersion $\io_{0,\vr^*}$ is given~by
\begin{gather*}
\na\wt\cD_0\big|_{\pi_P^*L_{0,(l,i)}\otimes\pi_B^*\pi_{(l,i)}^*L_0}
=\pi_P^*s_{0,(l,i)}\!\otimes\!\pi_B^*\pi_{(l,i)}^*\cD_0
~~~\forall~ i\!\in\![m_l^*],\,l\!\in\!\ale_P(\vr^*),\\
\hbox{and}\qquad
\na\wt\cD_0\big|_{\N_{\io_{0,\vr^*}}^{\top}}
=\pi_P^*\id\!\otimes\!\pi_B^*\cD_{(\ale_B(\vr^*),J_B(\vr^*))},
\end{gather*}
where $s_{0,(l,i)}$ is the homomorphism defined in Subsection~\ref{curve0bl_subs}.
\end{lmm}

\noindent
This lemma follows immediately from Lemma~\ref{deriv0str_lmm}.

\subsection{Inductive Construction}
\label{map0blconstr_subs}

\noindent
We are now ready to describe the inductive assumptions for our construction of the blowup
$$\pi_{0,(\ale,J)}\!: 
\wt\M_{(\ale,J)}^{\vr_{\max}} \!\equiv\! \wt\M_{0,(\ale,J)}(\Pn,d)  
\lra   \wt\M_{(\ale,J)}^0 \!\equiv\! \wt\M_{0,(\ale,J)}^0(\Pn,d).$$
Suppose $\vr\!\in\!\A_0(\ale;d,J)$ and we have constructed\\
${}\quad$ ($I1$) a blowup $\pi_{\vr-1}\!: 
\wt\M_{(\ale,J)}^{\vr-1}\!\lra\!\wt\M_{(\ale,J)}^0$
such that $\pi_{\vr-1}$ is an isomorphism outside of the preimage of the spaces 
$\wt\M_{0,\vr'}^0$ with $\vr'\!\le\!\vr\!-\!1$;\\
${}\quad$ ($I2$) a line bundle $\E_{\vr-1}\!\lra\!\wt\M_{(\ale,J)}^{\vr-1}$;\\
${}\quad$ ($I3$) a section  $\wt\cD_{\vr-1} \in 
 \Ga\big( \wt\M_{(\ale,J)}^{\vr-1}; 
  \E_{\vr-1}^*\!\otimes\!\pi_{\vr-1}^{~*}\pi_{\P F_{(\ale,J)}}^*\ev_0^*T\Pn\big)$.\\

\noindent
For each $\vr^*\!>\!\vr\!-\!1$, let
$$\wt\M_{0,\vr^*}^{\vr-1}\equiv\wt\M_{0,\vr^*}^{\vr-1}(\Pn,d)
\subset \wt\M_{(\ale,J)}^{\vr-1}$$
be the proper transform of $\wt\M_{0,\vr^*}^0$ in $\wt\M_{(\ale,J)}^{\vr-1}$.
We assume that\\
${}\quad$ ($I4$) the section $\wt\cD_{\vr-1}$ is transverse to the zero set
on the complement of the subvarieties $\wt\M_{0,\vr^*}^{\vr-1}$ 
with $\vr^*\!>\!\vr\!-\!1$ and vanishes identically along these subvarieties;\\
${}\quad$ ($I5$) if $\vr_1,\vr_2\!\in\!\A_0(\ale;d,J)$ are such that
$\vr_1\!\neq\!\vr_2$, $\vr_1\!\not\prec\!\vr_2$, $\vr_2\!\not\prec\!\vr_1$,
and $\vr\!-\!1\!<\!\vr_1,\vr_2$, then
$$\wt\M_{0,\vr_1}^{\vr-1} \cap \wt\M_{0,\vr_2}^{\vr-1} ~
\begin{cases}
\subset \wt\M_{0,\ti\vr(\vr_1,\vr_2)}^{\vr-1}, &\hbox{if}~ \ti\vr(\vr_1,\vr_2)\!>\!\vr\!-\!1;\\
=\eset, &\hbox{otherwise},
\end{cases}$$
where $\ti\vr(\vr_1,\vr_2)$ is as in Lemma~\ref{map0bl_lmm1}.\\

\noindent
We also assume that for all $\vr^*\!\in\!\A_0(\ale;d,J)$ such that
$\vr^*\!>\!\vr\!-\!1$:\\
${}\quad$ ($I6$) the domain of the $G_{\vr^*}$-invariant immersion $\io_{\vr-1,\vr^*}$
induced by $\io_{0,\vr^*}$~is
\begin{gather*}
 \wt\cM_{0,\vr_P^*}^{\rho_{\vr^*}(\vr-1)} \times
\ov\M_{0,(\ale_B(\vr^*),J_B(\vr^*))},
\qquad\hbox{where}\\
\rho_{\vr^*}(\vr\!-\!1)=
\begin{cases}
\max\big\{\rho\!\in\!\A_0(\vr^*;\vr')\!: 
\vr'\!\le\!\vr\!-\!1,\vr'\!\prec\!\vr^*\big\}, 
&\begin{aligned}
&\hbox{if}~\exists \vr'\!\in\!\A_0(\ale;d,J)\\
&~\hbox{s.t.}~\vr'\!\le\!\vr\!-\!1,\vr'\!\prec\!\vr^*;
\end{aligned}\\
0,& \hbox{otherwise};
\end{cases}
\end{gather*}
${}\quad$ ($I7$) if $\vr'\!\in\!\A_0(\ale;d,J)$ is such that
$\vr\!-\!1\!<\!\vr'\!\prec\!\vr^*$, then
$$\io_{\vr-1,\vr^*}^{~-1} \big( \wt\M_{0,\vr'}^{\vr-1}\big)
=\Big(\bigcup_{\rho\in\A_0(\vr^*;\vr')} \!\!\!\!\!\!\!  
\wt\cM_{0,\rho}^{\rho_{\vr^*}(\vr-1)} \Big)
\times \ov\M_{0,(\ale_B(\vr^*),J_B(\vr^*))};$$
${}\quad$ ($I8$) $\Im^s\,\bar\io_{\vr-1,\vr^*}\subset\bigcup_{\vr-1<\vr'\prec\vr^*}
\wt\M_{0,\vr'}^{\vr-1}$, where
$$\bar\io_{\vr-1,\vr^*}\!:
\Big(\wt\cM_{0,\vr_P^*}^{\rho_{\vr^*}(\vr-1)}\!\times\!
\ov\M_{0,(\ale_B(\vr^*),J_B(\vr^*))}\Big)\big/G_{\vr^*}
\lra \wt\M_{(\ale,J)}^{\vr-1},$$
is the immersion map induced by $\io_{\vr-1,\vr^*}$.\\

\noindent
Furthermore, we assume that\\
${}\quad$ ($I9$) the collections 
$\{\io_{\vr-1,\vr^*}\}_{\vr^*\in\A_0(\ale;d,J),\vr^*>\vr-1}$  
and  $\{\bar\io_{\vr-1,\vr^*}\}_{\vr-1\in\A_0(\ale;d,J),\vr^*>\vr-1}$
of immersions are properly self-intersecting.\\
Finally, for all $\vr^*\!\in\!\A_0(\ale;d,J)$ such that $\vr^*\!>\!\vr\!-\!1$:\\
${}\quad$ ($I10$) 
$\io_{\vr-1,\vr^*}^{\,*}\E_{\vr-1}=\pi_P^*\E_{\rho_{\vr^*}(\vr-1)}$, 
where
$$\pi_P, \, \pi_B\!: \wt\cM_{0,\vr_P^*}^{\rho_{\vr^*}(\vr-1)} \times
\ov\M_{0,(\ale_B(\vr^*),J_B(\vr^*))}
\lra \wt\cM_{0,\vr_P^*}^{\rho_{\vr^*}(\vr-1)}, \,
\ov\M_{0,(\ale_B(\vr^*),J_B(\vr^*))}$$
are the two projection maps;\\
${}\quad$ ($I11$) if $\vr^*$ is as in Lemma~\ref{map0bl_lmm1}, then
the normal bundle for the immersion $\io_{\vr-1,\vr^*}$ is given~by
\begin{equation*}\begin{split}
\N_{\io_{\vr-1,\vr^*}} 
&=  \N_{\io_{\vr-1,\vr^*}}^{\bot}\!\oplus\!\N_{\io_{\vr-1,\vr^*}}^{\top} \\
&\equiv\bigoplus_{l\in\ale_P(\vr^*)}\bigoplus_{i\in[m_l^*]}
  \!\!\pi_P^*L_{\rho_{\vr^*}(\vr-1),(l,i)}\!\otimes\!\pi_B^*\pi_{(l,i)}^*L_0~\oplus~  
\pi_P^*\E_{\rho_{\vr^*}(\vr-1)}^*\!\otimes\! 
\bigoplus_{l\in\aleph_S(\vr^*)}\!\!\!\!\!\pi_l^*L_0,
\end{split}\end{equation*}
where $L_{\rho_{\vr^*}(\vr-1),(l,i)},\E_{\rho_{\vr^*}(\vr-1)}\!\lra\!
\wt\cM_{0,\vr_P^*}^{\rho_{\vr^*}(\vr-1)}$ are the line bundles constructed 
in Subsection~\ref{curve0bl_subs};\\
${}\quad$ ($I12$) the differential of~$\wt\cD_{\vr-1}$,
$$\na\wt\cD_{\vr-1}\!:  \N_{\io_{\vr-1,\vr^*}} 
\lra  \io_{\vr-1,\vr^*}^{\,*}\big(\E_{\vr-1}^{\,*}\!\otimes\!
\pi_{\vr-1}^{~*}\pi_{\P F_{(\ale,J)}}^*\ev_0^*T\Pn\big)
=\pi_P^*\E_{\rho_{\vr^*}(\vr-1)}^*\!\otimes\!\pi_B^*\ev_0^*T\Pn,$$
in the normal direction to the immersion $\io_{\vr-1,\vr^*}$ is given~by
\begin{gather*}
\na\wt\cD_{\vr-1}
 \big|_{\pi_P^*L_{\rho_{\vr^*}(\vr-1),(l,i)}\otimes\pi_B^*\pi_{(l,i)}^*L_0}
 =\pi_P^*s_{\rho_{\vr^*}(\vr-1),(l,i)}\!\otimes\!\pi_B^*\pi_{(l,i)}^*\cD_0
\quad\forall~i\!\in\![m_l^*],\,l\!\in\!\ale_P(\vr^*)\\
\hbox{and}\qquad
\na\wt\cD_{\vr-1}\big|_{\N_{\io_{\vr-1,\vr^*}}^{\top}}
=\pi_P^*\id\!\otimes\!\pi_B^*\cD_{(\ale_B(\vr^*),J_B(\vr^*))},
\end{gather*}
where $s_{\rho_{\vr^*}(\vr-1),(l,i)}$ is 
the homomorphism defined in Subsection~\ref{curve0bl_subs}.\\

\noindent
By the inductive assumption~($I4$), the loci on which the sections~$\wt\cD_{\vr}$ 
fail to be transverse to the zero set shrink and eventually disappear.
For each~$\vr$, the behavior of~$\wt\cD_{\vr}$ in the directions normal to 
the ``bad'' locus is described by~($I12$).
By the inductive assumption~($I5$), if $\vr_1$ and $\vr_2$ are non-comparable elements
of $(\A_0(\ale;d,J),\prec)$, the proper transforms of $\wt\M_{0,\vr_1}^0$ and 
$\wt\M_{0,\vr_2}^0$ become disjoint by the time either is ready to be blown up
for any ordering~$<$ extending the partial ordering~$\prec$.
Similarly to Subsections~\ref{curve1bl_subs} and~\ref{curve0bl_subs},
($I5$) will imply that the end result of the present blowup construction is 
independent of the choice of an extension~$<$.
By~($I6$), our blowup construction modifies each immersion~$\io_{0,\si^*}$
by changing the first factor of the domain according to the blowup
construction of Subsection~\ref{curve0bl_subs}, until a proper transform
of the image of~$\io_{0,\si^*}$ is to be blown~up; see below.
By~($I8$), by the time this happens the immersion~$\bar\io_{0,\si^*}$
induced by~$\io_{0,\si^*}$ transforms into an embedding.
Thus, all blowup loci are smooth.\\

\noindent
We note that all of the assumptions ($I1$)-($I12$) are satisfied if $\vr\!-\!1$ 
is replaced by~$0$.
In particular, ($I5$) is a restatement of Lemma~\ref{map0bl_lmm1a},
while~($I4$) and~($I12$) are the two parts of Lemma~\ref{map0bl_lmm2}. 
The statements ($I7$)-($I11$), with $\vr\!-\!1$ replaced by~$0$,
are contained in Lemma~\ref{map0bl_lmm1}.\\

\noindent
If $\vr\!\in\!\A_0(\ale;d,J)$ is as above, let
$$\ti\pi_{\vr}\!: \wt\M_{(\ale,J)}^{\vr} \lra \wt\M_{(\ale,J)}^{\vr-1}$$
be the blowup of $\wt\M_{(\ale,J)}^{\vr-1}$ along $\wt\M_{0,\vr}^{\vr-1}$,
which is a smooth subvariety by the inductive assumption~($I8$).
We denote the exceptional divisor for this blowup by $\wt\M_{0,\vr}^{\vr}$.
If $\vr^*\!>\!\vr$, let $\wt\M_{0,\vr^*}^{\vr}\!\subset\!\wt\M_{(\ale,J)}^{\vr}$
be the proper transform of $\wt\M_{0,\vr^*}^{\vr-1}$.
We~put 
\begin{equation}\label{bundletwist_e3}
\E_{\vr}= \ti\pi_{\vr}^*\E_{\vr-1}\otimes\O\big(\wt\M_{0,\vr}^{\vr}\big).
\end{equation}
The section $\ti\pi_{\vr}^*\wt\cD_{\vr-1}$ vanishes identically
along the divisor $\wt\M_{0,\vr}^{\vr}$. Thus, it induces a section
$$\wt\cD_{\vr}\in \Ga
\big( \wt\M_{(\ale,J)}^{\vr}; \E_{\vr}^*\!\otimes\!
  \pi_{\vr}^*\pi_{\P F_{(\ale,J)}}^*\ev_0^*T\Pn\big),$$
where $\pi_{\vr}=\pi_{\vr-1}\!\circ\ti\pi_{\vr}$.\\

\noindent
The inductive assumptions ($I1$)-($I3$), 
with $\vr\!-\!1$ replaced by $\vr$, are clearly satisfied,
while~($I5$), ($I8$), and ($I9$) follow from~(2) of Lemma~\ref{ag_lmm2a}
and Corollary~\ref{immercoll_crl}.
On the other hand, by~($I6$),  the domain of the immersion $\io_{\vr-1,\vr}$~is
$$\wt\cM_{0,\vr_P}^{\rho_{\vr}(\vr-1)} \times
\ov\M_{0,(\ale_B(\vr),J_B(\vr))}
=\wt\cM_{0,\vr_P} \times \ov\M_{0,(\ale_B(\vr),J_B(\vr))},$$
where $\wt\cM_{0,\vr_P}\!\lra\!\wt\cM_{0,\vr_P}^0$ is the blowup constructed
in Subsection~\ref{curve0bl_subs}.
By~($I11$), the normal bundle for the immersion $\io_{\vr-1,\vr}$ is given~by
\begin{equation*}\begin{split}
\N_{\io_{\vr-1,\vr}}  &= \bigoplus_{l\in\ale_P(\vr)}\bigoplus_{i\in[m_l]} \!\pi_P^*L_{\rho_{\vr}(\vr-1),(l,i)}\!\otimes\!\pi_B^*\pi_{(l,i)}^*L_0~\oplus~  
\pi_P^*\E_{\rho_{\vr}(\vr-1)}^*\!\otimes\! 
\bigoplus_{l\in\aleph_S(\vr)}\!\!\!\!\pi_l^*L_0\\
&= \bigoplus_{l\in\ale_P(\vr)}\bigoplus_{i\in[m_l]}
\!\pi_P^*\L\!\otimes\!\pi_B^*\pi_{(l,i)}^*L_0~\oplus~  
\pi_P^*\L\!\otimes\! 
\bigoplus_{l\in\aleph_S(\vr)}\!\!\!\!\pi_l^*L_0\\
&=\pi_P^*\L \otimes  \pi_B^*F_{(\ale_B(\vr),J_B(\vr))},
\end{split}\end{equation*} 
where $\L\!\lra\!\wt\cM_{0,\vr_P}$ is the universal tangent line bundle 
constructed in Subsection~\ref{curve0bl_subs}.
We also note~that by~($I10$),
$$\io_{\vr-1,\vr}^{\,*}
\big(\E_{\vr-1}^{\,*}\!\otimes\!\pi_{\P F_{(\ale,J)}}^*\ev_0^*T\Pn\big)
=\pi_P^*\L \otimes  \pi_B^*\pi_{\P F_{(\ale_B(\vr),J_B(\vr))}}^*\ev_0^*T\Pn.$$
By~($I12$), the differential of~$\wt\cD_{\vr-1}$ 
in the normal direction to the immersion $\io_{\vr-1,\vr}$ is given~by
$$ \na\wt\cD_{\vr-1} =\pi_P^*\id\otimes \pi_B^*\cD_{(\ale_B(\vr),J_B(\vr))}.$$
Thus, if
$$ \io_{\vr,\vr}\!: \P\N_{\io_{\vr-1,\vr}}
\approx \wt\cM_{0,\vr_P} \times \wt\M_{0,(\ale_P(\vr),J_B(\vr))}^0
\lra \wt\M_{0,\vr}^{\vr}\subset \wt\M_{(\ale,J)}^{\vr}$$
is the immersion induced by $\io_{\vr-1,\vr}$, then 
\begin{equation*}\begin{split}
\io_{\vr,\vr}^{\,*}\wt\cD_{\vr} \!=\! \pi_B^*\wt\cD_0 
&\in \Ga\big(\wt\cM_{0,\vr_P}\!\times\!\wt\M_{0,(\ale_P(\vr),J_B(\vr))}^0;
\io_{\vr,\vr}^{~*}\big(\E_{\vr}^* \!\otimes\!\pi_{\vr}^*\pi_{\P F_{(\ale,J)}}^*\ev_0^*T\Pn
\big)\big)\\
&=\Ga\big(\wt\cM_{0,\vr_P}\!\times\!\wt\M_{0,(\ale_B(\vr),J_B(\vr))}^0;
\pi_B^*(\ga_{(\ale_B(\vr),J_B(\vr))}^*\!\otimes\!
\pi_{\P F_{(\ale_B(\vr),J_B(\vr))}}^*\ev_0^*T\Pn)\big).
\end{split}\end{equation*}
Lemmas~\ref{map0bl_lmm1} and~\ref{map0bl_lmm2} thus imply that the restriction of the section 
$\wt\cD_{\vr}$ to the exceptional divisor $\wt\M_{0,\vr}^{\vr}$
is transverse to the zero set away from the subvarieties
$\wt\M_{0,\vr^*}^{\vr}$ with $\vr^*\!>\!\vr$.
Thus, by the inductive assumption~($I4$) as stated above, 
($I4$) is satisfied with $\vr\!-\!1$ replaced by~$\vr$.\\

\noindent
We now verify that the remaining inductive assumptions are satisfied.
If $\vr\!<\!\vr^*$, but $\vr\!\not\prec\!\vr^*$, 
$$\rho_{\vr^*}(\vr) = \rho_{\vr^*}(\vr\!-\!1)
\qquad\hbox{and}\qquad
\ov\M_{0,\vr^*}^{\vr-1}\cap\ov\M_{0,\vr}^{\vr-1} =\eset,$$
by definition and by~($I5$), respectively.
It then follows that
\begin{gather*}
\io_{\vr,\vr^*}=\io_{\vr-1,\vr^*}, \qquad
\wt\M_{0,\vr^*}^{\vr}\!\cap\!\wt\M_{0,\vr'}^{\vr}
=\wt\M_{0,\vr^*}^{\vr-1}\!\cap\!\wt\M_{0,\vr'}^{\vr-1} \quad\forall \vr'\!>\!\vr, \\ 
\io_{\vr,\vr^*}^{\,*}\E_{\vr}=\io_{\vr-1,\vr^*}^{\,*}\E_{\vr-1},
\qquad \N_{\io_{\vr,\vr^*}}=\N_{\io_{\vr-1,\vr^*}},
\quad\hbox{and}\quad  \na\wt\cD_{\vr}=\na\wt\cD_{\vr-1}.
\end{gather*}
Thus, the inductive assumptions ($I6$), ($I7$), and ($I10$)-($I12$), as stated above,
imply the corresponding statements with $\vr\!-\!1$ replaced by~$\vr$.\\

\noindent
Suppose that $\vr\!\prec\!\vr^*$.
By~($I6$) and~(1) of Lemma~\ref{ag_lmm2a}, the domain of the immersion~$\io_{\vr,\vr^*}$ 
induced by  the immersion~$\io_{\vr-1,\vr^*}$ is the blowup~of
$$\wt\cM_{0,\vr_P^*}^{\rho_{\vr^*}(\vr-1)} \times
\ov\M_{0,(\ale_B(\vr^*),J_B(\vr^*))}$$
along the preimage of $\ov\M_{0,\vr}^{\vr-1}$ under~$\io_{\vr-1,\vr^*}$~in
$$\pi_{\rho_{\vr^*}(\vr-1)}\!\times\!\id\!:
\wt\cM_{0,\vr_P^*}^{\rho_{\vr^*}(\vr-1)} \!\times\!
\ov\M_{0,(\ale_B(\vr^*),J_B(\vr^*))}  \lra
\wt\cM_{0,\vr_P^*}^0 \!\times\! \ov\M_{0,(\ale_B(\vr^*),J_B(\vr^*))}.$$
By ($I7$), this preimage is
$$\Big(\bigcup_{\rho\in\A_0(\vr^*;\vr)} \!\!\!\!\!\!\! 
\wt\cM_{0,\rho}^{\rho_{\vr^*}(\vr-1)}\Big)
\times \ov\M_{0,(\ale_B(\vr^*),J_B(\vr^*))}.$$
By the last paragraph of Subsection~\ref{curve0bl_subs} and
the second paragraph after Lemma~\ref{map0bl_lmm1}, 
$$\wt\cM_{0,\rho_1}^{\rho_{\vr^*}(\vr-1)}
\cap \wt\cM_{0,\rho_2}^{\rho_{\vr^*}(\vr-1)} = \eset
\qquad \forall\, \rho_1,\rho_2\!\in\!\A_0(\vr^*;\vr), \, \rho_1\!\neq\!\rho_2.$$
Thus, by the construction of Subsection~\ref{curve0bl_subs},
the blowup of  $\wt\cM_{0,\vr_P^*}^{\rho_{\vr^*}(\vr-1)}$
along 
$$\bigcup_{\rho\in\A_0(\vr^*;\vr)}\!\!\!\!\!\!\wt\cM_{0,\rho}^{\rho_{\vr^*}(\vr-1)}$$
is $\wt\cM_{0,\vr_P^*}^{\rho_{\vr^*}(\vr)}$, as needed for
the inductive statement~($I6$), with $\vr\!-\!1$ replaced by~$\vr$.
For the same reasons, ($I10$),  \e_ref{bundletwist_e2b}, and~\e_ref{bundletwist_e3}
imply~that
\begin{equation*}\begin{split}
\io_{\vr,\vr^*}^{\,*}\E_{\vr}
&=\io_{\vr-1,\vr^*}^{\,*}\E_{\vr-1}
\otimes \io_{\vr,\vr^*}^{\,*}\O\big(\wt\M_{0,\vr}^{\vr}\big)\\
&=\pi_P^* \E_{\rho_{\vr^*}(\vr-1)}
\otimes \bigotimes_{\rho\in\A_0(\vr^*;\vr)} \!\!\!\!\!\!\!
\pi_P^* \O\big(\wt\cM_{0,\rho}^{\rho}\big)
=\pi_P^* \E_{\rho_{\vr^*}(\vr)}.
\end{split}\end{equation*}
Thus, the inductive statement~($I10$), with $\vr\!-\!1$ replaced by~$\vr$,
is satisfied. The assumption~($I7$) is checked similarly, using~(3) of 
Lemma~\ref{ag_lmm2a}.\\

\noindent
We next determine the normal bundle for the immersion $\io_{\vr,\vr^*}$.
The restrictions of the line bundles  $L_{\rho_{\vr^*}(\vr-1),(l,i)}$ 
and $\E_{\rho_{\vr^*}(\vr-1)}$ to the complement of the exceptional divisors in
$\wt\cM_{0,\vr_P^*}^{\rho_{\vr^*}(\vr-1)}$ are $\pi_{\rho_{\vr^*}(\vr-1)}^*L_{0,(l,i)}$ and 
$\pi_{\rho_{\vr^*}(\vr-1)}^*\E_0$, by the construction of Subsection~\ref{curve0bl_subs}.
Thus, by the last statement of Lemma~\ref{map0bl_lmm1}, ($I11$), 
and the inductive assumptions~($I1$) above and in Subsection~\ref{curve0bl_subs},
\begin{equation*}\begin{split}
&\io_{\vr-1,\vr^*}\big|_{\wt\cM_{0,\rho}^{\rho_{\vr^*}(\vr-1)}
                          \times\ov\M_{0,(\ale_B(\vr^*),J_B(\vr^*))}}^* 
T\ov\M_{0,\vr}^{\vr-1}
\big/T\big(\wt\cM_{0,\rho}^{\rho_{\vr^*}(\vr-1)}
             \!\times\!\ov\M_{0,(\ale_B(\vr^*),J_B(\vr^*))}\big)\\ 
 &~=\!\bigoplus_{l\in\ale_P(\vr^*)-\ale_P(\vr)} \bigoplus_{i\in[m_l^*]}
\!\!\pi_P^*L_{\rho_{\vr^*}(\vr-1),(l,i)}\!\otimes\!\pi_B^*\pi_{(l,i)}^*L_0
~\oplus~\!\bigoplus_{l\in\ale_P(\vr)} \bigoplus_{i\in[m_l^*]-I_{l,P}}
\!\!\!\!\!\!\!\!\!\pi_P^*L_{\rho_{\vr^*}(\vr-1),(l,i)}\!\otimes\!\pi_B^*\pi_{(l,i)}^*L_0
\end{split}\end{equation*}\\
for all $\rho\!\in\!\A_0(\vr^*;\vr)$ as in the statement of Lemma~\ref{map0bl_lmm1}.
Let
$$I_P(\rho)=\big\{(l,i)\!:l\!\in\!\ale_P(\vr),\, i\!\in\!I_{l,P}\big\}.$$
From Lemma~\ref{ag_lmm2b}, we then obtain
\begin{equation*}\begin{split}
\N_{\io_{\vr,\vr^*}} 
&=\bigoplus_{l\in\ale_P(\vr^*)}\bigoplus_{i\in[m_l^*]} \!\!
\bigg(\big(\pi_P^*L_{\rho_{\vr^*}(\vr-1),(l,i)}\!\otimes\!\pi_B^*\pi_{(l,i)}^*L_0\big)
  \!\otimes\!\pi_P^*
  \O\Big(~~~~~-\!\!\!\!\!\!\!\!\!\!\!\!\!\!\!\!\!\!\!
  \sum_{\rho\in\A_0(\vr^*;\vr),(l,i)\in I_P(\rho)}\!\!\!\!\!\!\!\!\!\!\!\!\!\!\!\!\!\!
   \wt\cM_{0,\rho}^{\rho_{\vr^*}(\vr-1)}\Big)\bigg)\\
&\qquad\qquad\qquad\qquad\qquad\oplus~  
\pi_P^*\E_{\rho_{\vr^*}(\vr-1)}^*\!\otimes\! 
\bigoplus_{l\in\aleph_S(\vr^*)}\!\!\!\!\!\pi_l^*L_0\!\otimes\!\pi_P^*
  \O\Big(-\!\!\!\!\!\!\!\!
  \sum_{\rho\in\A_0(\vr^*;\vr)}\!\!\!\!\!\!\!\wt\cM_{0,\rho}^{\rho_{\vr^*}(\vr-1)}\Big)\bigg)\\
&=\bigoplus_{l\in\ale_P(\vr^*)}\bigoplus_{i\in[m_l^*]} \!\!
\big(\pi_P^*L_{\rho_{\vr^*}(\vr),(l,i)}\!\otimes\!\pi_B^*\pi_{(l,i)}^*L_0\big)
~\oplus~\pi_P^*\E_{\rho_{\vr^*}(\vr)}^*\!\otimes\! 
\bigoplus_{l\in\aleph_S(\vr^*)}\!\!\!\!\!\pi_l^*L_0.
\end{split}\end{equation*}
The last equality above follows from \e_ref{bundletwist_e2a} by the same argument
as in the previous paragraph.
We have thus verified that the inductive assumption~($I11$), with $\vr\!-\!1$ replaced
by $\vr$, is satisfied.
Finally, the inductive assumption~($I12$) and the continuity of the two bundle
sections involved in the identity in~($I12$), with $\vr\!-\!1$ replaced by~$\vr$,
imply ($I12$) with $\vr\!-\!1$ replaced by~$\vr$.\\

\noindent
We conclude this construction after the blowup at the $\vr_{\max}$ step. Let
$$\wt\M_{0,(\ale,J)}(\Pn,d)=\wt\M_{(\ale,J)}^{\vr_{\max}}, \qquad
\ti\E=\E_{\vr_{\max}}, \qquad \wt\cD_{(\ale,J)}=\wt\cD_{\vr_{\max}}.$$
By the inductive assumption~($I4$), applied with $\vr\!-\!1$ replaced by $\vr_{\max}$,
the section $\wt\cD_{(\ale,J)}$ is transverse to the zero set.
As in the previous two subsections, the final result of this blowup construction
is independent of the order~$<$ chosen to extend the partial ordering~$\prec$
on~$\A_0^*(\ale;d,J)$, as can be seen from~($I5$).

\section{A Blowup of a Moduli Space of Genus-One Maps}
\label{map1bl_sec}

\subsection{Idealized Blowups and Immersions}
\label{map1prelim_subs1}

\noindent
In this section we describe the main blowup construction of this paper.
This is the sequential idealized blowup construction for $\ov\M_{1,k}(\Pn,d)$
with the initial data and the inductive step specified in 
Subsections~\ref{map1prelim_subs} and~\ref{map1blconstr_subs}, respectively.
This construction is outlined in Subsections~\ref{descr_subs} and~\ref{outline_subs}.\\

\noindent
In contrast to the situations in Sections~\ref{curvebl_sec} and~\ref{map0bl_sec},
the variety $\ov\M_{1,k}(\Pn,d)$ is singular.
In order to describe the structure of  $\ov\M_{1,k}(\Pn,d)$,
we introduce the notion of {\it idealized normal bundle} for an immersion.
Recall that the domain of an immersion
is assumed to be a smooth variety.

\begin{dfn}
\label{virtvar_e}
Suppose $\ov\M$ is a variety and $\io_X\!:X\!\lra\!M$ is an immersion.
An {\tt idealized normal bundle} for the immersion $\io_X$ is a vector bundle 
$\N_{\io_X}^{\ide}$ over $X$ such that $\N_{\io_X}\!\subset\!\N_{\io_X}^{\ide}$.
\end{dfn}

\noindent
{\it Remark:} An idealized normal bundle is of course not unique;
an idealized normal bundle plus any other vector bundle is still an 
idealized normal bundle. 
If the image of $\io_X$ is an irreducible component of~$\ov\M$,
an idealized normal bundle of the smallest possible rank still need not be unique;
it can be twisted by any divisor in $X$ disjoint from the preimage under~$\io_X$
of the other components of~$\ov\M$.
For each of the immersions we encounter in the next subsection, 
there is a natural choice for~$\N_{\io_X}^{\ide}$.
These idealized normal bundles also transform in a natural way under blowups 
and proper immersions, as described in Lemma~\ref{virimmer_lmm2} below.\\

\noindent
Suppose $\ov\M$ is a variety, $Z$ is a smooth subvariety of $\ov\M$,
and $\N_{\io_Z}^{\ide}$ is an idealized normal bundle for the embedding $\io_Z$ 
of $Z$ into~$\ov\M$.
Let 
$$\cE_Z\!\equiv\!\P\N_{\io_Z} \subset \Bl_Z\ov\M$$
be the exceptional divisor for the blowup of $\ov\M$ along $Z$.
We denote by $\Bl_Z^{\ide}\ov\M$ the variety obtained by identifying $\Bl_Z\ov\M$ 
with 
$$\cE_Z^{\ide}\equiv\P\N_{\io_Z}^{\ide}$$ 
along~$\cE_Z$. We will call 
$$\pi^{\ide}\!:  \Bl_Z^{\ide}\ov\M\lra \ov\M \qquad\hbox{and}\qquad
\cE_Z^{\ide}\subset\Bl_Z^{\ide}\ov\M$$
{\tt the idealized blowup of $\ov\M$ along~$Z$} and 
{\tt the idealized exceptional divisor for~$\pi^{\ide}$}, respectively.
More generally, we will call 
$$\pi\!: \wt\M \lra \ov\M$$
an {\tt idealized blowup of $\ov\M$} if $\pi$ is a composition
of idealized blowups along smooth subvarieties.
In practice, {\it idealized blowup} is simply a convenient term.
In the inductive assumption~($I1$) in Subsection~\ref{map1blconstr_subs}
below, it can be replaced by {\it morphism of varieties}, 
as the remaining inductive assumptions describe all the relevant properties 
of this morphism. Let 
$$\ga_Z\lra \cE_Z^{\ide}$$
be the tautological line bundle. Note that the normal bundle of 
$\cE_Z\!\subset\!\cE_Z^{\ide}$ in
$$\Pr_Z\ov\M=\Bl_Z\ov\M$$
is $\ga_Z|_{\cE_Z}$. 
(This observation implies the first statement of Lemma~\ref{virimmer_lmm2}.)\\

\noindent
Our strategy is as follows.  
We begin with a space with a properly self-intersecting collection of immersions, 
each with an idealized normal bundle.  
These are the immersions~$\io_{\si}$ with $\si\!\in\!\A_1(d,k)$ defined in Subsection~\ref{descr_subs};
their images are the subvarieties $\ov\M_{1,\si}(\Pn,d)$ of $\ov\M_{1,k}(\Pn,d)$.
The idealized normal bundle for the immersion~$\io_{\si}$ is the direct sum of 
the deformation spaces of the nodes between the contracted genus-one curve and 
the non-contracted genus-zero curves that are identified by~$\io_{\si}$.
At each stage, one of our immersions is an embedding, and we blow
it up, replacing it with its idealized exceptional divisor.  The
exceptional divisor of the blowup of the main component is the
intersection of the new main component with the idealized exceptional divisor. 
Then after each step, we have a new properly
self-intersecting collection of immersions.  Moreover, there
is a natural idealized normal bundle to each of the proper transforms
of the immersions we have ``yet to blow up''.\\

\noindent
We now say this more explicitly. 
The following two lemmas are direct extensions of
Corollary~\ref{immercoll_crl} and Lemma~\ref{ag_lmm2b}. 
The first lemma states that if
we have a properly self-intersecting collection of immersions, one of
which is an embedding, then upon blowing up the embedding, we still
have a properly self-intersecting collection of immersions. It is
immediate from the definition of ``properly self-intersecting'', by
checking in local coordinates.\\

\noindent
The second part of the second lemma follows from Lemma~\ref{ag_lmm2b} with only
one change.
Instead of writing 
$$\N_{\io_X} =\bigoplus_{i\in I}L_i \qquad\hbox{and}\qquad
\N_{\Pr_Z\io_X}=\bigoplus_{i\in I}
\Big(\pi^*L_i\otimes \bigotimes_{i\in I_{\vr}}\O(-E_{\vr})\Big)$$
as in the statement of Lemma~\ref{ag_lmm2b}, 
we are saying that if we have a natural inclusion
$\N_{\io_X}\!\subset\!\bigoplus_{i \in I}L_i$, then we get a natural conclusion
$$\N_{\Pr_Z\io_X} \subset \bigoplus_{i\in I}
\Big(\pi^*L_i\otimes \bigotimes_{i\in I_{\vr}}\O(-E_{\vr})\Big).$$
The vector bundles on the right are the original idealized normal
bundle and the new idealized normal bundle, respectively.

\begin{lmm}
\label{virimmer_lmm}
Suppose $\ov\M$ is a variety, 
$\{\io_{\si}\!:X_{\si}\!\lra\!\ov\M\}_{\si\in\A}$
is a properly self-intersecting collection of immersions, and
$\si\!\in\!\A$ is such that $\io_{\si}$ is an embedding.
If $\N_{\io_{\si}}^{\ide}$ is an idealized normal bundle for~$\io_{\si}$, then
$$\big\{\Pr_{\Im\,\io_{\si}}\io_{\si'}\big\}_{\si'\in\A-\{\si\}}
\cup\big\{\io_{\cE_{\Im\io_{\si}}^{\ide}}\big\}$$
is a properly self-intersecting collection of immersions 
into $\Bl_{\Im\,\io_{\si}}^{\ide}\ov\M$. 
\end{lmm}

\begin{lmm}
\label{virimmer_lmm2}
If $\ov\M$ is a variety, $Z$ is a smooth subvariety of $\ov\M$, 
and $\N_{\io_Z}^{\ide}$ is an idealized normal bundle for~$\io_Z$, then
$$\N_{\io_{\cE_Z^{\ide}}}^{\ide}=\ga_Z$$
is an idealized normal bundle for the immersion $\io_{\cE_Z^{\ide}}$.
Suppose in addition that $\io_X$, $\A$, $Z_{\vr}$, and $E_{\vr}$ are as 
in Lemma~\ref{ag_lmm2b} and 
$\N_{\io_X}^{\ide}$ is an idealized normal bundle for~$\io_X$.
If there exist a splitting
$$\N_{\io_X}^{\ide}=\bigoplus_{i\in I}L_i \lra X$$
and a subset $I_{\vr}$ of $I$ for each $\vr\!\in\!\A$ such that \e_ref{ag_lmm2b_e} holds,
then 
$$\N_{\Pr_Z\io_X}^{\ide}=\bigoplus_{i\in I}
\Big(\pi^*L_i\otimes \bigotimes_{i\in I_{\vr}}\O(-E_{\vr})\Big)$$
is an idealized normal bundle for the immersion $\Pr_Z\io_X$.
\end{lmm}

\begin{dfn}
\label{prsub_dfn}
Suppose $\ov\M$ is a variety, $\io_X\!:X\!\lra\!\ov\M$ is an immersion,
$\ov\M^0$ is a subvariety in $\ov\M$, and $TC\ov\M^0\!\subset\!T\ov\M$
is the tangent cone of $\ov\M^0$ in~$\ov\M$ ($TC\ov\M^0$ not necessarily reduced).
The subvariety $\ov\M^0$ is {\tt proper relative to~$\io_X$} if 
$$d\io_X\,TC\io_X^{-1}(\ov\M^0)
= \io_X^*TC\ov\M^0\cap \Im\, d\io_X\subset \io_X^*T\ov\M$$
and the~map
\begin{equation}\label{prsub_dfn_e}
\io_X^*TC\ov\M^0|_{\io_X^{-1}(\ov\M^0)}\big/\Im\,d\io_X|_{TC\io_X^{-1}(\ov\M^0)}
\lra \io_X^*T\ov\M/\Im\,d\io_X \subset \N_{\io_X}^{\ide}
\end{equation}
induced by inclusions is injective, with its image being reduced.
\end{dfn}

\noindent
The left-hand side of~\e_ref{prsub_dfn_e} denotes the family of cones over 
$\io_X^{-1}(\ov\M^0)$ such that for each $x\!\in\!\io_X^{-1}(\ov\M^0)$
$$\io_X^*TC\ov\M^0|_{\io_X^{-1}(\ov\M^0)}\big/\Im\,d\io_X|_{TC\io_X^{-1}(\ov\M^0)}
\Big|_x$$
is the quotient by the minimal vector subspace of $\Im\,d\io_X|_x\!=\!d\io_X(T_xX)$ 
containing the cone $\Im\,d\io_X|_{T_xC\io_X^{-1}(\ov\M^0)}$.
If $TC\io_X^{-1}(\ov\M^0)$ is a vector bundle, the two conditions 
in~(2) of Definition~\ref{prsub_dfn} are equivalent.\\

\noindent
If $\ov\M^0$ is a subvariety of $\ov\M$ which is proper relative
to an immersion $\io_X\!:X\!\lra\!\ov\M$,
we denote~by 
$$\N_{\io_X|\ov\M^0}\subset \io_X^*T\ov\M/\Im\,d\io_X
\subset \N_{\io_X}^{\ide}$$
the image of the homomorphism~\e_ref{prsub_dfn_e}.
We will call $\N_{\io_X|\ov\M^0}$ {\tt the normal cone of 
$\io_X|_{\io_X^{-1}(\ov\M^0)}$ in~$\ov\M^0$}.

\begin{lmm}
\label{virimmer_lmm3}
Suppose $\ov\M$ is a variety, $\io_X\!:X\!\lra\!\ov\M$ is an immersion
with an idealized normal bundle~$\N_{\io_X}^{\ide}$,
$\ov\M^0$ is a subvariety of $\ov\M$ which is proper relative to~$\io_X$, 
and 
$$\cZ\subset\ov\cZ\!\equiv\!\io_X^{-1}(\ov\M^0)$$
is such that $\N_{\io_X|\ov\M^0}$ is the closure of $\N_{\io_X|\ov\M^0}|_{\cZ}$ 
in $\N_{\io_X}^{\ide}$.\\
(1) If $X$ is a smooth subvariety of $\ov\M$, then
$\Pr_X\ov\M^0$ is proper relative to the immersion~$\io_{\cE_X^{\ide}}$,
$$\cE_X^{\ide}\cap\Pr_X\ov\M^0 \subset \cE_X$$
is the closure of $\P\N_{\io_X|\ov\M^0}|_{\cZ}$ in $\cE_X^{\ide}$, and 
$$\N_{\io_{\cE_X^{\ide}}|\Pr_X\ov\M^0}=\ga_X|_{\cE_X^{\ide}\cap\Pr_X\ov\M^0}.$$
(2) If $Z$ is a smooth subvariety of $\ov\M$ disjoint from $\io_X(\cZ)$
and $\N_{\io_Z}^{\ide}$ is an idealized normal bundle for~$\io_Z$, then
$\Pr_Z\ov\M^0$ is a proper subvariety of $\Bl_Z^{\ide}\ov\M$ relative 
to the immersion~$\Pr_Z\io_X$ and
$\N_{\Pr_Z\io_X|\Pr_Z\ov\M^0}$ is the closure of $\N_{\io_X|\ov\M^0}|_{\cZ}$ 
in $\N_{\Pr_Z\io_X}^{\ide}$.
\end{lmm}
 
\noindent
The first part of (1) essentially follows from the universal property of blowing
up: if $\ov \M$ is blown up along~$Z$, then the proper transform of $\ov \M^0$ in
$\ov \M$ (the scheme-theoretic closure of $\ov\M^0\!-\!Z$ in the
blowup) is the blowup of $\ov \M^0$ along $\ov\M^0\!\cap\!Z$,
and the normal bundle to the exceptional divisor in $Bl_{\ov \M^0 \cap Z} \ov \M^0$
is the restriction of the normal bundle of the exceptional
divisor in $Bl_Z \ov \M$.  The statement (1) is the etale-local version of this.
Part~(2) is clear by working in local coordinates.

\subsection{Preliminaries}
\label{map1prelim_subs}

\noindent
In this subsection, we state a number of known facts concerning 
the moduli space $\ov\M_{1,k}(\Pn,d)$ that insure that the inductive requirements
of the next subsection are satisfied at the initial stage of the inductive construction.
Lemmas~\ref{map1bl_lmm1a}-\ref{map1bl_lmm1b3}, with the exception of one statement, 
are straightforward to check from the definitions (and~\cite{P} in some cases).
We show that the last statement of Lemma~\ref{map1bl_lmm1b1} is a reinterpretation
of a standard fact concerning moduli spaces of stable maps.\\

\noindent
Let $(\A_1(d,k),\prec)$ be the partially ordered set of triples described
in Subsection~\ref{descr_subs}.
It has a unique minimal element and a unique maximal element:
$$\si_{\min}=(1;\eset,[k])  \qquad\hbox{and}\qquad \si_{\max}=(d;[k],\eset).$$
Let $<$ be an order on $\A_1(d,k)$ extending the partial ordering $\prec$.
For every $\si\!\in\!A_1(d,k)$, we define 
$$\si\!-\!1 \in \{0\}\!\sqcup\!\A_1(d,k)$$
as in~\e_ref{minusdfn_e}.
For each element $\si\!=\!(m;J_P,J_B)$ of $\A_1(d,k)$, let
$$\ov\M_{1,\si}^0 \!\equiv\! \ov\M_{1,\si}(\Pn,d)
\subset \ov\M_{1,k}^0 \!\equiv\! \ov\M_{1,k}(\Pn,d)$$
be the subvarieties defined in Subsection~\ref{descr_subs}.\\

\noindent
{\it Warning:} Note that $\ov\M_{1,k}^0$ denotes the entire moduli space
$\ov\M_{1,k}(\Pn,d)$ and not the main component $\ov\M_{1,k}^0(\Pn,d)$.
Similarly to Sections~\ref{curvebl_sec} and~\ref{map0bl_sec}, 
the superscript $0$ indicates the $0$th stage in the blowup process.

\begin{lmm}
\label{map1bl_lmm1a}
If $\si_1\!=\!(m_1;J_{1;P},J_{1;B})$ and $\si_2\!=\!(m_2;J_{2;P},J_{2;B})$
are elements of $\A_1(d,k)$, $\si_1\!\neq\!\si_2$,
$\si_1\!\not\prec\!\si_2$, and $\si_2\!\not\prec\!\si_1$, then
\begin{gather*}
\ov\M_{1,\si_1}^0 \cap \ov\M_{1,\si_2}^0\subset \ov\M_{1,\ti\si(\si_1,\si_2)}^0,
\qquad\hbox{where}\\
\ti\si(\si_1,\si_2)=\big(\!\min(m_1,m_2);J_{1;P}\!\cap\!J_{2;P},J_{1;B}\!\cup\!J_{2;B}\big).
\end{gather*}\\
\end{lmm}

\noindent
With $\si$ as above, we define
$$I_P(\si)=\ale_B(\si)=[m], \qquad J_P(\si)=J_P, \qquad J_B(\si)=J_B,
\qquad G_{\si}=S_m.$$
As in Subsection~\ref{descr_subs}, we denote by
\begin{gather*}
\io_{0,\si}\!: \ov\cM_{1,(I_P(\si),J_P(\si))}^0 \times 
\ov\M_{0,(\ale_B(\si),J_B(\si))} \lra 
\ov\M_{1,\si}^0\subset\ov\M_{1,k}^0,\\
\hbox{where}\qquad \ov\M_{0,(\ale_B(\si),J_B(\si))}=\ov\M_{0,(\ale_B(\si),J_B(\si))}(\Pn,d),
\end{gather*}
the natural node-identifying map and by
$$\bar\io_{0,\si}\!: 
\big(\ov\cM_{1,(I_P(\si),J_P(\si))}^0 \!\times\! 
\ov\M_{0,(\ale_B(\si),J_B(\si))}\big)\big/G_{\si}\lra \ov\M_{1,k}^0$$
the induced immersion.
Let 
$$\pi_P, \, \pi_B\!: 
\ov\cM_{1,(I_P(\si),J_P(\si))}^0 \!\times\!  
\ov\M_{0,(\ale_B(\si),J_B(\si))}\lra
\ov\cM_{1,(I_P(\si),J_P(\si))}^0,  \ov\M_{0,(\ale_P(\si),J_B(\si))}$$
be the two projection maps.

\begin{lmm}
\label{map1bl_lmm1b1}
If $d,n\!\in\!\Z^+$ and $k\!\in\!\bar\Z^+$, the collections
$\{\io_{0,\si}\}_{\si\in\A_1(d,k)}$ and $\{\bar\io_{0,\si}\}_{\si\in\A_1(d,k)}$ 
of immersions are properly self-intersecting.
If $\si^*\!=\!(m^*;J_P^*,J_B^*)\!\in\!\A_1(d,k)$,
$$\Im^s\,\bar\io_{0,\si^*} \subset \bigcup_{\si'\prec\si^*}\! \ov\M_{1,\si'}
\qquad\hbox{and}\qquad 
\N_{\io_{0,\si^*}}^{\ide}=\bigoplus_{i\in[m^*]}\!\! \pi_P^*L_i\!\otimes\!\pi_B^*\pi_i^*L_0$$
is an idealized normal bundle for $\io_{0,\si^*}$.
\end{lmm}

\noindent
We deduce the last claim of this lemma from the deformation-obstruction exact 
sequence~(24.2) in~\cite{H} as follows.
Suppose
\begin{gather*}
[\Si,u]=\io_{0,\si^*}\big([\Si_P]\!\times\![\Si_B,u_B]\big)
\in \ov\M_{1,\si^*}^0, \qquad\hbox{where}\\
[\Si_B,u_B]=\big([\Si_i,u_i]\big)_{i\in[m^*]}\in\ov\M_{0,(\ale_B(\si^*),J_B(\si^*))}.
\end{gather*}
By \cite[(24.2)]{H}, there exists a natural homomorphism
$$j_{\Si,u}\!: T\ov\M_{1,k}(\Pn,d)\big|_{[\Si,u]}=\Def(\Si,u) \lra \Def(\Si),$$
where $\Def(\Si,u)$ and $\Def(\Si)$ denote the deformations of 
the stable-map pair $(\Si,u)$ and the deformations of 
the curve~$\Si$ (with its marked points), respectively.
As $[\Si,u]$ is considered as the image of $[\Si_P]\!\times\![\Si_B,u_B]$
under $\io_{0,\si^*}$, there are $m^*$ distinguished nodes of~$\Si$.
These are the nodes of $\Si$ that do not correspond to either the nodes of
$\Si_P$ or the nodes of any of the curves $\Si_i$ with $i\!\in\![m^*]$;
see Figure~\ref{idebund_fig}.
Let 
$$\Def(\Si_P,\Si_B)\subset \Def(\Si)$$
be the deformations of $\Si$ that do not smooth out the distinguished nodes of~$\Si$.
Since the smoothing of a given node of $\Si$ is parametrized by 
the tensor product of the tangent lines to the two branches of $\Si$ at the node,
we have an exact sequence
$$0 \lra \Def(\Si_P,\Si_B)\lra \Def(\Si) \stackrel{j_{\Si}}{\lra} 
\N_{\io_{0,\si^*}}^{\ide}\big|_{[\Si,u]} \lra 0.$$
We denote by 
$$\Def\big(\Si_P,(\Si_B,u_B)\big)\subset
T\ov\M_{1,k}(\Pn,d)\big|_{[\Si,u]}=\Def(\Si,u)$$
the kernel of the map
$$j_{\Si}\!\circ j_{\Si,u}\!: \Def(\Si,u) \lra \N_{\io_{0,\si^*}}^{\ide}\big|_{[\Si_P]\times[\Si_B,u_B]}.$$
The space $\Def\big(\Si_P,(\Si_B,u_B)\big)$ consists of deformations of $(\Si,u)$
that do not smooth out the $m^*$ distinguished nodes of~$\Si$.
Thus,
\begin{equation*}\begin{split}
\Def\big(\Si_P,(\Si_B,u_B)\big) &\approx\Def(\Si_P)\oplus\Def(\Si_B,u_B)\\
&=T\ov\cM_{1,(I_P(\si^*),J_P(\si^*))}^0|_{[\Si_P]} \oplus  
T\ov\M_{0,(\ale_B(\si^*),J_B(\si^*))}|_{[\Si_B,u_B]}.
\end{split}\end{equation*}
The isomorphism from the right-hand side to the left-hand side is 
given by~$d\io_{0,\si^*}$.
Thus, the homomorphism  $j_{\Si}\!\circ\!j_{\Si,u}$ induces an injection
$$\N_{\io_{0,\si^*}}|_{[\Si,u]} \equiv 
TC\ov\M_{1,k}(\Pn,d)\big|_{[\Si,u]}\big/\Im\,d\io_{0,\si^*}
\lra \N_{\io_{0,\si^*}}^{\ide}\big|_{[\Si_P]\times[\Si_B,u_B]},$$
as needed.

\begin{figure}
\begin{pspicture}(-1.1,-3.1)(10,1.4)
\psset{unit=.4cm}
\psellipse(5,-1.5)(1.5,2.5)
\psarc[linewidth=.05](6.8,-1.5){2}{150}{210}\psarc[linewidth=.05](3.2,-1.5){2}{330}{30}
\pscircle(2.5,-1.5){1}\pscircle*(3.5,-1.5){.2}
\pscircle*(1.79,-.79){.25}\rput(1.25,-.8){\smsize{$1$}}
\pscircle*(1.79,-2.21){.25}\rput(1.23,-2.3){\smsize{$2$}}
\pscircle*(6.5,-1.5){.25}\rput(7,-1.5){\smsize{$3$}}
\psline(10,2)(10,-5)
\pscircle[fillstyle=solid,fillcolor=gray](11,2){1}\pscircle*(10,2){.25}
\pscircle[fillstyle=solid,fillcolor=gray](11,-1.5){1}\pscircle*(10,-1.5){.25}
\pscircle(11,-5){1}\pscircle*(10,-5){.25}
\pscircle[fillstyle=solid,fillcolor=gray](12.42,-3.58){1}\pscircle*(11.71,-4.29){.15}
\pscircle[fillstyle=solid,fillcolor=gray](12.42,-6.42){1}\pscircle*(11.71,-5.71){.15}
\rput(13,2){\smsize{$1$}}\rput(13,-1.5){\smsize{$2$}}\rput(14,-5){\smsize{$3$}}
\rput(8.4,-1.5){\Lgsize{$\times$}}
\psline{->}(15.5,-1.5)(20.5,-1.5)\rput(18,-1){\smsize{$\io_{0,\si^*}$}}
\psellipse(27,-1.5)(1.5,2.5)
\psarc[linewidth=.05](28.8,-1.5){2}{150}{210}\psarc[linewidth=.05](25.2,-1.5){2}{330}{30}
\pscircle(24.5,-1.5){1}\pscircle*(25.5,-1.5){.15}
\pscircle(29.5,-1.5){1}\pscircle*(28.5,-1.5){.25}\pnode(28.3,-1.7){A3}
\pscircle[fillstyle=solid,fillcolor=gray](23.09,-.09){1}
\pscircle[fillstyle=solid,fillcolor=gray](23.09,-2.91){1}
\pscircle[fillstyle=solid,fillcolor=gray](30.91,-.09){1}\pscircle*(30.21,-.79){.15}
\pscircle[fillstyle=solid,fillcolor=gray](30.91,-2.91){1}\pscircle*(30.21,-2.21){.15}
\pscircle*(23.79,-.79){.25}\pnode(23.95,-.95){A1}
\pscircle*(23.79,-2.21){.25}\pnode(24,-2.4){A2}
\rput(27,-7){\rnode{B}{\smsize{\begin{tabular}{c}distinguished\\nodes\end{tabular}}}}
\nccurve[nodesep=.1,angleA=92,angleB=-45,ncurv=.4]{->}{B}{A1}
\nccurve[nodesep=.1,angleA=94,angleB=-45,ncurv=.8]{->}{B}{A2}
\nccurve[nodesep=.1,angleA=88,angleB=-135,ncurv=.4]{->}{B}{A3}
\end{pspicture}
\caption{A Point in the Domain of $\io_{0,\si^*}$ and Its Image in $\ov\M_{1,k}(\Pn,d)$}
\label{idebund_fig}
\end{figure}

\begin{lmm}
\label{map1bl_lmm1b2}
If $d$, $n$, $k$, and $\si^*$ are as in Lemma~\ref{map1bl_lmm1b1},
$\si\!\in\!\A_1(d,k)$ is as above, and $\si\!\prec\!\si^*$, then
\begin{gather*}
\io_{0,\si^*}^{~-1} \big( \ov\M_{1,\si}^0\big)
=\Big(\bigcup_{\rho\in\A_P(\si^*;\si)} \!\!\!\!\!\!\!\!  \ov\cM_{1,\rho}^0 \Big)
\times \ov\M_{0,(\ale_B(\si^*),J_B(\si^*))},
\qquad\hbox{where}\\
\A_P(\si^*;\si)=\Big\{ 
\rho\!=\!\big(I_P\!\sqcup\!J_P,\{I_k\!\sqcup\!J_k\!: k\!\in\!K\}\big)
\!\in\!\A_1\big(I_P(\si^*),J_P(\si^*)\big)\!: |K|\!+\!|I_P|\!=\!m\Big\}
\end{gather*}
and $\A_1(I_P(\si^*),J_P(\si^*))$ and $\ov\cM_{1,\rho}^0\!\equiv\!\ov\cM_{1,\rho}$
are as in Subsection~\ref{curvebldata_subs}.
Furthermore, if $\rho\!\in\!\A_P(\si^*;\si)$ is as above,
$$\io_{0,\si^*}  \big|_{\ov\cM_{1,\rho}^0\times
\ov\M_{0,(\ale_B(\si^*),J_B(\si^*))}}^* T\ov\M_{1,\si}^0
\big/T\big(\ov\cM_{1,\rho}^0\!\times\!\ov\M_{0,(\ale_B(\si^*),J_B(\si^*))}\big)
 =\bigoplus_{i\in[m]-I_P}\!\!\!\!\!\!\pi_P^*L_i\!\otimes\!\pi_B^*\pi_i^*L_0.$$
\end{lmm}

\begin{lmm}
\label{map1bl_lmm1b3}
If $d$, $n$, $k$, $\si$, and $\si^*$ are as above, then
\begin{gather*}
\io_{0,\si}^{~-1} \big( \ov\M_{1,\si^*}^0\big)
=\ov\cM_{1,(I_P(\si),J_P(\si))}^0
\times \Big(\bigcup_{\vr\in\A_B(\si;\si^*)} \!\!\!\!\!\!\!\!  \ov\M_{0,\vr} \Big),
\qquad\hbox{where}\\
\A_B(\si;\si^*)=\big\{ \vr\!=\!\big((\si_l)_{l\in\ale_B(\si)},J_B\big)
\!\in\!\A_0\big(\ale_B(\si);d,J_B(\si)\big)\!:
\big|\ale_B(\vr)\big|\!=\!m^*\big\},
\end{gather*}
and $\A_0(\ale_B(\si);d,J_B(\si))$, $\ale_B(\vr)$, and
$\ov\M_{0,\vr}\!\equiv\!\ov\M_{0,\vr}(\Pn,d)$ are as in Subsection~\ref{map0prelim_subs}.
Furthermore, if $\vr\!\in\!\A_B(\si;\si^*)$ is as above,
$$\io_{0,\si}  \big|_{\ov\cM_{1,(I_P(\si),J_P(\si))}^0\times
\ov\M_{0,\vr}}^* T\ov\M_{1,\si^*}^0
\big/T\big(\ov\cM_{1,(I_P(\si),J_P(\si))}^0\!\times\!\ov\M_{0,\vr}\big)
 =\bigoplus_{i\in\aleph_P(\vr)}\!\!\!\!\pi_P^*L_i\!\otimes\!\pi_B^*\pi_i^*L_0,$$
where $\aleph_P(\vr)\!\subset\!\aleph_B(\si)$ is as in Subsection~\ref{map0prelim_subs}.
\end{lmm}

\noindent
We note that for every $\si^*\!\in\!\A_1(d,k)$,
$$\A_1\big(I_P(\si^*),J_P(\si^*)\big)= \bigsqcup_{\si\prec\si^*}\!\A_P(\si^*;\si).$$
Furthermore, if $\si_1,\si_2\!\in\!\A_1(d,k)$ are such that
$\si_1,\si_2\!\prec\!\si^*$, then
$$\rho_1\!\in\!\A_P(\si^*;\si_1), \quad \rho_2\!\in\!\A_P(\si^*;\si_2), \quad
\rho_1\!\prec\!\rho_2  \qquad\Lra\qquad \si_1\!\prec\!\si_2.$$
Thus, we can choose an ordering $<$ on $\A_1(I_P(\si^*),J_P(\si^*))$ 
extending the partial ordering~$\prec$ of Subsection~\ref{curve1bl_subs} such~that
$$\si_1\!<\!\si_2, \quad \rho_1\!\in\!\A_P(\si^*;\si_1), 
\quad \rho_2\!\in\!\A_P(\si^*;\si_2) \qquad\Lra\qquad \rho_1\!<\!\rho_2,$$
whenever $\si_1,\si_2\!\in\!\A_1(d,k)$ are such that $\si_1,\si_2\!\prec\!\si^*$.
In the next subsection, we will refer to the blowup construction of
Subsection~\ref{curve1bl_subs} corresponding to such an ordering.\\

\noindent
Similarly, if $\si'\!\in\!\A_1(d,k)$,
$$\A_0\big(\ale_B(\si');d,J_B(\si')\big)= \bigsqcup_{\si'\prec\si}\!\A_B(\si';\si).$$
Furthermore, if $\si_1,\si_2\!\in\!\A_1(d,k)$ are such that
$\si'\!\prec\!\si_1,\si_2$, then
$$\vr_1\!\in\!\A_B(\si';\si_1), \quad \vr_2\!\in\!\A_B(\si';\si_2), \quad
\vr_1\!\prec\!\vr_2  \qquad\Lra\qquad \si_1\!\prec\!\si_2.$$
Thus, we can choose an ordering $<$ on $\A_0(\ale_B(\si');d,J_B(\si'))$ 
extending the partial ordering~$\prec$ of Subsection~\ref{map0prelim_subs} such~that
$$\si_1\!<\!\si_2, \quad \vr_1\!\in\!\A_B(\si';\si_1), 
\quad \vr_2\!\in\!\A_B(\si';\si_2) \qquad\Lra\qquad \vr_1\!<\!\vr_2,$$
whenever $\si_1,\si_2\!\in\!\A_1(d,k)$ are such that $\si'\!\prec\!\si_1,\si_2$.
In the next subsection, we will refer to the blowup construction of
Subsection~\ref{map0blconstr_subs} corresponding to such an ordering.\\

\noindent
We denote by $\ov\M_{1,(0)}^0$ the main component $\ov\M_{1,k}^0(\Pn,d)$
of the moduli space $\ov\M_{1,k}(\Pn,d)$.
If $\si\!\in\!\A_1(d,k)$, we~put
\begin{alignat*}{1}
&\ov\cZ_{\si}^0= \io_{0,\si}^{~-1} \big( \ov\M_{1,(0)}^0\big)
\equiv \io_{0,\si}^{~-1} \big( \ov\M_{1,(0)}^0\!\cap\ov\M_{1,\si}^0\big);\\
&\cZ_{\si}^0= \io_{0,\si}^{~-1} \big( \ov\M_{1,(0)}^0\!\cap\M_{1,\si}^0\big)
\subset\ov\cZ_{\si}^0, \quad\hbox{where}\quad 
\M_{1,\si}^0=\M_{1,\si}^0(\Pn,d).
\end{alignat*}
We denote by  $\N\ov\cZ_{\si}^0\!\subset\!\N_{\io_{0,\si}}^{\ide}$ 
the normal cone $\N_{\io_{0,\si}|\ov\M_{1,(0)}^0}$ for 
$\io_{0,\si}|_{\ov\cZ_{\si}^0}$ in $\ov\M_{1,(0)}^0$.
Its structure is described in Lemma~\ref{map1bl_lmm2} below.
Let
$$ \cD_{0,\si} \in \Ga\big(\ov\cM_{1,(I_P(\si),J_P(\si))}^0 \!\times\!  
\ov\M_{0,(\ale_B(\si),J_B(\si))};
\Hom(\N_{\io_{0,\si}}^{\ide},\pi_P^*\E_0^*\!\otimes\!\pi_B^*\ev_0^*T\Pn)\big)$$
be the section defined by
$$ \cD_{0,\si}\big|_{\pi_P^*L_i\otimes\pi_B^*\pi_i^*L_0}
=\pi_P^*s_{0,i}\!\otimes\!\pi_B^*\pi_i^*\cD_0, 
\qquad\forall\,i\!\in\![m],$$
where $s_{0,i}$ and $\cD_0$ are as in 
Subsections~\ref{curve1bl_subs} and~\ref{map0str_subs}, respectively.

\begin{lmm}
\label{map1bl_lmm2}
For all $\si\!\in\!\A_1(d,k)$, $\ov\M_{1,(0)}^0$ is a proper subvariety
of $\ov\M_{1,k}^0$ relative to the immersions $\io_{0,\si}$ and $\bar\io_{0,\si}$.
Furthermore, 
\begin{gather*}
\cZ_{\si}^0=\big\{b\!\in\!
\cM_{1,(I_P(\si),J_P(\si))}^0 \!\times\!\M_{0,(\ale_B(\si),J_B(\si))}\!:
\ker\cD_{0,\si}|_b\!\neq\!\{0\}\big\}\\
\hbox{and}\qquad
\N\ov\cZ_{\si}^0\big|_{\cZ_{\si}^0}=\ker\cD_{0,\si}\big|_{\cZ_{\si}^0}.
\end{gather*}
Finally, $\ov\cZ_{\si}^0$ is the closure of $\cZ_{\si}^0$ in 
$\ov\cM_{1,(I_P(\si),J_P(\si))}^0 \!\times\!\ov\M_{0,(\ale_B(\si),J_B(\si))}$
and  $\N\ov\cZ_{\si}^0$ is the closure of
$\N\ov\cZ_{\si}^0\big|_{\cZ_{\si}^0}$ in~$\N_{\io_{0,\si}}^{\ide}$.
\end{lmm}

\noindent
This lemma is a consequence of~\cite[Theorem~\ref{g1comp-str_thm}]{g1comp} 
and related results.
In particular,  the first claim in the second sentence of Lemma~\ref{map1bl_lmm2} 
is a special case of the first statement of~\cite[Theorem~\ref{g1comp-str_thm}]{g1comp}.  
The second claim is nearly a special case of the last statement
of~\cite[Theorem~\ref{g1comp-str_thm}]{g1comp}, but some additional
comments are required.  \cite[Theorem~\ref{g1comp-str_thm}]{g1comp} by
itself is a purely topological statement, as it describes the
topological structure of a neighborhood of each stratum of
$\io_{0,\si}(\ov\cZ_{\si}^0)$ in~$\ov\M_{1,(0)}^0$.  On the other
hand, by Subsection~4.1 in~\cite{g1},
$\N\ov\cZ_{\si}^0\big|_{\cZ_{\si}^0}$ is contained in~$\ker\cD_{0,\si}$.  
The second claim in the second sentence of Lemma~\ref{map1bl_lmm2} can
then be obtained from a dimension count and a comparison of the gluing
construction used in the proof of \cite[Theorem~\ref{g1comp-str_thm}]{g1comp} 
with the analysis of limiting behavior in \cite[Subsect.~4.1]{g1}.  
This comparison implies that the gluing parameter in
the analytic construction of~\cite{g1comp} agrees to the first two orders
in the zero limit with the smoothing parameter in algebraic geometry.
Thus, $\N\ov\cZ_{\si}^0\big|_{\cZ_{\si}^0}$ must be equal
to~$\ker\cD_{\si,0}$.  Alternatively, suppose that $d\!\le\!n$.  If
the moduli space $\ov\M_{1,\si}^0$ is nonempty, then $m\!\le\!n$ and
thus for a Zariski open subset~$\cZ_{\si;1}$ of~$\cZ_{\si}^0$
\begin{equation}\label{map1bl_lmm2e3}
1\le\dim\N\ov\cZ_{\si}^0\big|_{\cZ_{\si;1}}=1
=\dim\ker\cD_{0,\si}\big|_{\cZ_{\si;1}}
\qquad\Lra\qquad
\N\ov\cZ_{\si}^0\big|_{\cZ_{\si;1}}=\ker\cD_{0,\si}\big|_{\cZ_{\si;1}}.
\end{equation}
Since $\cD_{0,\si}$ is transverse to the zero set over $\cZ_{\si}^0$,
the second claim in the second sentence of the lemma follows from~\e_ref{map1bl_lmm2e3},
if $d\!\le\!n$.
The general case follows from the observation that
$$\ov\M_{1,\si}(\Pn,d)=\big\{[\Si,u]\!\in\!\ov\M_{1,\si}(\P^{n+d},d)\!: 
u(\Si)\!\subset\!\Pn\big\}$$
and the $d\!\le\!n$ case.\\

\noindent
The first claim in the last sentence of Lemma~\ref{map1bl_lmm2} can be obtained by combining
the first statement of \cite[Theorem~\ref{g1comp-str_thm}]{g1comp},
the $m\!=\!1$ case of \cite[Theorem~2.8]{g2n2and3}, and the Implicit Function Theorem.
It also follows immediately from the last claim of Lemma~\ref{map1bl_lmm2}.
The latter can be deduced from \cite[Theorem~\ref{g1comp-str_thm}]{g1comp} as follows.
Suppose first that $m\!\le\!n$.
In this case, \cite[Theorem~2.8]{g2n2and3} implies that $\ov\cZ_{\si}^0$ admits 
a stratification 
$$\ov\cZ_{\si}^0=\cZ_{\si;1}\sqcup \bigsqcup_{\al\in\A}\cZ_{\si;\al}$$
such that $\cZ_{\si;1}$ is a Zariski open subset of $\ov\cZ_{\si}^0$,
\begin{gather}
\cZ_{\si;1}\subset\cZ_{\si}^0,  \qquad
\dim\,\N\ov\cZ_{\si}^0|_b\!=\!1 ~~~\forall\, b\!\in\!\cZ_{\si;1},\notag\\
\label{nofuzz_e1}
\hbox{and}\qquad
\max\big\{\!\dim\,\N\ov\cZ_{\si}^0|_b\!: b\!\in\!\cZ_{\si;\al}\big\}
\le \codim_{\ov\cZ_{\si}^0}\cZ_{\si;\al} ~~~\forall\,\al\!\in\!\A;
\end{gather}
see the next paragraph. Let 
$$\wt\cZ_{\si}^0= \P \N\ov\cZ_{\si}^0 \subset 
\P\N_{\io_{0,\si}}^{\ide}\big|_{\ov\cZ_{\si}^0}$$ 
be the exceptional divisor for the blowup of $\ov\M_{1,(0)}^0$ along $\ov\M_{1,\si}^0$.
Since all irreducible components of $\wt\cZ_{\si}^0$ must be of 
the same dimension, $\wt\cZ_{\si}^0$ must be the closure of $\wt\cZ_{\si}^0|_{\cZ_{\si}^0}$
by~\e_ref{nofuzz_e1}.
This closure property remains valid even if we do not assume that $m\!\le\!n$
for the following reason.
Let $pt\!\in\!\P^{n+d}$ be any point not contained in~$\Pn$.
Let
$$\pi\!: \P^{n+d}-\{pt\}\lra\Pn$$
be the corresponding linear projection.
It induces projection maps
\begin{alignat*}{1}
&\vph\!: \big\{[\Si,u]\!\in\!\ov\M_{1,k}^0(\P^{n+d},d)\!: pt\!\not\in\!u(\Si)\big\}
\lra \ov\M_{1,k}^0(\Pn,d) \qquad\hbox{and}\\
&\ti\vph\!: \big\{[\Si,u;v]\!\in\!\wt\cZ_{\si}^0(\P^{n+d},d)\!: pt\!\not\in\!u(\Si)\big\}
\lra \wt\cZ_{\si}^0(\Pn,d).
\end{alignat*}
The latter map takes $\wt\cZ_{\si}^0(\P^{n+d},d)|_{\cZ_{\si}^0(\P^{n+d},d)}$ to
$\wt\cZ_{\si}^0(\Pn,d)|_{\cZ_{\si}^0(\Pn,d)}$.
Since the closure of 
$$\wt\cZ_{\si}^0(\P^{n+d},d)|_{\cZ_{\si}^0(\P^{n+d},d)}$$
contains $\wt\cZ_{\si}^0(\Pn,d)$, it follows that
so does the closure of $\wt\cZ_{\si}^0(\Pn,d)|_{\cZ_{\si}^0(\Pn,d)}$.
This observation implies the last claim of Lemma~\ref{map1bl_lmm2}.\\

\noindent 
We conclude this subsection by briefly describing the stratification mentioned above.
A stratum $\M_{\Ga_B}$ of $\ov\M_{0,(\ale_B(\si),J_B(\si))}$ corresponds 
to a tuple $\Ga_B\!\equiv\!(\Ga_{B;l})_{l\in\aleph_B(\si)}$ of dual graphs,
all of which are trees.
The vertices of $\Ga_{B;l}$ correspond to the irreducible components of the domain
of the stable map~$b_l$ in the definition of $\ov\M_{0,(\ale_B(\si),J_B(\si))}$
at the beginning of Subsection~\ref{map0prelim_subs}.
Each vertex~$v$ of $\Ga_{B;l}$ is labeled by a nonnegative integer, which specifies 
the degree of the stable map $b_l$ on the corresponding component~$\Si_v$.
There is an edge in~$\Ga_{B;l}$ between two vertices if and only if 
the two corresponding components of the domain share a node.
In addition, there are tails attached at some vertices of $\Ga_{B;l}$,
which are labeled by the indexing set for marked points of the map~$b_l$,
i.e.~$J_{l,P}$ in the notation of~Subsection~\ref{map0prelim_subs}. 
Let $v_l^*$ be the vertex of $\Ga_{B;l}$ to which the tail corresponding to 
the marked point~$0$ is attached. If the degree of $v_l^*$ is positive, let
$$\chi_l(\Ga_B)\equiv\chi_l(\Ga_{B;l})=\{v_l^*\}.$$
Otherwise, denote by $\chi_l(\Ga_B)$ the set of positive-degree vertices of $\Ga_{B;l}$ 
that are not separated from~$v_l^*$ by a positive-degree vertex.
Suppose
$$b\!\equiv\!(b_l)_{l\in\ale_B(\si)} \in  \M_{\Ga_B}\equiv 
 \ov\M_{0,(\ale_B(\si),J_B(\si))} \cap\prod_{l\in\ale_B(\si)}\!\!\!\!\!\M_{\Ga_{B;l}},
\qquad\hbox{with}\qquad b_l\!=\![\Si_l,u_l]$$
as in the paragraph preceding Lemma~\ref{deriv0str_lmm}.
If $l\!\in\!\ale_B(\si)$ and $v\!=\!v_l^*$, let 
$$\Im\,\cD_v|_b=\Im\,\cD_0|_{b_l}\equiv\Im\, du_l|_{x_0(b_l)}\subset T_{\ev_0(b)}\Pn.$$ 
If $v$ is a vertex of $\Ga_{B;l}$ different from~$v_l^*$, we denote by
$\Im\,\cD_v|_b$ the image of $d\{u_l|_{\Si_v}\}$ at the node of $\Si_v$ corresponding to 
the edge of $\Ga_{B;l}$ that leaves $v$ on the unique path from $v$ on $v_l^*$ in~$\Ga_{B;l}$.
Note that if $v\!\in\!\chi_l(\Ga_B)$, the image of this node under $u_l$ is~$\ev_0(b)$.
We~set
$$\chi(\Ga)=\bigsqcup_{l\in\ale_B(\si)}\!\!\!\!\chi_l(\Ga_B).$$
With $b$ as above, let
$$\codim\,\cD|_b=\big|\chi(\Ga_B)\big|-
\dim\hbox{Span}\big\{\Im\,\cD_v|_b\!:v\!\in\!\chi_l(\Ga_B),\,l\!\in\!\ale_B(\si)\big\}.$$
For each pair 
$\al\!=\!(\Ga_B,\mu)$, where $\mu\!\in\!\Z^+$ is such that
\begin{equation}\label{alcond_e}
\max\big(1,|\chi(\Ga_B)|\!-\!n\big)\le \mu\le|\chi(\Ga_B)|, 
\end{equation}
we put 
$$\cZ_{\Ga_B;\al}=\big\{b\!\in\!\M_{\Ga_B}\!: \codim\,\cD|_b\!=\!\mu\big\}.$$
By the first statement of \cite[Theorem~\ref{g1comp-str_thm}]{g1comp},
$$\ov\cZ_{\si}^0=\bigsqcup_{\al}\cZ_{\si;\al},
\qquad\hbox{where}\qquad
\cZ_{\si;\al}=\ov\cM_{1,(I_P(\si),J_P(\si))}^0\!\times\!\cZ_{\Ga_B;\al}.$$
The disjoint union is taken over all pairs $\al\!=\!(\Ga,\mu)$ as described above.
From transversality as in the first claim of Lemma~\ref{deriv0str_lmm}, 
it is easy to see that
\begin{equation}\label{codimest_e}\begin{split}
\codim_{\M_{\Ga_B}} \cZ_{\Ga_B;\al}
&=\big(n-(|\chi(\Ga_B)|\!-\!\mu)\big)\mu\\
&\ge n-(|\chi(\Ga_B)|\!-\!\mu);
\end{split}\end{equation}
see the end of \cite[Subsect.~\ref{g1cone-g1str_subs}]{g1cone}, for example.
The above inequality follows from the first inequality in~\e_ref{alcond_e}. 
By~\e_ref{codimest_e}, if $m\!=\!|\ale_B(\si)|\!\le\!n$, 
\begin{equation*}\begin{split}
\codim_{\ov\cZ_{\si}^0} \cZ_{\si;\al} 
&= \codim_{\M_{\Ga_B}} \cZ_{\Ga_B;\al} + 
\codim_{\ov\M_{0,(\ale_B(\si),J_B(\si))}}\M_{\Ga_B}\\
&\qquad\qquad\qquad\qquad\qquad
- \codim_{\ov\cM_{1,(I_P(\si),J_P(\si))}^0\times\ov\M_{0,(\ale_B(\si),J_B(\si))}} 
\ov\cZ_{\si}^0 \\
&\ge  \big(n\!-\!|\chi(\Ga_B)|\!+\!\mu\big)+\big(|\chi(\Ga_B)|\!-\!|\ale_B(\si)|\big)
-\big(n\!-\!|\ale_B(\si)|\!+\!1\big)
=\mu\!-\!1.
\end{split}\end{equation*}
On the other hand, by the last statement of \cite[Theorem~\ref{g1comp-str_thm}]{g1comp},
$$\max\big\{\!\dim\,\N\ov\cZ_{\si}^0|_b\!: b\!\in\!\cZ_{\si;\al}\big\}= \mu.$$
We conclude that 
$$\max\big\{\!\dim\,\N\ov\cZ_{\si}^0|_b\!: b\!\in\!\cZ_{\si;\al}\big\}
\le \codim_{\ov\cZ_{\si}^0} \cZ_{\si;\al} +1.$$
The equality holds if and only if $\mu\!=\!1$ and $\Ga_B$ is a tuple of one-vertex graphs, 
i.e.~the image of $\cM_{1,(I_P(\si),J_P(\si))}\!\times\!\M_{\Ga_B}$
under $\io_{0,\si}$ is contained in $\M_{1,\si}$.
This observation concludes the proof of the stratification claim made in the previous
paragraph.

\subsection{Inductive Construction}
\label{map1blconstr_subs}

\noindent
This subsection is the analogue of Subsection~\ref{map0blconstr_subs}
in the present situation.
Suppose $\si\!\in\!\A_1(d,k)$ and we have constructed\\
${}\quad$ ($I1$) an idealized blowup $\pi_{\si-1}\!:\ov\M_{1,k}^{\si-1}\!\lra\!\ov\M_{1,k}^0$
such that $\pi_{\si-1}$ is an isomorphism outside of the preimages of the subvarieties
$\ov\M_{1,\si'}^0$ with $\si'\!\le\!\si\!-\!1$;\\
${}\quad$ ($I2$) for each $\si'\!\in\!\{(0)\}\!\sqcup\!\A_1(d,k)$, 
a subvariety $\ov\M_{1,\si'}^{\si-1}$ of $\ov\M_{1,k}^{\si-1}$
such~that
$$\ov\M_{1,k}^{\si-1}=\ov\M_{1,(0)}^{\si-1}
\cup\bigcup_{\si'\in\A_1(d,k)} \!\!\!\!\!\!\! \ov\M_{1,\si'}^{\si-1}, \qquad
\pi_{\si-1}\big(\ov\M_{1,\si'}^{\si-1}\big)=\ov\M_{1,\si'}^0
~~~\forall\, \si'\!\in\!\{(0)\}\!\sqcup\!\A_1(d,k),$$
and $\ov\M_{1,\si^*}^{\si-1}$ is the proper transform of $\ov\M_{1,\si^*}^0$
for $\si^*\!=\!(0)$ and 
for all $\si^*\!\in\!\A_1(d,k)$ such that $\si^*\!>\!\si\!-\!1$.\\
We assume that\\ 
${}\quad$ ($I3$) for all $\si_1,\si_2\!\in\!\A_1(d,k)$ such that $\si_1\!\neq\!\si_2$,
$\si_1\!\not\prec\si_2$, $\si_2\!\not\prec\si_1$, and $\si\!-\!1\!<\!\si_1,\si_2$,
$$ \ov\M_{1,\si_1}^{\si-1}\cap \ov\M_{1,\si_2}^{\si-1} ~
\begin{cases}
\subset\ov\M_{1,\ti\si(\si_1,\si_2)}^{\si-1}, & 
\hbox{if}~ \ti\si(\si_1,\si_2)\!>\!\si\!-\!1;\\
=\eset,& \hbox{otherwise},\\
\end{cases}$$
where  $\ti\si(\si_1,\si_2)$ is as in Lemma~\ref{map1bl_lmm1a}.\\

\noindent
We also assume that for every $\si'\!\in\!A_1(d,k)$ such that $\si'\!\le\!\si\!-\!1$:\\
${}\quad$ ($I4$) $\ov\M_{1,\si'}^{\si-1}$ is the image of a $G_{\si'}$-invariant immersion
\begin{gather*}
\io_{\si-1,\si'}\!: \wt\cM_{1,(I_P(\si'),J_P(\si'))}
\times \wt\M_{0,(\ale_B(\si'),J_B(\si'))}^{\vr_{\si'}(\si-1)}
\lra \ov\M_{1,k}^{\si-1},
\qquad\hbox{where}\\
\vr_{\si'}(\si\!-\!1)=
\begin{cases}
\max\big\{\vr\!\in\!\A_B(\si';\si^*)\!: \si'\!\prec\!\si^*\!\le\!\si\!-\!1\big\},
&\hbox{if}~\exists\si^*\!\in\!\A_1(d,k)
~\hbox{s.t.}~\si'\!\prec\!\si^*\!\le\!\si\!-\!1;\\
0,& \hbox{otherwise},
\end{cases}
\end{gather*}
and $\wt\M_{0,(\ale_B(\si'),J_B(\si'))}^{\vr_{\si'}(\si-1)}
\!\equiv\!\wt\M_{0,(\ale_B(\si'),J_B(\si'))}^{\vr_{\si'}(\si-1)}(\Pn,d)$
is the blowup $\wt\M_{0,(\ale_B(\si'),J_B(\si'))}^0(\Pn,d)$
 constructed in Subsection~\ref{map0blconstr_subs};\\
${}\quad$ ($I5$) if $\si^*\!\in\!\A_1(d,k)$ is such that $\si\!-\!1\!<\!\si^*$ and 
$\si'\!\prec\!\si^*$, then
$$\io_{\si-1,\si'}^{\,-1}\big(\ov\M_{1,\si^*}^{\si-1}\big)
= \wt\cM_{1,(I_P(\si'),J_P(\si'))} \times \Big(\bigcup_{\vr\in\A_B(\si';\si^*)}
\!\!\!\!\!\!\!\!\!\wt\M_{0,\vr}^{\vr_{\si'}(\si-1)}\Big),$$
where $\wt\M_{0,\vr}^{\vr_{\si'}(\si-1)}\!\equiv\!\wt\M_{0,\vr}^{\vr_{\si'}(\si-1)}(\Pn,d)$
is the subvariety of $\wt\M_{0,(\ale_B(\si'),J_B(\si'))}^{\vr_{\si'}(\si-1)}$ 
described in Subsection~\ref{map0blconstr_subs};\\
${}\quad$ ($I6$) an idealized normal bundle for the immersion
$\io_{\si-1,\si'}$ is given by
$$\N_{\io_{\si-1,\si'}}^{\ide} = \pi_P^*\L \otimes
\pi_B^*\pi_{\vr_{\si'}(\si-1)}^{~*}\ga_{(\ale_B(\si'),J_B(\si'))},$$
where
$$\pi_B,\, \pi_P\!:
\wt\cM_{1,(I_P(\si'),J_P(\si'))}
\times \wt\M_{0,(\ale_B(\si'),J_B(\si'))}^{\vr_{\si'}(\si-1)}
\lra \wt\cM_{1,(I_P(\si'),J_P(\si'))}, \,
\wt\M_{0,(\ale_B(\si'),J_B(\si'))}^{\vr_{\si'}(\si-1)}$$
are the two projection maps and $\L\!\lra\!\wt\cM_{1,(I_P(\si'),J_P(\si'))}$
is the universal tangent line bundle of Subsection~\ref{curve1bl_subs};\\
${}\quad$ ($I7$) $\ov\cZ_{\si'}^{\si-1}\!\equiv\!
\io_{\si-1,\si'}^{~-1}(\ov\M_{1,(0)}^{\si-1})$ is the closure of
$$\cZ_{\si'}^{\si-1}\equiv \wt\cM_{1,(I_P(\si'),J_P(\si'))} \times 
\Big(\wt\cD_{\vr_{\si'}(\si-1)}^{\,-1}(0)
-\bigcup_{\stackrel{\vr\in\A_B(\ale_B(\si');d,J_B(\si'))}{\vr_{\si'}(\si-1)<\vr}}
\!\!\!\!\!\! \wt\M_{0,\vr}^{\vr_{\si'}(\si-1)}\Big)$$
in $\wt\cM_{1,(I_P(\si'),J_P(\si'))}
\!\times\! \wt\M_{0,(\ale_B(\si'),J_B(\si'))}^{\vr_{\si'}(\si-1)}$ and
$$\N\ov\cZ_{\si'}^{\si-1}\equiv\N_{\io_{\si-1,\si'}|\ov\M_{1,(0)}^{\si-1}}
=\N_{\io_{\si-1,\si'}}^{\ide}\big|_{\ov\cZ_{\si'}^{\si-1}}$$
is the normal cone for $\io_{\si-1,\si'}|\ov\M_{1,(0)}^{\si-1}$ in $\ov\M_{1,k}^{\si-1}$;\\
${}\quad$ ($I8$) the immersion map
$$\bar\io_{\si-1,\si'}\!:\big( \wt\cM_{1,(I_P(\si');J_P(\si'))}
\!\times\! \wt\M_{0,(\ale_B(\si'),J_B(\si'))}^{\vr_{\si'}(\si-1)}\big)\big/G_{\si'}
\lra \ov\M_{1,k}^{\si-1}$$
induced by $\io_{\si-1,\si'}$ is an embedding.\\

\noindent
Furthermore, 
we assume that for every $\si^*\!\in\!A_1(d,k)$ such that $\si^*\!>\!\si\!-\!1$:\\
${}\quad$ ($I9$)  the domain of the $G_{\si^*}$-invariant immersion $\io_{\si-1,\si^*}$
induced by $\io_{0,\si^*}$~is
\begin{gather*}
\ov\cM_{1,(I_P(\si^*),J_P(\si^*))}^{\rho_{\si^*}(\si-1)}
\times \ov\M_{0,(\ale_B(\si^*),J_B(\si^*))},
\qquad\hbox{where}\\
\rho_{\si^*}(\si\!-\!1)=
\begin{cases}
\max\big\{\rho\!\in\!\A_P(\si^*;\si')\!: 
\si'\!\le\!\si\!-\!1,\si'\!\prec\!\si^*\big\},
& \begin{split}
  &\hbox{if}~\exists\si'\!\in\!\A_1(d,k)\\
  &~\hbox{s.t.}~\si'\!\le\!\si\!-\!1,\si'\!\prec\!\si^*;
  \end{split}\\
0,& \hbox{otherwise},
\end{cases}
\end{gather*}
and $\ov\cM_{1,(I_P(\si^*),J_P(\si^*))}^{\rho_{\si^*}(\si-1)}
\!\lra\!\ov\cM_{1,(I_P(\si^*),J_P(\si^*))}$
is the blowup constructed in Subsection~\ref{curve1bl_subs};\\
${}\quad$ ($I10$)  if $\si'\!\in\!\A_1(d,k)$ is such that
$\si\!-\!1\!<\!\si'\!\prec\!\si^*$, then
$$\io_{\si-1,\si^*}^{~-1} \big( \wt\M_{1,\si'}^{\si-1} \big)
=\Big(\bigcup_{\rho\in\A_P(\si^*;\si')} \!\!\!\!\!\!\!  
\wt\cM_{1,\rho}^{\rho_{\si^*}(\si-1)} \Big)
\times \ov\M_{0,(\ale_B(\si^*),J_B(\si^*))};$$
${}\quad$ ($I11$) if $\si^*$ is as in Lemma~\ref{map1bl_lmm1b1},
an idealized normal bundle for the immersion $\io_{\si-1,\si^*}$ is given by
$$\N_{\io_{\si-1,\si^*}}^{\ide}= 
\bigoplus_{i\in[m^*]} \pi_P^*L_{\rho_{\si^*}(\si-1),i} \!\otimes\!
\pi_B^*\pi_i^*L_0,$$
where
$$\pi_P,\, \pi_B\!:
\ov\cM_{1,(I_P(\si^*),J_P(\si^*))}^{\rho_{\si^*}(\si-1)}
\!\times\! \ov\M_{0,(\ale_B(\si^*),J_B(\si^*))} \lra
\ov\cM_{1,(I_P(\si^*),J_P(\si^*))}^{\rho_{\si^*}(\si-1)}, \,
\ov\M_{0,(\ale_B(\si^*),J_B(\si^*))}$$
are the two projection maps and 
$L_{\rho_{\si^*}(\si-1),i}\!\lra\!\ov\cM_{1,(I_P(\si^*),J_P(\si^*))}^{\rho_{\si^*}(\si-1)}$
is the line bundle constructed in Subsection~\ref{curve1bl_subs};\\
${}\quad$ ($I12$) $\ov\cZ_{\si^*}^{\si-1}\!\equiv\!
\io_{\si-1,\si^*}^{~-1}(\ov\M_{1,(0)}^{\si-1})$ is the closure of $\cZ_{\si^*}^0$ in
$\ov\cM_{1,(I_P(\si^*),J_P(\si^*))}^{\rho_{\si^*}(\si-1)}
\!\times\!\ov\M_{0,(\ale_B(\si^*),J_B(\si^*))}$ and the normal cone
$$\N\ov\cZ_{\si^*}^{\si-1}\equiv\N_{\io_{\si-1,\si^*}|\ov\M_{1,(0)}^{\si-1}}$$
for  $\io_{\si-1,\si^*}|\ov\cZ_{\si^*}^{\si-1}$ is the closure of
$\N\ov\cZ_{\si^*}^0\big|_{\cZ_{\si^*}^0}$ in~$\N_{\io_{\si-1,\si^*}}^{\ide}$;\\
${}\quad$ ($I13$) $\Im^s\,\bar\io_{\si-1,\si^*}\subset\bigcup_{\si-1<\si'\prec\si^*}
\wt\M_{0,\si'}^{\si-1}$, where
$$\bar\io_{\si-1,\si^*}\!:\big(  
\ov\cM_{1,(I_P(\si^*),J_P(\si^*))}^{\rho_{\si^*}(\si-1)}
\!\times\! \ov\M_{0,(\ale_B(\si^*),J_B(\si^*))}\big)\big/G_{\si^*}
\lra \ov\M_{1,k}^{\si-1},$$
is the immersion map induced by $\io_{\si-1,\si^*}$.\\

\noindent
Finally, we assume that\\
${}\quad$ ($I14$) the collections
$\{\io_{\si-1,\si'}\}_{\si'\in\A_1(d,k)}$ and $\{\bar\io_{\si-1,\si'}\}_{\si'\in\A_1(d,k)}$ 
of immersions are properly self-intersecting;\\
${}\quad$ ($I15$) for all $\si'\!\in\!\A_1(d,k)$,
the subvariety $\ov\M_{1,(0)}^{\si-1}$ of $\ov\M_{1,k}^{\si-1}$ is proper
relative to the immersions~$\io_{\si-1,\si'}$ and~$\bar\io_{\si-1,\si'}$.\\

\noindent
By the inductive assumption~($I3$), if $\si_1$ and $\si_2$ are non-comparable elements
of $(\A_1(d,k),\prec)$, the proper transforms of $\ov\M_{1,\si_1}^0$ and 
$\ov\M_{1,\si_2}^0$ become disjoint by the time either is ready to be blown up
for any ordering~$<$ extending the partial ordering~$\prec$.
Similarly to the three blowup constructions encountered previously,
($I3$) will imply that the end result of the present blowup construction is 
independent of the choice of an extension~$<$.
By~($I9$), our blowup construction modifies each immersion~$\io_{0,\si^*}$
by changing the first factor of the domain according to the blowup
construction of Subsection~\ref{curve1bl_subs}, until a proper transform
of the image of~$\io_{0,\si^*}$ is to be blown~up; see below.
By~($I11$) and~($I13$), in the process, the singular locus of $\io_{0,\si^*}$
disappears and the first component in every summand of 
$\N_{\io_{0,\si^*}}^{\ide}$ gets twisted to~$\L$.
In particular, all blowup loci are smooth.
On the other hand, by the inductive assumptions~($I7$) and~($I8$),
for $\si'\!\le\!\si\!-\!1$ the intersection of the proper transform of $\ov\M_{1,(0)}^0$
with the proper transform of the exceptional divisor $\ov\M_{1,\si'}^{\si'}$
is an embedding of a subvariety of a smooth variety.
The singular locus of this subvariety is annihilated by the time the entire blowup
construction is complete, according to the inductive assumptions~($I7$) above
and the inductive assumption~($I4$) in Subsection~\ref{map0blconstr_subs}.
These assumptions imply that the proper transform of $\ov\M_{1,(0)}^0$ after
the final blowup step is smooth.\\

\noindent
We note that all of the assumptions ($I1$)-($I15$) are satisfied if $\si\!-\!1$ 
is replaced by~$0$.
In particular, ($I4$) and~($I12$) are restatements of Lemmas~\ref{map1bl_lmm1a}
and~\ref{map1bl_lmm2}, respectively.
The statements ($I10$), ($I11$), and ($I13$)-($I15$), with $\si\!-\!1$ replaced by~$0$,
are contained in Lemmas~\ref{map1bl_lmm1b1} and~\ref{map1bl_lmm1b2}.\\

\noindent
If $\si\!\in\!\A_1(d,k)$ is as above, let
$$\ti\pi_{\si}\!: \ov\M_{1,k}^{\si}\lra\ov\M_{1,k}^{\si-1}$$
be the idealized blowup of $\ov\M_{1,k}^{\si-1}$ along $\ov\M_{1,\si}^{\si-1}$,
which is a smooth subvariety by the inductive assumption~($I13$).
We denote the idealized exceptional divisor,
$$\cE_{\ov\M_{1,\si}^{\si-1}}^{\ide}=\P\N_{\bar\io_{\si-1,\si}}^{\ide},$$
by $\ov\M_{1,\si}^{\si}$.
For each $\si'\!\in\!\{(0)\}\!\sqcup\!(\A_1(d,k)\!-\!\{\si\})$, we denote~by
$$\ov\M_{1,\si'}^{\si} \subset \Bl_{\ov\M_{1,\si}^{\si-1}}\ov\M_{1,k}^{\si-1}
\subset \ov\M_{1,k}^{\si}$$
the proper transform of $\ov\M_{1,\si'}^{\si-1}$.
Let $\pi_{\si}\!=\!\pi_{\si-1}\!\circ\!\ti\pi_{\si}$.\\

\noindent
The inductive assumptions ($I1$) and ($I2$),
with $\si\!-\!1$ replaced by $\si$, are clearly satisfied, 
while ($I3$), ($I8$) for $\si'\!\neq\!\si$,
and ($I13$)-($I15$) follow from~(2) of Lemma~\ref{ag_lmm2a}, 
Corollary~\ref{immercoll_crl}, and Lemma~\ref{virimmer_lmm3}.
On the other hand, by~($I9$), the domain of the immersion $\io_{\si-1,\si}$~is
$$\ov\cM_{1,(I_P(\si),J_P(\si))}^{\rho_{\si}(\si-1)}
\!\times\! \ov\M_{0,(\ale_B(\si),J_B(\si))}
=\wt\cM_{1,(I_P(\si),J_P(\si))}
\!\times\! \ov\M_{0,(\ale_B(\si),J_B(\si))}.$$
By~($I11$), the chosen idealized normal bundle for the immersion~$\io_{\si-1,\si}$ is given~by
\begin{equation}\label{sinormal_e}
\N_{\io_{\si-1,\si}}^{\ide}= 
\bigoplus_{i\in[m]} \pi_P^*L_{\rho_{\si}(\si-1),i} \!\otimes\!
\pi_B^*\pi_i^*L_0
=\pi_P^*\L\otimes \pi_B^*F_{(\ale_B(\si),J_B(\si))}.
\end{equation}
Thus, the domain of the immersion~$\io_{\si,\si}$ induced by $\io_{\si-1,\si}$~is
$$\P\N_{\io_{\si-1,\si}}^{\ide} = \wt\cM_{1,(I_P(\si),J_P(\si))}
\!\times\! \wt\M_{0,(\ale_B(\si),J_B(\si))}^0
=\wt\cM_{1,(I_P(\si),J_P(\si))}
\!\times\! \wt\M_{0,(\ale_B(\si),J_B(\si))}^{\vr_{\si}(\si)}.$$
By the first statement of Lemma~\ref{virimmer_lmm2},
an idealized normal bundle for the embedding $\io_{\si,\si}$
is the tautological line bundle over $\P\N_{\io_{\si-1,\si}}^{\ide}$,~i.e.
$$\N_{\io_{\si,\si}}^{\ide}
= \pi_P^*\L \!\otimes\! \pi_B^*\ga_{(\ale_B(\si),J_B(\si))}
= \pi_P^*\L \!\otimes\! \pi_B^*\pi_{\vr_{\si}(\si)}^*\ga_{(\ale_B(\si),J_B(\si))}.$$
Thus, the inductive assumptions ($I4$) and~($I6$), with $\si'\!=\!\si$ and
$\si\!-\!1$ replaced by~$\si$, are satisfied.
The same is the case with~($I8$), since the map $\bar\io_{\si-1,\si}$
is an embedding by~($I13$).\\

\noindent
We also note that by the first statement of Lemma~\ref{map1bl_lmm1b3}, the inductive assumptions~($I1$) and~($I2$), and the last statement of Lemma~\ref{ag_lmm2a},
$$\io_{\si-1,\si}^{~-1} \big( \ov\M_{1,\si^*}^{\si-1}\big)
=\wt\cM_{1,(I_P(\si),J_P(\si))}
\times \Big(\bigcup_{\vr\in\A_B(\si;\si^*)} \!\!\!\!\!\!\!\!  \ov\M_{0,\vr} \Big)$$
for all $\si^*\!\in\!\A_1(d,k)$ such that $\si\!\prec\!\si^*$.
In addition, by the last statement of Lemma~\ref{map1bl_lmm1b3} 
$$\io_{\si-1,\si}  \big|_{\wt\cM_{1,(I_P(\si),J_P(\si))}\times
\ov\M_{0,\vr}}^* T\ov\M_{1,\si^*}^{\si-1}
\big/T\big(\wt\cM_{1,(I_P(\si),J_P(\si))}\!\times\!\ov\M_{0,\vr}\big)
 \subset\N_{\io_{\si-1,\si}}^{\ide}$$
is a vector bundle for all $\vr\!\in\!\A_B(\si;\si^*)$ and
$$\io_{\si-1,\si}  \big|_{\cM_{1,(I_P(\si),J_P(\si))}\times
\M_{0,\vr}}^* T\ov\M_{1,\si^*}^{\si-1}
\big/T\big(\cM_{1,(I_P(\si),J_P(\si))}\!\times\!\M_{0,\vr}\big)
=\bigoplus_{i\in\aleph_P(\vr)}\!\!\!\!\pi_P^*L_i\!\otimes\!\pi_B^*\pi_i^*L_0.$$
Thus, by the first equality in~\e_ref{sinormal_e},
\begin{equation*}\begin{split}
\io_{\si-1,\si}  \big|_{\wt\cM_{1,(I_P(\si),J_P(\si))}\times
\ov\M_{0,\vr}}^* T\ov\M_{1,\si^*}^{\si-1}
\big/T\big(\wt\cM_{1,(I_P(\si),J_P(\si))}\!\times\!\ov\M_{0,\vr}\big)
&=\pi_P^*\L\!\otimes\bigoplus_{i\in\aleph_P(\vr)}\!\!\!\!\pi_B^*\pi_i^*L_0\\
&=\pi_P^*\L \otimes\pi_B^*F_{\vr;P}.
\end{split}\end{equation*}
It follows that
\begin{equation*}\begin{split}
\io_{\si,\si}^{~-1} \big( \ov\M_{1,\si^*}^{\si}\big)
&=\bigcup_{\vr\in\A_B(\si;\si^*)} \!\!\!\!\!\!\!\! 
\P\big(\pi_P^*\L \otimes\pi_B^*F_{\vr;P}\big)\big|_{
\wt\cM_{1,(I_P(\si),J_P(\si))}\times\ov\M_{0,\vr}} \\
&=\wt\cM_{1,(I_P(\si),J_P(\si))}\times 
\Big(\bigcup_{\vr\in\A_B(\si;\si^*)} \!\!\!\!\!\!\!\!\P F_{\vr;P}\Big)\\
&\equiv\wt\cM_{1,(I_P(\si),J_P(\si))}\times 
\Big(\bigcup_{\vr\in\A_B(\si;\si^*)}\!\!\!\!\!\!\!\!\wt\M_{0,\vr}^0\Big)
=\wt\cM_{1,(I_P(\si),J_P(\si))}\times 
\Big(\bigcup_{\vr\in\A_B(\si;\si^*)}\!\!\!\!\!\!\!\!\wt\M_{0,\vr}^{\vr_{\si}(\si)}\Big),
\end{split}\end{equation*}
as needed for the inductive assumption ($I5$) with $\si\!-\!1$ replaced by $\si$
and $\si'\!=\!\si$.\\

\noindent
Furthermore, by~($I12$), $\N\ov\cZ_{\si}^{\si-1}$ is the closure~of
\begin{equation*}\begin{split}
\N\ov\cZ_{\si}^0\big|_{\cZ_{\si}^0}
&\equiv  \ker\cD_{0,\si}\big|_{
\cM_{1,(I_P(\si),J_P(\si))}^0 \times\M_{0,(\ale_B(\si),J_B(\si))}}\\
&= \pi_P^*\L \otimes
\pi_B^*\ker\cD_{(\ale_B(\si),J_B(\si))}\big|_{\M_{0,(\ale_B(\si),J_B(\si))}}
\end{split}\end{equation*}
in $\pi_P^*\L\!\otimes\!\pi_B^*F_{(\ale_B(\si),J_B(\si))}$, 
where $\cD_{(\ale_B(\si),J_B(\si))}$ 
is the bundle homomorphism described in Subsection~\ref{map0prelim_subs}.
Thus, by the first statement of Lemma~\ref{virimmer_lmm3},
$$\ov\cZ_{\si}^{\si}\equiv\io_{\si,\si}^{~-1}(\ov\M_{1,(0)}^{\si})$$
is the closure of 
$$\cM_{1,(I_P(\si),J_P(\si))}\times 
\big\{b\!\in\!\P F_{(\ale_B(\si),J_B(\si))}|_{\M_{0,(\ale_B(\si),J_B(\si))}}\!:
\wt\cD_0b\!=\!0\big\}$$
in $\wt\cM_{1,(I_P(\si),J_P(\si))}
    \!\times\! \wt\M_{0,(\ale_B(\si),J_B(\si))}^{\vr_{\si}(\si)}$.
The inductive assumption~($I7$), with $\si'\!=\!\si$ and
$\si\!-\!1$ replaced by~$\si$, now follows from the first statement of 
Lemma~\ref{map0bl_lmm2}.\\

\noindent
We next verify that the inductive assumptions ($I4$)-($I7$)
hold for $\si'\!<\!\si$, with $\si\!-\!1$ replaced by~$\si$.
If $\si'\!\not\prec\!\si$, then
$$\vr_{\si'}(\si)=\vr_{\si'}(\si\!-\!1)
\qquad\hbox{and}\qquad
\ov\M_{1,\si'}^{\si-1}\cap\ov\M_{1,\si}^{\si-1}=\eset,$$
by definition and ($I3$), respectively. It then follows that
\begin{gather*}
\io_{\si,\si'}=\io_{\si-1,\si'}, \qquad
\N_{\io_{\si,\si'}}^{\ide}=\N_{\io_{\si-1,\si'}}^{\ide},\\
\hbox{and}\qquad  \ov\M_{1,\si'}^{\si}\cap\ov\M_{1,\si^*}^{\si}
=\ov\M_{1,\si'}^{\si-1}\cap\ov\M_{1,\si^*}^{\si-1}
\quad\forall\, \si^*\!\in\!\{(0)\}\!\sqcup\!\A_1(d,k).
\end{gather*}
Thus, the inductive assumptions ($I4$)-($I7$),
as stated above, imply the corresponding statements with $\si\!-\!1$ replaced by~$\si$.\\

\noindent
Suppose that $\si'\!\prec\!\si$.
By~($I4$) and~(1) of Lemma~\ref{ag_lmm2a},
the domain of the immersion $\io_{\si,\si'}$ induced
by~$\io_{\si-1,\si'}$ is the blowup~of 
$$\wt\cM_{1,(I_P(\si'),J_P(\si'))}
\times \wt\M_{0,(\ale_B(\si'),J_B(\si'))}^{\vr_{\si'}(\si-1)}$$
along the preimage of $\ov\M_{1,\si}^{\si-1}$ under $\io_{\si-1,\si'}$ in
$$\id\!\times\!\pi_{\vr_{\si'}(\si-1)}\!:
\wt\cM_{1,(I_P(\si'),J_P(\si'))}
\!\times\! \wt\M_{0,(\ale_B(\si'),J_B(\si'))}^{\vr_{\si'}(\si-1)}
\lra
\wt\cM_{1,(I_P(\si'),J_P(\si'))}
\!\times\! \wt\M_{0,(\ale_B(\si'),J_B(\si'))}^0.$$
By~($I5$), this preimage~is 
$$\wt\cM_{1,(I_P(\si'),J_P(\si'))} \times
\Big(\bigcup_{\vr\in\A_B(\si';\si)}\!\!\!\!\!\!\!\!\wt\M_{0,\vr}^{\vr_{\si'}(\si-1)}\Big).$$
By the inductive assumption~($I5$) in Subsection~\ref{map0blconstr_subs} and
the second paragraph after Lemma~\ref{map1bl_lmm1b3}, 
$$\wt\M_{0,\vr_1}^{\vr_{\si'}(\si-1)} \cap \wt\M_{0,\vr_2}^{\vr_{\si'}(\si-1)}
=\eset
\qquad\forall\, \vr_1,\vr_2\!\in\!\A_B(\si';\si),\,  \vr_1\!\neq\!\vr_2.$$
Thus, by the construction of Subsection~\ref{map0blconstr_subs},
the blowup of $\wt\M_{0,(\ale_B(\si'),J_B(\si'))}^{\vr_{\si'}(\si-1)}$
along
$$\bigcup_{\vr\in\A_B(\si';\si)}\!\!\!\!\!\!\!\!\wt\M_{0,\vr}^{\vr_{\si'}(\si-1)}$$
is $\wt\M_{0,(\ale_B(\si'),J_B(\si'))}^{\vr_{\si'}(\si)}$, as needed
for the inductive statement~($I4$), with $\si\!-\!1$ replaced by~$\si$.
The inductive requirement~($I5$) is obtained by the same reasoning,
using the last statement of Lemma~\ref{ag_lmm2a}.\\

\noindent
Since $\ov\M_{1,\si}^{\si-1}$ is not contained in $\ov\M_{1,\si'}^{\si-1}$,
the bundle homomorphism
$$\io_{\si-1,\si'}^{\,*}T\ov\M_{1,\si}^{\si-1}\lra \N_{\io_{\si-1,\si'}}^{\ide}$$
must be surjective on every fiber over $\io_{\si-1,\si'}^{-1}(\ov\M_{1,\si}^{\si-1})$
by~($I14$).
Thus, the inductive assumption~($I6$), for $\si'\!<\!\si$, continues to hold. 
Furthermore, by~($I7$) and the last statement of Lemma~\ref{ag_lmm2a},
$\ov\cZ_{\si'}^{\si}$ is the closure of
\begin{equation*}\begin{split}
&\wt\cM_{1,(I_P(\si'),J_P(\si'))} \times
\Big(\wt\cD_{\vr_{\si'}(\si-1)}^{~-1}(0)
-\bigcup_{\stackrel{\vr\in\A_B(\ale_B(\si');d,J_B(\si'))}{\vr_{\si'}(\si-1)<\vr}}
\!\!\!\!\!\! \wt\M_{0,\vr}^{\vr_{\si'}(\si-1)}\Big)\\
&\qquad\qquad\qquad
= \wt\cM_{1,(I_P(\si'),J_P(\si'))} \times
\Big(\wt\cD_{\vr_{\si'}(\si-1)}^{~-1}(0)
-\bigcup_{\vr\in\A_B(\si';\si)}\!\!\!\!\!\! \wt\M_{0,\vr}^{\vr}
-\bigcup_{\stackrel{\vr\in\A_B(\ale_B(\si');d,J_B(\si'))}{\vr_{\si'}(\si)<\vr}}
\!\!\!\!\!\! \wt\M_{0,\vr}^{\vr_{\si'}(\si)}\Big)
\end{split}\end{equation*}
in $\wt\cM_{1,(I_P(\si'),J_P(\si'))}
\!\times\! \wt\M_{0,(\ale_B(\si'),J_B(\si'))}^{\vr_{\si'}(\si)}$.
By the construction of Subsection~\ref{map0blconstr_subs},
$$ \wt\cD_{\vr_{\si'}(\si-1)}
\big|_{\wt\M_{0,(\ale_B(\si'),J_B(\si'))}^{\vr_{\si'}(\si)}
-\bigcup_{\vr\in\A_B(\si';\si)}\! \wt\M_{0,\vr}^{\vr}}
=\wt\cD_{\vr_{\si'}(\si)}
\big|_{\wt\M_{0,(\ale_B(\si'),J_B(\si'))}^{\vr_{\si'}(\si)}
-\bigcup_{\vr\in\A_B(\si';\si)}\! \wt\M_{0,\vr}^{\vr}}.$$
Since $\wt\cD_{\vr_{\si'}(\si)}$ is transverse 
to the zero set outside of the subvarieties $\wt\M_{0,\vr}^{\vr_{\si'}(\si)}$
with $\vr\!>\!\vr_{\si'}(\si)$ by the inductive requirement~($I4$) in 
Subsection~\ref{map0blconstr_subs}, we conclude that the first part of
the inductive assumption ($I7$), with $\si\!-\!1$ replaced by~$\si$, is satisfied.
The second part follows from the last statement of Lemma~\ref{virimmer_lmm3}.\\

\noindent
It remains to verify the inductive assumption ($I9$)-($I12$),
with $\si\!-\!1$ replaced by~$\si$.
Suppose $\si^*\!\in\!\A_1(d,k)$ is such that $\si\!<\!\si^*$. 
If $\si\!\not\prec\!\si^*$, then
$$\rho_{\si^*}(\si)=\rho_{\si^*}(\si\!-\!1)
\qquad\hbox{and}\qquad
\ov\M_{1,\si^*}^{\si-1}\cap\ov\M_{1,\si}^{\si-1}=\eset,$$
by definition and ($I3$), respectively. It then follows that
\begin{gather*}
\io_{\si,\si^*}=\io_{\si-1,\si^*}, \qquad
\N_{\io_{\si,\si^*}}^{\ide}=\N_{\io_{\si-1,\si^*}}^{\ide},\\ 
\hbox{and}\qquad \ov\M_{1,\si^*}^{\si}\cap\ov\M_{1,\si'}^{\si}
=\ov\M_{1,\si^*}^{\si-1}\cap\ov\M_{1,\si'}^{\si-1}~~~\forall 
\si'\!\in\!\{(0)\}\!\sqcup\!\A_1(d,k).
\end{gather*}
Thus, the inductive assumptions ($I9$)-($I12$),
as stated above, imply the corresponding statements with $\si\!-\!1$ replaced by~$\si$.\\

\noindent
Suppose that $\si\!\prec\!\si^*$.
By~($I9$) and~(1) of Lemma~\ref{ag_lmm2a}, 
the domain of the immersion $\io_{\si,\si^*}$ induced
by~$\io_{\si-1,\si^*}$ is the blowup~of 
$$ \ov\cM_{1,(I_P(\si^*),J_P(\si^*))}^{\rho_{\si^*}(\si-1)}
\times \ov\M_{0,(\ale_B(\si^*),J_B(\si^*))}$$
along the preimage of $\ov\M_{1,\si}^{\si-1}$ under $\io_{\si-1,\si^*}$ in
$$\pi_{\rho_{\si^*}(\si-1)}\!\times\!\id\!:
\ov\cM_{1,(I_P(\si^*),J_P(\si^*))}^{\rho_{\si^*}(\si-1)}
\!\times\! \ov\M_{0,(\ale_B(\si^*),J_B(\si^*))}
\lra \ov\cM_{1,(I_P(\si^*),J_P(\si^*))}^0
\!\times\! \ov\M_{0,(\ale_B(\si^*),J_B(\si^*))}$$
By~($I10$), this preimage~is 
$$\Big(\bigcup_{\rho\in\A_P(\si^*;\si)}\!\!\!\!\!\!
\ov\cM_{1,\rho}^{\rho_{\si^*}(\si-1)}\Big) \times\ov\M_{0,(\ale_B(\si^*),J_B(\si^*))}.$$
By Lemma~\ref{curve1bl_lmm} and the paragraph after Lemma~\ref{map1bl_lmm1a}, 
$$\ov\cM_{1,\rho_1}^{\rho_{\si^*}(\si-1)}\cap\ov\cM_{1,\rho_2}^{\rho_{\si^*}(\si-1)}
=\eset
\qquad\forall\, \rho_1,\rho_2\!\in\!\A_P(\si^*;\si),\, 
\rho_1\!\neq\!\rho_2.$$
Thus, by the construction of Subsection~\ref{curve1bl_subs},
the blowup of $\ov\cM_{1,(I_P(\si^*),J_P(\si^*))}^{\rho_{\si^*}(\si-1)}$
along
$$\bigcup_{\rho\in\A_P(\si^*;\si)}\!\!\!\!\!\!
\ov\cM_{1,\rho}^{\rho_{\si^*}(\si-1)} $$
is $\ov\cM_{1,(I_P(\si^*),J_P(\si^*))}^{\rho_{\si^*}(\si-1)}$, as needed
for the inductive statement~($I9$), with $\si\!-\!1$ replaced by~$\si$.
The inductive assumptions~($I10$) and~($I11$) are verified similarly,
using the last statement of Lemma~\ref{ag_lmm2a} and Lemma~\ref{virimmer_lmm2}.
The argument for~($I11$) is nearly identical to the verification of
the inductive assumption~($I11$) in Subsection~\ref{map0blconstr_subs}.
Finally, the inductive requirement~($I12$), with $\si\!-\!1$ replaced by~$\si$,
follows from the last statement of Lemma~\ref{virimmer_lmm3}, 
along with the assumptions~($I1$) and~($I2$).\\

\noindent
We conclude this blowup construction after the $\si_{\max}$ step and put
$$\wt\M_{1,k}^0(\Pn,d)=\ov\M_{1,(0)}^{\si_{\max}}, \qquad
\ti\pi=\pi_{\si_{\max}}\big|_{\ov\M_{1,(0)}^{\si_{\max}}}, \quad\hbox{and}\quad
\wt\cZ_{\si}(\Pn,d)=\ov\cZ_{\si}^{\si_{\max}}.$$
The inductive assumptions ($I1$)-($I8$) imply that 
$$\ti\pi\!: \wt\M_{1,k}^0(\Pn,d) \lra \ov\M_{1,k}^0(\Pn,d)$$
is a desingularization as described in Subsection~\ref{descr_subs}.
By~($I3$), the final result of this blowup construction is independent
of the choice of full ordering~$<$ extending the natural partial ordering~$\prec$
on~$\A_1(d,k)$.

\section{Proof of Theorem~\ref{cone_thm}}
\label{cone_sec}

\subsection{Pushforwards of Vector Bundles}
\label{pushfor_subs}

\noindent
In this section we prove Theorem~\ref{cone_thm} by lifting the construction
of Section~\ref{map1bl_sec} from stable maps into $\Pn$ to stable maps into
(the total space of) the line bundle~$\cL$.\\

\noindent
Let $\tau\!:\cL\!\lra\!\Pn$ be the bundle projection map.
We denote by $\ov\M_{1,k}(\cL,d)$ the moduli space of degree-$d$ stable maps
from genus-one curves with $k$ marked points into~$\cL$.
The projection map $\tau$ induces a morphism,
$$p\!: \ov\M_{1,k}(\cL,d)\lra\ov\M_{1,k}(\Pn,d),
\qquad [\Si,u]\lra[\Si,\tau\circ u].$$
Since no fiber of $\cL$ contains the image of a non-constant holomorphic map,
the ghost components of $(\Si,\tau\circ u)$ are precisely the same as 
the ghost components of~$(\Si,u)$.
We note that
$$p^{-1}([\Si,u])=H^0(\Si;u^*\cL)\big/\Aut(\Si,u).$$
In particular, $p$ is a bundle of vector spaces, 
but of two possible ranks: $da$ and $da\!+\!1$.
Let $\S_{\cL}$ denote the sheaf of (holomorphic) sections of
$$p|_{\ov\M_{1,k}^0(\Pn,d)}\!: 
\ov\M_{1,k}(\cL,d)\big|_{\ov\M_{1,k}^0(\Pn,d)}\lra\ov\M_{1,k}^0(\Pn,d).$$
Similarly, denote by $\ti{S}_{\cL}$ the sheaf of sections of
$$p\!: \ti\pi^*\ov\M_{1,k}(\cL,d)\lra\wt\M_{1,k}^0(\Pn,d),$$
where $\ti\pi\!:\wt\M_{1,k}^0(\Pn,d)\!\lra\!\ov\M_{1,k}^0(\Pn,d)$
is the desingularization map of Theorem~\ref{main_thm}:\\

\begin{pspicture}(-3.5,-2.2)(10,.5)
\psset{unit=.4cm}
\rput(15,0){$\ov\M_{1,k}(\cL,d)$}\rput(15,-4){$\ov\M_{1,k}(\Pn,d)$}
\psline{->}(15,-1)(15,-3)\rput(15.6,-2){\small{$p$}}
\psline{->}(6,-1)(6,-3)\rput(6.6,-2){\small{$p$}}
\rput(6,0){$\wt\pi^*\ov\M_{1,k}(\cL,d)$}
\rput(6,-4){$\wt\M_{1,k}^0(\Pn,d)$}
\psline{->}(8.5,-4.1)(12.2,-4.1)\rput(10.2,-3.5){\smsize{$\ti\pi$}}
\psline{->}(8.9,-.1)(12.4,-.1)\rput(10.4,.5){\smsize{$\ti\pi$}}
\end{pspicture}

\begin{lmm}
\label{conesheaf_lmm}
With notation as in Theorem~\ref{cone_thm} and above,\\
${}\quad$ (1) the sheaves $\S_{\cL}$ and $\pi_*\ev_*\O_{\Pn}(a)$ over 
$\ov\M_{1,k}^0(\Pn,d)$ are isomorphic;\\
${}\quad$ (2) the sheaves $\ti{\S}_{\cL}$ and $\pi_*\ti\pi^*\ev_*\O_{\Pn}(a)$ over 
$\wt\M_{1,k}^0(\Pn,d)$ are isomorphic.
\end{lmm}

\noindent
Let $\U_{\cL}$ be the universal curve over $\ov\M_{1,k}(\cL,d)|_{\ov\M_{1,k}^0(\Pn,d)}$,
with structure map~$\pi_{\cL}$ and evaluation map~$\ev_{\cL}$.
The projection map~$\tau$ induces a morphism $\ti{p}_{\cL}$ on $\U_{\cL}$ so that 
the diagram\\

\begin{pspicture}(-1,-2.2)(10,1.2)
\psset{unit=.4cm}
\rput(10.2,0){$\U_{\cL}$}\rput(20,0){$\U$} 
\psline{->}(11,0)(19.2,0)\rput(15,.6){\small{$\ti{p}$}}
\rput(14,2){$\cL$}\psline{->}(10.5,0.5)(13.5,1.8)\rput{30}(11.5,1.5){\small{$\ev_{\cL}$}}
\rput(24,2){$\Pn$}\psline{->}(20.3,0.5)(23.3,1.8)\rput{25}(21,1.2){\small{$\ev$}}
\psline{->}(14.8,2)(23.2,2)\rput(19,2.5){\small{$\tau$}}
\rput(10,-4){$\ov\M_{1,k}(\cL,d)|_{\ov\M_{1,k}^0(\Pn,d)}$}\rput(20,-4){$\ov\M_{1,k}^0(\Pn,d)$}
\psline{->}(14.5,-4)(17.3,-4)\rput(16,-3.5){\small{$p$}}
\psline{->}(10,-1)(10,-3)\rput(9.3,-2){\small{$\pi_{\cL}$}}
\psline{->}(20,-1)(20,-3)\rput(20.5,-2){\small{$\pi$}}
\end{pspicture}

\noindent
commutes.
Suppose $W\!\subset\!\ov\M_{1,k}^0(\Pn,d)$ is an open subset.\\
(i) An element 
$$s\in \big\{\pi_*\ev_*\O_{\Pn}(a)\big\}(W)\equiv H^0\big(\pi^{-1}(W);\ev^*\cL\big)$$
induces a morphism $\ti{s}\!:\pi^{-1}(W)\!\lra\!\cL$  so that 
$\ev\!=\!\tau\!\circ\ti{s}$. 
In turn, $\ti{s}$ induces morphisms $f_s$ and $\ti{f}_s$ to 
$\ov\M_{1,k}(\cL,d)|_{\ov\M_{1,k}^0(\Pn,d)}$
and~$\U_{\cL}$,\\

\begin{pspicture}(-1,-2.2)(10,.3)
\psset{unit=.4cm}
\rput(10.2,0){$\pi^{-1}(W)$}\rput(20,0){$\U_{\cL}$} 
\psline{->}(12.2,0)(19.2,0)\rput(14.8,.6){\small{$\ti{f}_s$}}
\rput(10,-4){$W$}\rput(20,-4.2){$\ov\M_{1,k}(\cL,d)|_{\ov\M_{1,k}^0(\Pn,d)}$}
\psline{->}(11,-4)(15.3,-4)\rput(14,-3.5){\small{$f_s$}}
\psline{->}(10,-1)(10,-3)\rput(9.4,-2){\small{$\pi$}}
\psline{->}(20,-1)(20,-3)\rput(20.7,-2){\small{$\pi_{\cL}$}}
\rput(25,0){$\cL$} 
\psline{->}(20.8,0)(24.5,0)\rput(22.6,.6){\small{$\ev_{\cL}$}}
\end{pspicture}

\noindent
so that $\ti{s}\!=\!\ev_{\cL}\!\circ\ti{f}_s$.
Then,
$$\ev\circ\ti{p}\circ\ti{f}_s=\tau\circ\ev_{\cL}\!\circ\ti{f}_s
=\tau\circ\ti{s}=\ev\!: \pi^{-1}(W)\lra\Pn
\qquad\Lra\qquad p\circ f_{s}=\id_W,$$
since $\pi\circ\ti{p}\circ\ti{f}_s=p\circ f_{s}\!\circ\pi$.
Thus, $f_s\!\in\!\S_{\cL}(W)$.
It is immediate that the~map
$$\big\{\pi_*\ev_*\O_{\Pn}(a)\big\}(W)\lra\S_{\cL}(W),
\qquad s\lra f_s,$$
induces a sheaf homomorphism.\\
(ii) Conversely, let $\si\!\in\!\S_{\cL}(W)$, 
i.e.~$\si\!:W\!\lra\!\ov\M_{1,k}(\cL,d)$ is a morphism such that $p\circ\si=\id_W$.
Since $\U_{\cL}\!=\!p^*\U$,
$$\pi^{-1}(W)\equiv\U|_W=\si^*\U_{\cL}.$$
Thus, $\si$ lifts to a morphism
$$\ti\si\!: \pi^{-1}(W)=\si^*\U_{\cL}\lra \U_{\cL}.$$
Let $g_{\si}\!=\!\ev_{\cL}\!\circ\ti\si$. Then,
$$\tau\circ g_{\si}=\tau\circ\ev_{\cL}\!\circ\ti\si
=\ev\circ\ti{p}\circ\ti\si=\ev,$$
i.e.~$g_{\si}\!\in\!H^0(\pi^{-1}(W);\ev^*\cL)$.
It is immediate that the~map
$$\S_{\cL}(W)\lra\big\{\pi_*\ev_*\O_{\Pn}(a)\big\}(W),
\qquad \si\lra g_{\si},$$
induces a sheaf homomorphism.
Furthermore, 
$$g_{f_s}=s \quad\forall\,s\in\big\{\pi_*\ev_*\O_{\Pn}(a)\big\}(W)
\qquad\hbox{and}\qquad
f_{g_{\si}}=\si \quad\forall\,\si\in\S_{\cL}(W).$$
These observations imply the first statement of Lemma~\ref{conesheaf_lmm}.
The second claim is proved similarly.\\

\noindent
Let
$$\ov\M_{1,k}^0(\cL,d)\subset \ov\M_{1,k}(\cL,d)$$ 
be the closure of the locus of maps from smooth domains.
We show in Subsection~\ref{conebl_subs} that the proper transform 
$\wt\M_{1,k}^0(\cL,d)$ of $\ov\M_{1,k}^0(\cL,d)$ in
$$\ti\pi^*\ov\M_{1,k}(\cL,d)\lra \wt\M_{1,k}^0(\Pn,d)$$
is smooth.
Similarly to the case of $\wt\M_{1,k}^0(\Pn,d)$, the main stratum of 
$\wt\M_{1,k}^0(\cL,d)$,
$$\M_{1,k}^{\eff}(\cL,d)\equiv 
\ov\M_{1,k}(\cL,d)\big|_{\M_{1,k}^{\eff}(\Pn,d)}
=\wt\M_{1,k}^0(\cL,d)-
\bigcup_{\si\in\A_1(d,k)}\!\!\!\!\!\!p^{-1}\big(\Im\,\io_{\si_{\max},\si}\big),$$
is smooth.
On the other hand, by the inductive assumption ($I1$) and the last paragraph 
of Subsection~\ref{conebl_subs}, for each $\si\!\in\!\A_1(d,k)$
$$\wt\M_{1,k}^0(\cL,d)\cap p^{-1}\big(\Im\,\io_{\si_{\max},\si}\big)$$
is the image of a smooth variety under the bundle homomorphism $\bar{j}_{\si_{\max},\si}$
lifting the embedding~$\bar\io_{\si_{\max},\si}$ of Subsection~\ref{map1blconstr_subs}.
Thus, 
$$\wt\M_{1,k}^0(\cL,d)\cap p^{-1}\big(\Im\,\io_{\si_{\max},\si}\big)$$
is a smooth subvariety of $\wt\M_{1,k}^0(\cL,d)$.
As its normal cone in $\wt\M_{1,k}^0(\cL,d)$ is a line bundle by the inductive assumption ($I1$) 
of Subsection~\ref{conebl_subs} for every $\si\!\in\!\A_1(d,k)$,
we conclude that the entire space $\wt\M_{1,k}^0(\cL,d)$ is smooth.
Furthermore, the fibers of 
$$\ti{p}\!: \wt\M_{1,k}^0(\cL,d)\lra \wt\M_{1,k}^0(\Pn,d)$$
are vector spaces of the same rank and $\wt\M_{1,k}^0(\cL,d)$ contains
$\wt\M_{1,k}^0(\Pn,d)$ as the zero section.
Thus, $\ti{p}$ is a vector bundle.\\

\noindent
Lemma~\ref{conesheaf_lmm} and the previous paragraph imply (1) of 
Theorem~\ref{cone_thm}.
The second claim of this theorem is obtained in the last paragraph of
Subsection~\ref{conebl_subs}.
Finally, (3)~of Theorem~\ref{cone_thm} follows from~(4) of Theorem~\ref{main_thm}
and the following lemma.

\begin{lmm}
\label{pushfor_lmm3}
Suppose $\ti\pi\!:\wt\M\!\lra\!\ov\M$ is a morphism between varieties, 
$\U\!\lra\!\ov\M$ is a flat family of curves, $\cL\!\lra\!\U$
is a line bundle, and $\pi\!:\ti\U\!\lra\!\wt\M$ and $\wt\cL\!\lra\!\wt\U$
are the pullbacks of $\U$ and $\cL$ via $\ti\pi$:\\

\begin{pspicture}(-3,-2.2)(10,1)
\psset{unit=.4cm}
\rput(12.3,0){$\U$}
\rput(15.7,2){$\cL$}\psline{<-}(12.8,0.5)(15.2,1.8)\rput(12.3,-4){$\ov\M$}
\psline{->}(12.3,-1)(12.3,-3)\rput(12.9,-2){\small{$\pi$}}
\psline{->}(6.8,-1)(6.8,-3)\rput(7.4,-2){\small{$\pi$}}
\rput(7.8,0){$\wt\U\!=\!\wt\pi^*\U$}
\rput(11.2,2){$\wt\cL\!=\!\ti\pi^*\cL$}\psline{<-}(7,0.5)(9.4,1.8)
\rput(6.8,-4){$\wt\M$}
\psline{->}(7.8,-4.1)(11.3,-4.1)\rput(9.6,-3.5){\smsize{$\ti\pi$}}
\psline{->}(9.8,-.1)(11.6,-.1)\rput(10.7,.5){\smsize{$\ti\pi$}}
\end{pspicture}

\noindent
If the morphism $\ti\pi$ is surjective and its fibers are compact and connected,
then
$$\ti\pi_*\pi_*\wt\cL=\pi_*\cL.$$\\
\end{lmm}

\noindent
Since $\cL$ is locally trivial, Lemma~\ref{pushfor_lmm3} follows from
$$\ti\pi_*\O_{\wt\M}=\O_{\ov\M}.$$
In turn, this identity follows from the fact that every holomorphic function
on a compact connected variety is constant.
Thus, if $W\!\subset\!\ov\M$ is any open subset and $\ti{f}$ is a holomorphic function on 
$\ti\pi^{-1}(W)\!\subset\!\wt\M$, then $\ti{f}$ is constant on the fibers of $\ti\pi$,
i.e.~$\ti{f}\!=\!\ti\pi^*f$ for some holomorphic function $f$ on~$W$.

\subsection{Construction of Bundle Homomorphism}
\label{conehomomor_subs}

\noindent
In this subsection we describe the surjective bundle homomorphism 
that appears in the second statement of Theorem~\ref{cone_thm};
see Proposition~\ref{conehomomor_prp}.
The construction of this homomorphism is similar to 
the construction of the homomorphism
$\wt\cD_{(\ale,J)}$ in Subsections~\ref{map0prelim_subs} and~\ref{map0blconstr_subs}.\\

\noindent
Let $\cL\!\lra\!\Pn$ be a line bundle as in Subsection~\ref{descr_subs}.
If $J$ is a finite set, let
$$\V_0=\ov\M_{0,\{0\}\sqcup J}(\cL,d) \lra \ov\M_{0,\{0\}\sqcup J}(\Pn,d)$$
be the corresponding cone.
In particular, if $[\Si,u]\!\in\!\ov\M_{0,\{0\}\sqcup J}(\Pn,d)$, then
$$\V_0\big|_{[\Si,u]}=H^0(\Si;u^*\cL)\big/\Aut(\Si,u).$$
In this, genus-{\it zero}, case, this is a vector bundle of the expected rank.
Let
$$\na^u\!: \Ga(\Si;u^*\cL)\lra\Ga(\Si;T^*\Si\!\otimes\!u^*\cL)$$
be the pullback of the standard Hermitian connection in $\cL$ by $u$.
We define
\begin{gather*}
\begin{split}
\D_0&\in\Ga\big(\ov\M_{0,\{0\}\sqcup J}(\Pn,d);
\Hom(L_0\!\otimes\!\V_0,\ev_0^*\cL)\big)\\
&\quad=\Ga\big(\ov\M_{0,\{0\}\sqcup J}(\Pn,d);
\Hom(L_0,\Hom(\V_0,\ev_0^*\cL))\big)\\
&\quad=\Ga\big(\ov\M_{0,\{0\}\sqcup J}(\Pn,d);
\Hom(\V_0,\Hom(L_0,\ev_0^*\cL))\big)
\end{split}\\ 
\hbox{by}\qquad \D_0\xi = \na^u\xi|_{x_0(\Si,u)} \quad\forall\,\xi\!\in\!H^0(\Si;u^*\cL),
\end{gather*}
where $x_0(\Si,u)\!\in\!\Si$ is the marked point labeled by~$0$ as before.
We note that $\D_0$ vanishes identically on the subvarieties $\ov\M_{0,\si}(\Pn,d)$ 
with $\si\!\in\!\A_0(d,J)$ defined in Subsection~\ref{map0str_subs}.\\

\noindent
If $\aleph$ and $J$ are finite sets, let
$$p_{(\aleph,J)}\!:\V_{(\aleph,J)} \lra \ov\M_{0,(\aleph,J)}(\Pn,d)$$
be the vector bundle induced by $\cL$, where $\ov\M_{0,(\aleph,J)}(\Pn,d)$
is as in Subsection~\ref{map0prelim_subs}.
It is immediate that
$$\V_{(\ale,J)} = \big\{
(\xi_i)_{i\in\ale}\!\in\!\bigoplus_{i\in\ale}\pi_i^*\V_0\!:
\ev_0(\xi_i)\!=\!\ev_0(\xi_{i'})~\forall\,i,i'\!\in\!\ale\big\}
=\ov\M_{0,(\ale,J)}(\cL,d).$$
Note that for every $\si\!=\!(m;J_P,J_B)\!\in\!\A_0(d,J)$,
$$\io_{\si}^*\V_0\!=\!\pi_B^*\V_{([m],J_B)}
\lra\ov\cM_{0,\{0\}\sqcup[m]\sqcup J_B} \!\times\!\ov\M_{0,([m],J)}(\Pn,d),$$
where $\io_{\si}$ is as in Subsection~\ref{map0str_subs}.

\begin{lmm}
\label{conederiv_lmm1}
If $d\!\in\!\Z^+$, $J$, $\cL$, and $\V_0$ are as above,
the bundle homomorphism 
$$\D_0\in\Ga\big(\ov\M_{0,\{0\}\sqcup J}(\Pn,d);
\Hom(\V_0,L_0^*\!\otimes\!\ev_0^*\cL)\big)$$
is surjective on the complement of 
the subvarieties $\ov\M_{0,\si}(\Pn,d)$ with  $\si\!\in\!\A_0(d,J)$.
Furthermore, for every 
$$\si\!\equiv\!(m;J_P,J_B) \in \A_0(d,J),$$
the differential of~$\D_0$,
$$\na\D_0\!: \N_{\io_{\si}} \lra \io_{\si}^*\,\Hom(\V_0,L_0^*\!\otimes\!\ev_0^*\cL)
=\pi_P^*L_0^*\!\otimes\!\pi_B^*\Hom(\V_{([m],J_B)},\ev_0^*\cL),$$
in the normal direction to the immersion $\io_{\si}$ is given~by
$$ \na\D_0\big|_{\pi_P^*L_i\otimes\pi_B^*\pi_i^*L_0}
= \pi_P^*s_i\!\otimes\!\pi_B^*\pi_i^*\D_0 \qquad\forall\,i\in[m],$$
where $s_i$ is the homomorphism defined in Subsection~\ref{curvebldata_subs}.
\end{lmm}

\noindent
Lemma~\ref{conederiv_lmm1} can viewed as the analogue of Lemma~\ref{deriv0str_lmm} 
for vector bundle sections.
The first claim of Lemma~\ref{conederiv_lmm1} is an immediate consequence of the fact that
$$H^1\big(\Si;u^*\cL\!\otimes\!\O(-2z)\big)=\{0\}$$
for every genus-zero stable map $(\Si,u)$ and a smooth point $z\!\in\!\Si$
such that the restriction of $u$ to the irreducible component of $\Si$
containing $z$ is not constant.
The second statement follows from \cite[Lemma~\ref{g1cone-derivest_lmm}]{g1cone}.\\

\noindent
With notation as in Subsection~\ref{map0prelim_subs}, let 
\begin{equation*}\begin{split}
\D_{(\ale,J)} &\in \Ga\big(\ov\M_{0,(\ale,J)}(\Pn,d);
\Hom(\V_{(\ale,J)},\Hom(F_{(\ale,J)},\ev_0^*\cL))\big)\\
&\quad = \Ga\big(\ov\M_{0,(\ale,J)}(\Pn,d);
\Hom(F_{(\ale,J)},\Hom(\V_{(\ale,J)},\ev_0^*\cL))\big)\\
&\quad =\Ga\big(\ov\M_{0,(\ale,J)}(\Pn,d);
\Hom(F_{(\ale,J)}\!\otimes\!\V_{(\ale,J)},\ev_0^*\cL)\big)
\end{split}\end{equation*}
be the homomorphism defined by
$$\D_{(\ale,J)}\big|_{\pi_i^*L_0\otimes\pi_j^*\V_0}=
\begin{cases}
\pi_i^*\D_0,& 
\hbox{if}~j\!=\!i;\\
0,&\hbox{otherwise};
\end{cases}
\qquad\forall\, i,j\!\in\!\ale.$$\\
It induces a section
\begin{equation*}\begin{split}
\wt\D_0 \in& \Ga\big(\wt\M_{0,(\ale,J)}^0(\Pn,d);
\Hom(\ga_{(\ale,J)},\pi_{\P F_{(\ale,J)}}^*\Hom(\V_{(\ale,J)},\ev_0^*\cL))\big)\\
& =\Ga\big(\wt\M_{0,(\ale,J)}^0(\Pn,d); \Hom(\pi_{\P F_{(\ale,J)}}^*\V_{(\ale,J)},
\E_0^*\!\otimes\!\pi_{F_{(\ale,J)}}^*\ev_0^*\cL)\big).
\end{split}\end{equation*}
This section vanishes identically on the subvarieties 
$\wt\M_{0,\vr}^0(\Pn,d)$ of $\wt\M_{0,(\ale,J)}^0(\Pn,d)$ 
with $\vr\!\in\!\A_0(\ale;d,J)$, defined in Subsection~\ref{map0prelim_subs}.

\begin{lmm}
\label{conederiv_lmm2}
The bundle homomorphism 
$$\wt\D_0\in\Ga\big(\wt\M_{0,(\ale,J)}^0(\Pn,d); \Hom(\pi_{\P F_{(\ale,J)}}^*\V_{(\ale,J)},
\E_0^*\!\otimes\!\pi_{F_{(\ale,J)}}^*\ev_0^*\cL)\big)$$
is surjective on the complement of the subvarieties 
$\wt\M_{0,\vr^*}^0(\Pn,d)$ with $\vr^*\!\in\!\A_0(\ale;d,J)$.
Furthermore, for every $\vr^*\!\in\!\A_0(\ale;d,J)$ as in Lemma~\ref{map0bl_lmm1},
the differential of~$\wt\D_0$,
\begin{equation*}\begin{split}
\na\wt\D_0\!: \N_{\io_{0,\vr^*}} \lra&  \io_{0,\vr^*}^{\,*}
\Hom\big(\pi_{\P F_{(\ale,J)}}^*\V_{(\ale,J)},
\E_0^*\!\otimes\!\pi_{F_{(\ale,J)}}^*\ev_0^*\cL\big)\\
& =\pi_P^*\E_0^*\!\otimes\!
\pi_B^*\Hom(\V_{(\ale_B(\vr^*),J_B(\vr^*))},\ev_0^*\cL),
\end{split}\end{equation*}
in the normal direction to the immersion $\io_{0,\vr^*}$ is given~by
\begin{gather*}
\na\wt\D_0\big|_{\pi_P^*L_{0,(l,i)}\otimes\pi_B^*\pi_{(l,i)}^*L_0} 
= \pi_P^*s_{0,(l,i)}\!\otimes\!\pi_B^*\pi_{(l,i)}^*\D_0
\qquad\forall\,i\!\in\![m_l^*],\, l\!\in\!\ale_P(\vr^*),\\
\hbox{and}\qquad
\na\wt\D_0\big|_{\N_{\io_{0,\vr^*}}^{\top}}
=\pi_P^*\id\!\otimes\!\pi_B^*\D_{(\ale_B(\vr^*),J_B(\vr^*))},
\end{gather*}
where $s_{0,(l,i)}$ is the homomorphism defined in Subsection~\ref{curve0bl_subs}.
\end{lmm}

\noindent
This lemma follows immediately from Lemma~\ref{conederiv_lmm1}.

\begin{prp}
\label{conehomomor_prp}
With notation as above, there exists a surjective bundle homomorphism
$$\wt\D_{(\ale,J)} \in \Ga\big(\wt\M_{0,(\ale,J)}(\Pn,d);
\Hom(\pi_{0,(\ale,J)}^*\pi_{\P F_{(\ale,J)}}^*\V_{(\ale,J)},\wt\E^*\!\otimes\!
\pi_{0,(\ale,J)}^*\pi_{\P F_{(\ale,J)}}^*\ev_0^*\cL)\big)$$
such that 
\begin{gather*}
\wt\D_{(\ale,J)}\big|_{\P F_{(\ale,J)}^0}
=\wt\D_0\big|_{\P F_{(\ale,J)}^0}, \qquad\hbox{where}\\
\P F_{(\ale,J)}^0 =\P F_{(\ale,J)}-
\bigcup_{\vr\in\A_0(\ale;d,J)}\!\!\!\!\!\!\!\!\wt\M_{0,\vr}^0(\Pn,d)
\subset \wt\M_{0,(\ale,J)}^0(\Pn,d),\wt\M_{0,(\ale,J)}(\Pn,d).
\end{gather*}\\
\end{prp}

\noindent
In fact, in the notation of Subsection~\ref{map0blconstr_subs}, for every
$\vr\!\in\!\{0\}\!\sqcup\!\A_0(\ale;d,J)$ there exists a bundle homomorphism
$$\wt\D_{\vr} \in \Ga\big(\wt\M_{0,(\ale,J)}^{\vr};
\Hom(\pi_{\vr}^*\pi_{\P F_{(\ale,J)}}^*\V_{(\ale,J)},\E_{\vr}^*\!\otimes\!
\pi_{\vr}^*\pi_{\P F_{(\ale,J)}}^*\ev_0^*\cL)\big)$$
such that\\
${}\quad$ (i) the restrictions of $\wt\D_{\vr}$ and $\wt\D_0$ to
$\P F_{(\ale,J)}^0$ agree;\\
${}\quad$ (ii) $\wt\D_{\vr}$ is surjective outside of the subvarieties 
$\wt\M_{0,\vr^*}^{\vr}$ with $\vr^*\!>\!\vr$;\\
${}\quad$ (iii) $\wt\D_{\vr}$ vanishes identically on the subvarieties 
$\wt\M_{0,\vr^*}^{\vr}$ with $\vr^*\!>\!\vr$;\\
${}\quad$ (iv) for each $\vr^*\!>\!\vr$, the differential of $\wt\D_{\vr}$
in the normal direction to the immersion $\io_{\vr,\vr^*}$\\
${}\qquad\quad$ is given as in the statement of Lemma~\ref{conederiv_lmm2}, 
but with $s_{0,(l,i)}$ replaced by~$s_{\rho_{\vr^*}(\vr),(l,i)}$.\\
Similarly to the construction of the bundle sections $\wt\cD_{\vr}$ in
Subsection~\ref{map0blconstr_subs}, we construct the bundle homomorphisms~$\wt\D_{\vr}$
inductively starting with $\wt\D_0$ and twisting by the exceptional divisor at each step.
The inductive assumptions (i)-(iv) are analogous to ($I3$), ($I4$), and~($I12$)
in Subsection~\ref{map0blconstr_subs} and are verified similarly.
Of course, we take
$$\wt\D_{(\ale,J)}=\wt\D_{\vr_{\max}}.$$

\subsection{Structure of the Cone $\V_{1,k}^d$}
\label{conestr_subs}

\noindent
In this subsection we describe the structure of the cone
$$p_0\!:\ov\M_{1,k}(\cL,d)\lra\ov\M_{1,k}(\Pn,d),$$ 
restating the primary structural result of~\cite{g1cone}.\\

\noindent
For each element $\si\!=\!(m;J_P,J_B)$ of $\A_1(d,k)$, let
$$\V_{1,\si}^0 \!\equiv\! \ov\M_{1,\si}(\cL,d)
=p_0^{-1}\big(\ov\M_{1,\si}^0\big)\equiv p_0^{-1}\big(\ov\M_{1,\si}(\Pn,d)\big)
\subset \V_{1,k}^0 \!\equiv\! \ov\M_{1,k}(\cL,d).$$
The subvarieties\footnote{typo fixed} $\ov\M_{1,\si}(\cL,d)$ of $\ov\M_{1,k}(\cL,d)$
can also be defined analogously to the subvarieties $\ov\M_{1,\si}(\Pn,d)$ of 
$\ov\M_{1,k}(\Pn,d)$; see the beginning of Subsection~\ref{descr_subs}.
Similarly to Subsection~\ref{map1prelim_subs}, let
$$j_{0,\si}\!: \ov\cM_{1,(I_P(\si),J_P(\si))}^0 \times \V_{(\ale_B(\si),J_B(\si))} \lra 
\V_{1,\si}^0\subset\V_{1,k}^0$$
be the natural node-identifying immersion so that the diagram\\

\begin{pspicture}(-1.1,-1.8)(10,.5)
\psset{unit=.4cm}
\rput(10,0){$\ov\cM_{1,(I_P(\si),J_P(\si))}^0~\times~\V_{(\ale_B(\si),J_B(\si))}$}
\rput(25,0){$\V_{1,\si}^0~\subset~\V_{1,k}^0$}
\rput(10.4,-4){$\ov\cM_{1,(I_P(\si),J_P(\si))}^0~~\times~~\ov\M_{0,(\ale_B(\si),J_B(\si))}$}
\rput(25,-4){$\ov\M_{1,\si}^0~\subset~\ov\M_{1,k}^0$}
\psline{->}(17.2,0)(22,0)\psline{->}(18.5,-4)(21.5,-4)
\psline{->}(6,-1)(6,-3.2)\psline{->}(14,-1)(14,-3.2)\psline{->}(22.7,-1)(22.7,-3)
\rput(19.5,.7){\small{$j_{0,\si}$}}\rput(20,-3.4){\small{$\io_{0,\si}$}}
\rput(6.6,-2){\small{$\id$}}\rput(15,-2){\small{$p_{0,\si}$}}
\rput(23.4,-2){\small{$p_0$}}
\end{pspicture}

\noindent
commutes. 

\begin{lmm}
\label{cone1bl_lmm1}
If $d,n\!\in\!\Z^+$ and $k\!\in\!\bar\Z^+$, the collection
$\{j_{0,\si}\}_{\si\in\A_1(d,k)}$
of immersions is properly self-intersecting.
For every $\si\!\in\!\A_1(d,k)$,
$$\N_{j_{0,\si}}^{\ide}=\big\{\id\!\times\!p_{0,\si}\big\}^*
\N_{\io_{0,\si}}^{\ide}$$
is an idealized normal bundle for $j_{0,\si}$.
\end{lmm}

\noindent
The differential $dp_0$ of $p_0$ induces a surjective linear map
$$\Im\, dj_{0,\si}\lra \Im\,d\io_{0,\si}.$$
Since the fibers of $p_0$ are vector spaces, it follows that $dp_0$ induces an injection
$$j_{0,\si}^*TC\V_{1,k}^0\big/\Im\, dj_{0,\si}\lra
\io_{0,\si}^*TC\ov\M_{1,k}^0\big/\Im\, d\io_{0,\si}.$$
Thus, Lemma~\ref{cone1bl_lmm1} follows from Lemma~\ref{map1bl_lmm1b1}.\\

\noindent
We denote by $\V_{1,(0)}^0$ the main component $\ov\M_{1,k}^0(\cL,d)$
of the moduli space $\ov\M_{1,k}(\cL,d)$.
If $\si\!\in\!\A_1(d,k)$, we~put
$$\W_{\si}^0= j_{0,\si}^{~-1} \big(\V_{1,(0)}^0\big)
\equiv j_{0,\si}^{~-1} \big( \V_{1,(0)}^0\!\cap\V_{1,\si}^0\big).$$
Note that
$$\big\{\id\!\times\!p_{0,\si}\big\}\big(\W_{\si}^0\big)
=\bar\cZ_{\si}^0\equiv\io_{0,\si}^{-1}\big(\ov\M_{1,(0)}^0\big).$$
Let  $\N\W_{\si}^0\!\subset\!\N_{j_{0,\si}}^{\ide}$ be
the normal cone $\N_{j_{0,\si}|\V_{1,(0)}^0}$ for 
$j_{0,\si}|_{\W_{\si}^0}$ in $\V_{1,(0)}^0$.
Its structure is described in Lemma~\ref{cone1bl_lmm2} below.
Let
\begin{equation*}\begin{split}
\D_{0,\si} &\in  \Ga\big(\ov\cM_{1,(I_P(\si),J_P(\si))}^0 \!\times\!  
\ov\M_{0,(\ale_B(\si),J_B(\si))};
\Hom(\pi_B^*\V_{(\ale_B(\si),J_B(\si))},\Hom(\N_{\io_{0,\si}}^{\ide},
\pi_P^*\E_0^*\!\otimes\!\pi_B^*\ev_0^*\cL))\big)\\
&\quad=\Ga\big(\ov\cM_{1,(I_P(\si),J_P(\si))}^0 \!\times\!  
\ov\M_{0,(\ale_B(\si),J_B(\si))};
\Hom(\N_{\io_{0,\si}}^{\ide},\pi_P^*\E_0^*\!\otimes\!
\pi_B^*\Hom(\V_{(\ale_B(\si),J_B(\si))},\ev_0^*\cL))\big)
\end{split}\end{equation*}
be the section defined by
$$ \D_{0,\si}\big|_{\pi_P^*L_i\otimes\pi_B^*\pi_i^*L_0}
=\pi_P^*s_i\!\otimes\!\pi_B^*\pi_i^*\D_0, 
\qquad\forall\,i\!\in\![m],$$
where $s_i$ and $\D_0$ are as in 
Subsections~\ref{curvebldata_subs} and~\ref{conehomomor_subs}, respectively.
If $\xi\!\in\!\pi_B^*\V_{(\ale_B(\si),J_B(\si))}$, we will view $\D_{0,\si}\xi$
as a homomorphism
$$\D_{0,\si}\xi\!: 
\N_{j_{0,\si}}^{\ide}\big|_{\xi}\!=\!
\N_{\io_{0,\si}}^{\ide}\big|_{\{\id\times p_{0,\si}\}(\xi)}
\lra \pi_P^*\E_0^*\!\otimes\!\pi_B^*\ev_0^*\cL
\big|_{\{\id\times p_{0,\si}\}(\xi)}.$$

\begin{lmm}
\label{cone1bl_lmm2}
For all $\si\!\in\!\A_1(d,k)$, $\V_{1,(0)}^0$ is a proper subvariety
of $\V_{1,k}^0$ relative to the immersion~$j_{0,\si}$.
The homomorphism 
$$\N\W_{\si}^0 \lra \{\id\!\times\!p_{0,\si}\}^*\N\bar\cZ_{\si}^0$$
induced by $dp_0$ is injective.
Furthermore, 
\begin{gather*}
\W_{\si}^0\big|_{\cZ_{\si}^0}=\big\{
\xi\!\in\!\pi_B^*\V_{(\ale_B(\si),J_B(\si))}|_{\cZ_{\si}^0}\!:
\ker\,\{\D_{0,\si}\xi\}
\big|_{\N\bar\cZ_{\si}^0|_{\{\id\times p_{0,\si}\}(\xi)}}\!\neq\!\{0\}
\big\}\\
\hbox{and}\qquad
\N\W_{\si}^0\big|_{\xi}=
\ker\big\{\D_{0,\si}\xi\big\}
\big|_{\N\bar\cZ_{\si}^0|_{\{\id\times p_{0,\si}\}(\xi)}}
\subset\N_{j_{0,\si}}^{\ide} \qquad\forall\,\xi\!\in\!\W_{\si}^0\big|_{\cZ_{\si}^0}.
\end{gather*}
Finally, $\W_{\si}^0$ is the closure of $\W_{\si}^0|_{\cZ_{\si}^0}$ in
$\ov\cM_{1,(I_P(\si),J_P(\si))}^0\!\times\!\V_{(\ale_B(\si),J_B(\si))}$ 
and  $\N\W_{\si}^0$ is the closure of
$\N\W_{\si}^0\big|_{\W_{\si}^0|_{\cZ_{\si}^0}}$ in~$\N_{j_{0,\si}}^{\ide}$.
\end{lmm}

\noindent
Since the fibers of $p_0$ are vector spaces, 
the first two sentences of this lemma follow from Lemma~\ref{map1bl_lmm2}.
The middle claim of Lemma~\ref{cone1bl_lmm2} is a restatement of
\cite[Lemma~\ref{g1cone-g1conebdstr_lmm}]{g1cone}.
The remaining claims of the lemma follow
from \cite[Lemma~\ref{g1cone-g1conebdstr_lmm}]{g1cone} by dimension counting,
similarly to the argument following Lemma~\ref{map1bl_lmm2}.\\

\noindent
{\it Remark:} It may appear that the statement of Lemma~\ref{cone1bl_lmm2} 
depends on the choice of a hermitian connection (or metric) in the line bundle $\cL\!\lra\!\Pn$.
As explained in detail in \cite[Subsect.~\ref{g1cone-g1conelocalstr_subs2}]{g1cone}, 
the dependence is only on the holomorphic structure of~$\cL$, as the case should~be.

\subsection{Desingularization Construction}
\label{conebl_subs}

\noindent
In this subsection we lift the inductive blowup construction of 
Subsection~\ref{map1blconstr_subs} to the cone 
$$p_0\!:\V_{1,k}^0\lra\ov\M_{1,k}^0.$$
For each $\si\!\in\!\A_1(d,k)$, let $\V_{1,k}^{\si}\!\equiv\!\pi_{\si}^*\V_{1,k}^0$ 
be the pullback of $\V_{1,k}^0$ to~$\ov\M_{1,k}^{\si}$:\\

\begin{pspicture}(-1.5,-2)(10,.3)
\psset{unit=.4cm}
\rput(10,0){$\V_{1,k}^{\si}\!\equiv\!\pi_{\si}^*\V_{1,k}^0$}\rput(18,0){$\V_{1,k}^0$}
\rput(10,-4){$\ov\M_{1,k}^{\si}$}\rput(18,-4){$\ov\M_{1,k}^0$}
\psline{->}(13,0)(17,0)\psline{->}(11.5,-4)(16.5,-4)
\psline{->}(10,-.8)(10,-3)\psline{->}(18,-.8)(18,-3)
\rput(15,.6){\small{$\pi_{\si}$}}\rput(14,-3.5){\small{$\pi_{\si}$}}
\rput(10.7,-2){\small{$p_{\si}$}}\rput(18.7,-2){\small{$p_0$}}
\end{pspicture}

\noindent
For each $\si'\!\in\!\A_1(d,k)$, let
$$\V_{1,\si'}^{\si}=\V_{1,k}^{\si}\big|_{\ov\M_{1,\si'}^{\si}}
=\pi_{\si}^{-1}\big(\V_{1,k}^0\big|_{\ov\M_{1,\si'}^0}\big).$$
The bundle homomorphisms $j_{0,\si'}$ lift to bundle homomorphisms onto $\V_{1,\si'}^{\si}$
covering the immersion~$\io_{\si,\si'}$ of Subsection~\ref{map1blconstr_subs}:\\

\begin{pspicture}(-3,-2.4)(10,.5)
\psset{unit=.4cm}
\rput(8,0){$\wt\cM_{1,(I_P(\si'),J_P(\si'))}\times
\pi_{\vr_{\si'}(\si)}^*\pi_{\P F_{(\ale,J)}}^*\V_{(\ale_B(\si'),J_B(\si'))}$}
\rput(25,0){$\V_{1,\si'}^{\si}~\subset~\V_{1,k}^{\si}$}
\rput(6.5,-4){$\wt\cM_{1,(I_P(\si'),J_P(\si'))}\times
~~~~~\wt\M_{0,(\ale_B(\si'),J_B(\si'))}^{\vr_{\si'}(\si)}$}
\rput(25,-4){$\ov\M_{1,\si'}^{\si}~\subset~\ov\M_{1,k}^{\si}$}
\psline{->}(18,0)(22,0)\psline{->}(15.3,-4)(21.2,-4)
\psline{->}(2,-1)(2,-3.2)\psline{->}(12,-1)(12,-3.2)\psline{->}(22.7,-1)(22.7,-3)
\rput(20,.7){\small{$j_{\si,\si'}$}}\rput(18,-3.4){\small{$\io_{\si,\si'}$}}
\rput(2.6,-2){\small{$\id$}}\rput(13.1,-2){\small{$p_{\si,\si'}$}}
\rput(23.4,-2){\small{$p_{\si}$}}
\rput(-6,-1.5){$\si'\!\le\!\si\!:$}
\end{pspicture}

\begin{pspicture}(-3,-2.5)(10,.5)
\psset{unit=.4cm}
\rput(8,0){$\ov\cM_{1,(I_P(\si'),J_P(\si'))}^{\rho_{\si'}(\si)}
\times\V_{(\ale_B(\si'),J_B(\si'))}$}
\rput(25,0){$\V_{1,\si'}^{\si}~\subset~\V_{1,k}^{\si}$}
\rput(8.6,-4){$\ov\cM_{1,(I_P(\si'),J_P(\si'))}^{\rho_{\si'}(\si)}\times
\ov\M_{0,(\ale_B(\si'),J_B(\si'))}$}
\rput(25,-4){$\ov\M_{1,\si'}^{\si}~\subset~\ov\M_{1,k}^{\si}$}
\psline{->}(15.3,0)(22,0)\psline{->}(16.3,-4)(21.2,-4)
\psline{->}(5,-1)(5,-3.2)\psline{->}(12,-1)(12,-3.2)\psline{->}(22.7,-1)(22.7,-3)
\rput(19,.7){\small{$j_{\si,\si'}$}}\rput(19,-3.4){\small{$\io_{\si,\si'}$}}
\rput(5.6,-2){\small{$\id$}}\rput(13.1,-2){\small{$p_{\si,\si'}$}}
\rput(23.4,-2){\small{$p_{\si}$}}
\rput(-6,-1.5){$\si'\!>\!\si\!:$}
\end{pspicture}

\noindent
The collection $\{\io_{\si,\si'}\}_{\si'\in\A_1(d,k)}$ of immersions 
is properly self-intersecting by the inductive assumption ($I14$)
of Subsection~\ref{map1blconstr_subs}.
Thus, by the same argument as in the paragraph following Lemma~\ref{cone1bl_lmm1},
so is the collection $\{j_{\si,\si'}\}_{\si'\in\A_1(d,k)}$.
Furthermore, 
\begin{equation}\label{idebundle_e}
\N_{j_{\si,\si'}}^{\ide}=\big\{\id\!\times\!p_{\si,\si'}\big\}^*
\N_{\io_{\si,\si'}}^{\ide}
\end{equation}
is an idealized normal bundle for $j_{\si,\si'}$.
These two observations also follow from Lemma~\ref{cone1bl_lmm1}
by induction using Lemmas~\ref{virimmer_lmm} and~\ref{virimmer_lmm2}.

\begin{lmm}
\label{coneblstr_lmm}
If $\si\!\in\!\A_1(d,k)$, $\V_{1,\si}^{\si-1}$ is a smooth subvariety of $\V_{1,k}^{\si-1}$
and 
$$p_{\si}\!: \V_{1,k}^{\si}\lra\ov\M_{1,k}^{\si}$$
is the idealized blowup of $\V_{1,k}^{\si-1}$ along $\V_{1,\si}^{\si-1}$.
\end{lmm}

\noindent
Recall from Subsection~\ref{map1blconstr_subs} that the immersion
$$\bar\io_{\si-1,\si}\!:
\big(\wt\cM_{1,(I_P(\si),J_P(\si))}\!\times\!\ov\M_{0,(\ale_B(\si),J_B(\si))}\big)\big/G_{\si}
\lra \ov\M_{1,\si}^{\si-1}\subset\ov\M_{1,k}^{\si-1}$$
induced by $\io_{\si-1,\si}$ is an embedding and 
$$\ti\pi_{\si}\!: \ov\M_{1,k}^{\si}\lra\ov\M_{1,k}^{\si-1}$$ 
is the idealized blowup along $\ov\M_{1,\si}^{\si-1}$.
Thus, the immersion 
$$\bar{j}_{\si-1,\si}\!:
\big(\wt\cM_{1,(I_P(\si),J_P(\si))}\!\times\!\V_{(\ale_B(\si),J_B(\si))}\big)\big/G_{\si}
\lra \V_{1,\si}^{\si-1}\subset \V_{1,k}^{\si-1}$$
induced by $j_{\si-1,\si}$ is also an embedding and 
$\V_{1,\si}^{\si-1}$ is a smooth subvariety of $\V_{1,k}^{\si-1}$.
Let
$$\ti\pi_{\si}\!: \V\lra\V_{1,k}^{\si-1}$$
be the idealized blowup along $\V_{1,k}^{\si-1}$.
Since
$$\N_{j_{\si-1,\si}}^{\ide}=\big\{\id\!\times\!p_{\si-1,\si}\big\}^*
\N_{\io_{\si-1,\si}}^{\ide}$$
and the linear map 
$$j_{\si-1,\si}^*TC\V_{1,k}^{\si-1}\big/\Im\, dj_{\si-1,\si}\lra
\io_{\si-1,\si}^*TC\ov\M_{1,k}^{\si-1}\big/\Im\, d\io_{\si-1,\si}$$
induced by $dp_{\si-1}$ is injective, $p_{\si-1}$ lifts to a map~$p$
over the blowdown maps~$\ti\pi_{\si}$:\\

\begin{pspicture}(-1.5,-2)(10,.3)
\psset{unit=.4cm}
\rput(10,0){$\V$}\rput(18,0){$\V_{1,k}^{\si-1}$}
\rput(10,-4){$\ov\M_{1,k}^{\si}$}\rput(18,-4){$\ov\M_{1,k}^{\si-1}$}
\psline{->}(10.8,0)(16.6,0)\psline{->}(11.5,-4)(16.5,-4)
\psline{->}(10,-.8)(10,-3)\psline{->}(18,-.8)(18,-3)
\rput(13.5,.6){\small{$\ti\pi_{\si}$}}\rput(14,-3.5){\small{$\ti\pi_{\si}$}}
\rput(10.5,-2){\small{$p$}}\rput(19.1,-2){\small{$p_{\si-1}$}}
\end{pspicture}

\noindent
Then $p$ and the top arrow\footnote{added ``the top arrow'' to make
clear which $\ti \pi_{\si}$ you mean.  --- R} $\ti\pi_{\si}$ factor through a morphism $f$ to $\ti\pi_{\si}^*\V_{1,k}^{\si-1}$:\\

\begin{pspicture}(-2.5,-2)(10,1.2)
\psset{unit=.4cm}
\rput(5,3){$\V$}
\psline{->}(5.5,2.9)(16.7,.4)\psline{->}(5,2.3)(9,-3.2)\psline{->}(5.2,2.6)(9,.7)
\rput(6.8,-1){\small{$p$}}\rput{-12}(12,2){\small{$\ti\pi_{\si}$}}
\rput{-35}(8,1.8){\small{$f$}}
\rput(10,0){$\ti\pi_{\si}^*\V_{1,k}^{\si-1}$}\rput(18,0){$\V_{1,k}^{\si-1}$}
\rput(10,-4){$\ov\M_{1,k}^{\si}$}\rput(18,-4){$\ov\M_{1,k}^{\si-1}$}
\psline{->}(11.8,0)(16.6,0)\psline{->}(11.5,-4)(16.5,-4)
\psline{->}(10,-.8)(10,-3)\psline{->}(18,-.8)(18,-3)
\rput(13.5,.5){\small{$\ti\pi_{\si}$}}\rput(14,-3.5){\small{$\ti\pi_{\si}$}}
\rput(10.7,-2){\small{$p_{\si}$}}\rput(19.1,-2){\small{$p_{\si-1}$}}
\end{pspicture}

\noindent
We show in the next paragraph that $f$ is an isomorphism.
Since $\ti\pi_{\si}^*\V_{1,k}^{\si-1}\!=\!\V_{1,k}^{\si}$,
this implies the second statement of Lemma~\ref{coneblstr_lmm}.\\

\noindent
By construction, the maps 
$$\ti\pi_{\si}\!: \ov\M_{1,k}^{\si}\lra\ov\M_{1,k}^{\si-1} \qquad\hbox{and}\qquad
\ti\pi_{\si}\!: \V\lra\V_{1,k}^{\si-1}$$
are isomorphisms on the complements of the idealized exceptional divisors
$$\ov\M_{1,\si}^{\si}\equiv\cE_{\ov\M_{1,\si}^{\si-1}}^{\ide}\subset\ov\M_{1,k}^{\si}
\qquad\hbox{and}\qquad
\V_{1,\si}^{\si}\equiv\cE_{\V_{1,\si}^{\si-1}}^{\ide}\subset\V.$$
Thus, $f\!: \V\!\lra\!\ti\pi_{\si}^*\V_{1,k}$ is an isomorphism over the complement of
$\ov\M_{1,\si}^{\si}$ in $\ov\M_{1,k}^{\si}$.
In particular, $f$ is linear on all fibers of~$p$.
Furthermore,
$$\ti\pi_{\si}^*\V_{1,k}^{\si-1}\big|_{\ov\M_{1,\si}^{\si}}
=\big\{(\ell,v)\!\in\!\P\N_{\ov\M_{1,\si}^{\si-1}}^{\ide}\!\times\!\V_{1,k}^{\si-1}\!:
\ti\pi_{\si}(\ell)\!=\!p_{\si-1}(v)\big\}.$$
On the other hand, since 
$$\N_{\V_{1,\si}^{\si-1}}^{\ide}=p_{\si-1}^*\N_{\ov\M_{1,\si}^{\si-1}}^{\ide}$$
by~\e_ref{idebundle_e}, we have 
$$\V_{1,\si}^{\si}=p_{\si-1}^*\P\N_{\ov\M_{1,\si}^{\si-1}}^{\ide}
=\big\{(v,\ell)\!\in\!\V_{1,k}^{\si-1}\!\times\!\P\N_{\ov\M_{1,\si}^{\si-1}}^{\ide}\!:
p_{\si-1}(v)\!=\!\ti\pi_{\si}(\ell)\big\}.$$
Thus, the restriction of $f$ to $\V_{1,\si}^{\si}$ must interchange $v$ and $\ell$,
i.e.~it is a vector bundle isomorphism over $\ov\M_{1,k}^{\si}$.
Finally, $\ov\M_{1,\si}^{\si-1}$ is a smooth subvariety of~$\V_{1,k}^{\si-1}$ and
$$T\big(\ov\M_{1,k}^{\si-1}\!\cap\!\V_{1,\si}^{\si-1}\big)=
T\ov\M_{1,\si}^{\si-1}
=TC\ov\M_{1,k}^{\si-1}\cap T\V_{1,\si}^{\si-1}\subset TC\V_{1,k}^{\si-1}.$$
Thus, similarly to (1) of Lemma~\ref{ag_lmm2a}, 
the proper transform of $\ov\M_{1,k}^{\si-1}$ in $\V$
is the blowup of $\ov\M_{1,k}^{\si-1}$ along 
$$\ov\M_{1,k}^{\si-1}\!\cap\!\V_{1,\si}^{\si-1}=\ov\M_{1,\si}^{\si-1},$$
i.e.~$\V$ contains $\ov\M_{1,k}^{\si}$ as the zero section.
The map~$f$ must be the identity on~$\ov\M_{1,k}^{\si}$.
Since $f$ is a linear isomorphism on all fibers of~$p$ by the above,
it then follows that $f$ is an isomorphism everywhere.\\

\noindent
{\it Remark:} If $\V_{1,k}^{\si-1}$ is a vector bundle over $\ov\M_{1,k}^{\si-1}$,
the second statement of Lemma~\ref{coneblstr_lmm} applies to standard blowups of
$\ov\M_{1,k}^{\si-1}$ and $\V_{1,k}^{\si-1}$ as well.
However, the second statement does not generally apply to standard blowups 
in the setting of Lemma~\ref{coneblstr_lmm}, as 
the analogue of the morphism $f$ may not be surjective.\\

\noindent
By the inductive assumption ($I1$) of Subsection~\ref{map1blconstr_subs}, 
the projection map~$\pi_{\si}$ is an isomorphism outside of the subvarieties
$\V_{1,\si'}^{\si}$ with $\si'\!\le\!\si$.
We denote by 
$$\V_{1,(0)}^{\si}\subset \V_{1,k}^{\si}$$ 
the proper transform of $\V_{1,(0)}^0$.
For each $\si'\!\in\!\A_1(d,k)$, let 
$$\W_{\si'}^{\si}=j_{\si,\si'}^{-1}\big(\V_{1,(0)}^{\si}\big)
=j_{\si,\si'}^{-1}\big(\V_{1,(0)}^{\si}\!\cap\!\V_{1,\si'}^{\si}\big).$$
By the inductive assumption ($I15$) of Subsection~\ref{map1blconstr_subs},
$\ov\M_{1,(0)}^{\si}$ is a proper subvariety of $\ov\M_{1,k}^{\si}$
with respect to the immersion~$\io_{\si,\si'}$.
Thus, by the same argument as in the paragraph following Lemma~\ref{cone1bl_lmm2},
the subvariety $\V_{1,(0)}^{\si}$ of $\V_{1,k}^{\si}$ is proper with respect 
to the immersion~$j_{\si,\si'}$.
Furthermore, if 
$$\N\W_{\si'}^{\si} \equiv \N\W_{j_{\si,\si'}|\V_{1,(0)}^{\si}}
\subset \N_{j_{\si',\si}}^{\ide}$$ 
denotes the normal cone for $j_{\si',\si}|_{\W_{\si'}^{\si}}$ in $\V_{1,(0)}^{\si}$,
then the homomorphism 
$$\N\W_{\si'}^{\si} \lra \{\id\!\times\!p_{\si,\si'}\}^*\N\bar\cZ_{\si'}^{\si}$$
induced by $dp_{\si}$ is injective.
These two observations also follow from Lemma~\ref{cone1bl_lmm2}
by induction using Lemma~\ref{virimmer_lmm3}.\\

\noindent
If $\si'\!\in\!\A_1(d,k)$, let 
$$\wt\cZ_{\si';B}^0=\wt\cD_0^{-1}(0)\cap
\P F_{(I_P(\si'),J_P(\si'))}\big|_{\M_{0,(\ale_B(\si),J_B(\si))}}.$$
By the inductive assumptions ($I7$) in Subsection~\ref{map1blconstr_subs} and
($I4$) in Subsection~\ref{map0blconstr_subs},
$$\ov\cZ_{\si'}^{\si}\equiv\io_{\si,\si'}^{\,-1}\big(\ov\M_{1,(0)}^{\si}\big)$$
is the closure of $\wt\cM_{1,(I_P(\si'),J_P(\si'))}\!\times\!\wt\cZ_{\si';B}^0$
in 
$$\wt\cM_{1,(I_P(\si'),J_P(\si'))}\times\wt\M_{0,(\ale_B(\si'),J_B(\si'))}^{\vr_{\si'}(\si)}$$
for all $\si\!\in\!\A_1(k,d)$ such that $\si'\!\le\!\si$.\\

\noindent
Suppose $\si\!\in\!\{0\}\!\cup\!\A_1(d,k)$ and $\si'\!\in\!\A_1(d,k)$. We claim that\\
${}\quad$ ($I1$) if $\si'\!\le\!\si$, then $\W_{\si'}^{\si}$ is the closure of
\begin{equation*}\begin{split}
\wt\cM_{1,(I_P(\si'),J_P(\si'))}\times \ker\wt\D_0|_{\wt\cZ_{\si';B}^0}
&\subset  \wt\cM_{1,(I_P(\si'),J_P(\si'))}\times
\pi_{\P F_{(\ale,J)}}^*\V_{(\ale_B(\si'),J_B(\si'))},\\
&\qquad \wt\cM_{1,(I_P(\si'),J_P(\si'))}\times
\pi_{\vr_{\si'}(\si)}^*\pi_{\P F_{(\ale,J)}}^*\V_{(\ale_B(\si'),J_B(\si'))}
\end{split}\end{equation*}
in $\wt\cM_{1,(I_P(\si'),J_P(\si'))}\times
\pi_{\vr_{\si'}(\si)}^*\pi_{\P F_{(\ale,J)}}^*\V_{(\ale_B(\si'),J_B(\si'))}$ and
$$\N\W_{\si'}^{\si}=\N_{j_{\si,\si'}}^{\ide}\big|_{\W_{\si'}^{\si}};$$
${}\quad$ ($I2$) if $\si'\!>\!\si$, then $\W_{\si'}^{\si}$ and $\N\W_{\si'}^{\si}$
are the closures~of
\begin{gather*}
\W_{\si'}^0|_{\cZ_{\si'}^0}\subset
\ov\cM_{1,(I_P(\si'),J_P(\si'))}\!\times\!\V_{(\ale_B(\si'),J_B(\si'))},
\ov\cM_{1,(I_P(\si'),J_P(\si'))}^{\rho_{\si'}(\si)}\!\times\!\V_{(\ale_B(\si'),J_B(\si'))}\\
\hbox{and}\qquad\
\N\W_{\si'}^0\big|_{\W_{\si'}^0|_{\cZ_{\si'}^0}}
\subset \N_{j_{0,\si'}}^{\ide}|_{\W_{\si'}^0|_{\cZ_{\si'}^0}}
\subset \N_{j_{0,\si'}}^{\ide},\N_{j_{\si,\si'}}^{\ide}
\end{gather*}
in $\ov\cM_{1,(I_P(\si'),J_P(\si'))}^{\rho_{\si'}(\si)}\!\times\!\V_{(\ale_B(\si'),J_B(\si'))}$
and in $\N_{j_{\si,\si'}}^{\ide}$, respectively.\\

\noindent
If $\si\!=\!0$, the assumption ($I1$) is trivially satisfied, 
while ($I2$) constitutes part of Lemma~\ref{cone1bl_lmm2}.
Suppose $\si\!\in\!\A_1(d,k)$ and the two assumptions hold with $\si$
replaced by~$\si\!-\!1$. 
By Lemma~\ref{coneblstr_lmm}, $\V_{1,k}^{\si}$ is the idealized blowup
of $\V_{1,k}^{\si-1}$ along~$\V_{1,\si}^{\si-1}$.
Thus, by the last statement of Lemma~\ref{virimmer_lmm3} both of 
the inductive assumptions continue to hold for~$\si'\!\neq\!\si$.\\

\noindent
On the other hand, let
\begin{equation*}\begin{split}
\cZ_{\si;B}&= \big\{b\!\in\!\M_{0,(\ale_B(\si),J_B(\si))}\!: 
\ker\cD_{(\ale_B(\si),J_B(\si))}\!\neq\!0\big\}, \\
\W_{\si;B}^0&=\big\{
\xi\!\in\!\V_{(\ale_B(\si),J_B(\si))}|_{\cZ_{\si;B}}\!:
\ker\,\{\D_{(\ale_B(\si),J_B(\si))}\xi\}
\big|_{\ker\cD_{(\ale_B(\si),J_B(\si))}|_{p_{\si-1,\si}(\xi)}}\!\neq\!\{0\}\big\},
\qquad\hbox{and}\\
\N\W_{\si;B}^0&=
\big\{(\xi,\ups)\!:\xi\!\in\!\W_{\si;B}^0, \,
\ups\!\in\!\ker\{\D_{0,\si}\xi\}\big|_{\ker\cD_{(\ale_B(\si),J_B(\si))}|_{p_{\si-1,\si}(\xi)}}\big\}
\subset p_{\si-1,\si}^{\,*}F_{(\ale_B(\si),J_B(\si))}.
\end{split}\end{equation*}
By the inductive assumption ($I12$) in Subsection~\ref{map1blconstr_subs},
\begin{gather*}
\ov\cZ_{\si}^{\si-1}\cap\big(
\wt\cM_{1,(I_P(\si),J_P(\si))}\!\times\!\M_{0,(\ale_B(\si),J_B(\si))}\big)
=\wt\cM_{1,(I_P(\si),J_P(\si))}\!\times\!\cZ_{\si;B} \qquad\hbox{and}\\
\N\bar\cZ_{\si}^0|_{\wt\cM_{1,(I_P(\si),J_P(\si))}\times\cZ_{\si;B}}
=\pi_P^*\L\otimes\pi_B^*\ker\cD_{0,(\ale_B(\si),J_B(\si))}.
\end{gather*}
By the inductive assumption~($I2$) above, Lemma~\ref{cone1bl_lmm2},
and the inductive assumption~($I11$) in Subsection~\ref{map1blconstr_subs},
$\W_{\si}^{\si-1}$ and $\N\W_{\si}^{\si-1}$ are the closures~of 
\begin{gather*}
\W_{\si}^0\big|_{\cZ_{\si}^0}=\big\{
\xi\!\in\!\pi_B^*\V_{(\ale_B(\si),J_B(\si))}|_{\cZ_{\si}^0}\!:
\ker\,\{\D_{0,\si}\xi\}
\big|_{\N\bar\cZ_{\si}^0|_{\{\id\times p_{0,\si}\}(\xi)}}\!\neq\!\{0\}
\big\}\\
\hbox{and}\qquad
\N\W_{\si}^0\big|_{\W_{\si}^0|_{\cZ_{\si}^0}}=
\big\{(\xi,\ups)\!:\xi\!\in\!\W_{\si}^0\big|_{\cZ_{\si}^0},~
\ups\!\in\!\ker\{\D_{0,\si}\xi\}
\big|_{\N\bar\cZ_{\si}^0|_{\{\id\times p_{0,\si}\}(\xi)}}\big\}
\subset\N_{j_{0,\si}}^{\ide}, \N_{j_{\si-1,\si}}^{\ide}
\end{gather*}
in  $\wt\cM_{1,(I_P(\si),J_P(\si))}\!\times\!\V_{(\ale_B(\si),J_B(\si))}$
and in 
$$\N_{j_{\si-1,\si}}^{\ide}= \pi_P^*\L\otimes 
\pi_B^*p_{\si-1,\si}^*F_{(\ale_B(\si),J_B(\si))}.$$
As before,
$$\pi_P,\pi_B\!: \wt\cM_{1,(I_P(\si),J_P(\si))}\!\times\!\V_{(\ale_B(\si),J_B(\si))}
\lra  \wt\cM_{1,(I_P(\si),J_P(\si))},\V_{(\ale_B(\si),J_B(\si))}$$
are the projections onto the principle and bubble components.
The bundle homomorphisms $s_i$ and $\ti{s}_i$ of Subsection~\ref{curve1bl_subs} 
agree on 
$$\cM_{1,(I_P(\si),J_P(\si))}\subset \ov\cM_{1,(I_P(\si),J_P(\si))},
\wt\cM_{1,(I_P(\si),J_P(\si))}.$$
The homomorphism $\ti{s}_i$ is an isomorphism from $\ti{L}_i$ to $\ti{E}^*$
over $\wt\cM_{1,(I_P(\si),J_P(\si))}$, and both line bundles are isomorphic to~$\L$.
It follows that $\W_{\si}^{\si-1}$ and $\N\W_{\si}^{\si-1}$ are the closures~of 
$$\wt\cM_{1,(I_P(\si),J_P(\si))}\times\W_{\si;B}^0 \qquad\hbox{and}\qquad
\pi_P^*\L\otimes\pi_B^*\N\W_{\si;B}^0$$
in  $\wt\cM_{1,(I_P(\si),J_P(\si))}\!\times\!\V_{(\ale_B(\si),J_B(\si))}$
and in 
$$\N_{j_{\si-1,\si}}^{\ide}=\pi_P^*\L\otimes
\pi_B^*p_{\si-1,\si}^{\,*}F_{(\ale_B(\si),J_B(\si))}.$$
Thus, by the first statement of Lemma~\ref{virimmer_lmm3}, 
$$\W_{\si}^{\si}\equiv j_{\si,\si}^{\,-1}\big(\V_{1,(0)}^{\si}\big)$$
is the closure of 
\begin{equation*}\begin{split}
\P\big(\pi_P^*\L\!\otimes\!\pi_B^*\N\W_{\si;B}^0\big)
&=\wt\cM_{1,(I_P(\si),J_P(\si))}\times \P\N\W_{\si;B}^0
=\wt\cM_{1,(I_P(\si),J_P(\si))}\times
\ker\wt\D_0|_{\wt\cZ_{\si;B}^0}\\
&\quad\subset\P \N_{j_{\si-1,\si}}^{\ide}
=\wt\cM_{1,(I_P(\si),J_P(\si))}\!\times\!
\pi_{\P F_{(\ale,J)}}^*\V_{(\ale_B(\si),J_B(\si))},
\end{split}\end{equation*}
i.e.~the first part of the inductive assumption ($I1$) for $\si'\!=\!\si$ is satisfied.
Furthermore, by the second part of~(1) of Lemma~\ref{virimmer_lmm3},
$$\N\W_{\si}^{\si}=\ga_{\V_{1,\si}^{\si-1}}\big|_{\W_{\si}^{\si}}
=\big\{\id\!\times\!p_{\si,\si}\big\}^*\ga_{\ov\M_{1,\si}^{\si-1}}\big|_{\W_{\si}^{\si}}
=\big\{\id\!\times\!p_{\si,\si}\big\}^*\N_{\io_{\si,\si}}^{\ide}\big|_{\W_{\si}^{\si}}
=\N_{j_{\si,\si}}^{\ide}\big|_{\W_{\si}^{\si}}.$$
We have thus verified the second  part of the inductive assumption ($I1$) 
for $\si'\!=\!\si$.\\

\noindent
Since the immersions $\bar{i}_{\si,\si'}$ with $\si'\!\le\!\si$ are embeddings by
the inductive assumption~($I8$) in Subsection~\ref{map1blconstr_subs}, so are 
the immersions
$$\bar{j}_{\si,\si'}\!: \big(\wt\cM_{1,(I_P(\si'),J_P(\si'))}\!\times\!
\pi_{\vr_{\si'}(\si)}^*\pi_{\P F_{(\ale,J)}}^*\V_{(\ale_B(\si'),J_B(\si'))}\big)
\big/G_{\si'}\lra \V_{1,\si'}^{\si}\subset\V_{1,k}^{\si}$$
induced by $j_{\si,\si'}$.
In particular, all of the morphisms 
\begin{equation*}\begin{split}
\bar{j}_{\si_{\max},\si'}\!: \big(\wt\cM_{1,(I_P(\si'),J_P(\si'))}\!\times\!
\pi_{0,(\ale_B(\si'),J_B(\si'))}^*\pi_{\P F_{(\ale,J)}}^*\V_{(\ale_B(\si'),J_B(\si'))}\big)
\big/G_{\si'} \qquad\qquad\qquad&\\
\lra \V_{1,\si'}^{\si_{\max}}\subset\V_{1,k}^{\si_{\max}}=
\ti\pi^*\ov\M_{1,k}(\cL,d)&
\end{split}\end{equation*}
are embeddings.
On the other hand, by the inductive assumption~($I1$),
$$\wt\W_{\si'}\equiv\W_{\si'}^{\si_{\max}}\equiv
j_{\si_{\max},\si'}^{\,-1}\big(\wt\M_{1,k}^0(\cL,d)\big)
\equiv j_{\si_{\max},\si'}^{\,-1}\big(\V_{1,(0)}^{\si_{\max}}\big)$$
is the closure of 
$$\wt\cM_{1,(I_P(\si'),J_P(\si'))}\times
\ker\wt\D_0|_{\wt\cZ_{\si';B}^0}
\subset \wt\cM_{1,(I_P(\si'),J_P(\si'))}\!\times\!
\pi_{0,(\ale_B(\si'),J_B(\si'))}^*\pi_{\P F_{(\ale,J)}}^*\V_{(\ale_B(\si'),J_B(\si'))}.$$
By Proposition~\ref{conehomomor_prp} and the inductive assumption~($I8$) 
in Subsection~\ref{map1blconstr_subs}, this closure~is
\begin{gather*}
\wt\cM_{1,(I_P(\si'),J_P(\si'))}\times\ker\wt\D_{(I_P(\si'),J_P(\si'))}|_{\wt\cZ_{\si';B}},
\qquad\hbox{where}\\
\wt\cZ_{\si';B}=\wt\cD_{(I_P(\si'),J_P(\si'))}^{-1}(0).
\end{gather*}
Since the bundle section $\wt\cD_{(I_P(\si'),J_P(\si'))}$ is transverse to the zero set,
$\wt\cZ_{\si';B}$ is a smooth subvariety of $\wt\M_{0,(I_P(\si'),J_P(\si'))}(\Pn,d)$
and
$$\wt\W_{\si'}\lra \wt\cM_{1,(I_P(\si'),J_P(\si'))}\times\wt\cZ_{\si';B}$$
is a smooth vector bundle by Proposition~\ref{conehomomor_prp}.
We conclude that 
$$\wt\M_{1,k}^0(\cL,d)\cap \V_{1,\si'}^{\si_{\max}}$$
is a smooth subvariety of $\wt\M_{1,k}^0(\cL,d)$ for all $\si'\!\in\!\A_{1,k}(k,d)$.
Its normal cone is a line bundle by the inductive assumption~($I1$).\\

\vspace{.2in}

\noindent
{\it Department of Mathematics, Stanford University, Stanford, CA 94305-2125}\\
vakil@math.stanford.edu\\

\noindent
{\it Department of Mathematics, SUNY, Stony Brook, NY 11794-3651}\\
azinger@math.sunysb.edu\\

\end{document}